\documentclass[twoside]{article}
\usepackage{etex}
\usepackage{lipsum} 
\usepackage{mathrsfs}

\usepackage[sc]{mathpazo}
\usepackage[T1]{fontenc} 
\usepackage[hmarginratio=1:1,top=32mm,columnsep=20pt]{geometry} 
\usepackage{multicol} 
\usepackage[hang, small,labelfont=bf,up,textfont=it,up]{caption} 
\usepackage{booktabs} 
\usepackage{paralist} 

\usepackage{abstract} 

\usepackage{titlesec} 
\titleformat{\section}[block]{\large\scshape}{\thesection.}{1em}{} 
\titleformat{\subsection}[block]{\large}{\thesubsection.}{1em}{} 

\usepackage{fancyhdr} 
\usepackage{amsthm}
\usepackage{amsmath,amstext,graphicx,epstopdf}
\usepackage{multirow,tikz,tipa,caption}
\usepackage{lineno,lipsum}
\usepackage{algorithmic}
\usepackage[ruled,vlined]{algorithm2e}
\usepackage{bm,adjustbox,mathtools,bigints}
\usepackage{subcaption}

\makeatletter
\newcounter{megaalgorithm}
\newenvironment{megaalgorithm}[1][htb]
  {
   \let\c@algocf\c@megaalgorithm
   \begin{algorithm}[#1]%
  }{\end{algorithm}}
\makeatother

\makeatletter
\newcounter{algorithmone}

\makeatother

\makeatletter
\newcounter{algorithmtwo}

\makeatother

\newtheorem{dummy}{***}
\newtheorem{thm}[dummy]{Theorem}

\newtheorem{rmk}[dummy]{Remark}
\newtheorem{cor}[dummy]{Corollary}

\usepackage{xparse}

\NewDocumentCommand{\ceil}{s O{} m}{%
  \IfBooleanTF{#1} 
    {\left\lceil#3\right\rceil} 
    {#2\lceil#3#2\rceil} 
}

\usepackage{parskip}
\usepackage{lipsum}
\makeatletter
\def\thm@space@setup{%
  \thm@preskip=12pt plus 0pt minus 8pt
}
\makeatother

\pgfmathsetmacro{\minsize}{0.2cm}

\newlength\myindent
\usetikzlibrary{intersections}

\usepackage{amsfonts}
\usepackage{amsmath,amssymb,latexsym}

\title{\vspace{-15mm}\fontsize{16pt}{10pt}\selectfont\textbf{A Convex Geodesic Selective Model for Image Segmentation}}

\author{
\large
\textsc{Michael Roberts$^\dagger$, Ke Chen$^\dagger$\thanks{Email k.chen@liverpool.ac.uk,\ Web: www.liv.ac.uk/cmit \ (corresponding author). Work supported by UK EPSRC   grant EP/N014499/1.}\,  and
Klaus L. Irion$^\ddagger$
}
\\[2mm]
\normalsize $^\dagger$Centre for Mathematical Imaging Techniques and\\
\normalsize   Department of Mathematical Sciences,\\
 \normalsize  University of Liverpool,
\normalsize   United Kingdom\\
\normalsize $^\ddagger$Department of Radiology, Manchester University NHS Foundation Trust,
\normalsize  Manchester, UK \\
\vspace{-5mm}
}
\date{}
\newcommand*\patchAmsMathEnvironmentForLineno[1]{%
  \expandafter\let\csname old#1\expandafter\endcsname\csname #1\endcsname
  \expandafter\let\csname oldend#1\expandafter\endcsname\csname end#1\endcsname
  \renewenvironment{#1}%
     {\linenomath\csname old#1\endcsname}%
     {\csname oldend#1\endcsname\endlinenomath}}%
\newcommand*\patchBothAmsMathEnvironmentsForLineno[1]{%
  \patchAmsMathEnvironmentForLineno{#1}%
  \patchAmsMathEnvironmentForLineno{#1*}}%
\AtBeginDocument{%
\patchBothAmsMathEnvironmentsForLineno{equation}%
\patchBothAmsMathEnvironmentsForLineno{align}%
\patchBothAmsMathEnvironmentsForLineno{flalign}%
\patchBothAmsMathEnvironmentsForLineno{alignat}%
\patchBothAmsMathEnvironmentsForLineno{gather}%
\patchBothAmsMathEnvironmentsForLineno{multline}%
}

\begin{document}

\maketitle

\thispagestyle{fancy} 


\begin{abstract}
\vspace{0.1in}
\noindent
\small
Selective segmentation is an important application of image processing. In contrast to global segmentation in which all objects are segmented, selective segmentation is used to isolate specific objects in an image and is of particular interest in medical imaging -- permitting segmentation and review of a single organ. An important consideration is to minimise the amount of user input to obtain the segmentation; this differs from interactive segmentation in which more user input is allowed than selective segmentation.
To achieve selection, we propose a selective segmentation model which uses the edge-weighted geodesic distance from a marker set as a penalty term.
It is demonstrated that this edge-weighted geodesic penalty term improves on previous selective penalty terms. A convex formulation of the model is also presented, allowing arbitrary initialisation. It is shown that the proposed model is less parameter dependent and requires less user input than previous models.
Further modifications are made to the edge-weighted geodesic distance term to ensure segmentation robustness to noise and blur.
We can show that the overall Euler-Lagrange equation admits a unique viscosity solution.
Numerical results show that the result is robust to user input and permits selective segmentations that are not possible with other models.

\end{abstract}



\section{Introduction}

Segmentation of an image into its individual objects is one incredibly important
application of image processing techniques.
 Segmentation can take two forms; firstly global segmentation for isolation of all foreground objects
 in an image from the background and secondly, selective segmentation for
 isolation of a subset of the objects in an image from the background.
 A comprehensive review of selective segmentation can be found in \cite{batra2011interactive,he2013interactive} and  in \cite{zhao2013overview}
 for  medical image segmentation
 where
selection refers to extraction of single organs.

 Approaches to image segmentation broadly fall into two classes; region-based and edge-based. Some region-based approaches are region growing \cite{RC1}, watershed algorithms \cite{BC31}, Mumford-Shah \cite{SC22} and Chan-Vese \cite{SC10}. The final two of these are partial differential equations (PDEs)-based variational approaches to the problem of segmentation. There are also models which mix the two classes to use the benefits of the region-based and edge-based approaches and will incorporate features of each. Edge-based methods aim to encourage an evolving contour towards the edges in an image and normally require an edge detector function \cite{SC5}.
 The first edge-based variational approach was devised by Kass et al.   \cite{SC17}
with the famous snakes model, this was further developed by Casselles et al.  \cite{SC5}
who introduced the Geodesic Active Contour (GAC) model.
Region-based global segmentation models include the well known works of
 Mumford-Shah \cite{SC22} and Chan-Vese \cite{SC10}.
 Importantly they are non-convex and hence a minimiser of these models may only be a local, not the global, minimum. Further work by Chan et al. \cite{SC8} gave rise to a method to find
the global minimiser for the Chan-Vese model under certain conditions.

This paper is mainly concerned with selective segmentation of objects in an image, given a set of points near the object or objects to be segmented.
It builds in such user input to a model
using a set $
\mathcal{M}=\{ (x_{i},y_{i})\in\Omega, 1\le i\le k\}
$
where $\Omega\subset\mathbb{R}^{2}$ is the image domain \cite{ SC1, SC2,RC16}.
 Nguyen et al. \cite{SC25}  considered marker sets $\mathcal{M}$ and $\mathcal{A}$ which consist of points inside and outside, respectively, the object or objects to be segmented.
Gout et al. \cite{RC16} combined the GAC approach with the geometrical constraint that
 the contour passes through the points of $\mathcal{M}$. This was enforced with a
 distance function which is zero at $\mathcal{M}$ and non-zero elsewhere.
Badshah and Chen \cite{SC1} then   combined the Gout et al. model with
  \cite{SC10}  to incorporate a constraint on the intensity in the
selected region, thereby encouraging the contour to segment homogeneouss regions.
Rada and Chen \cite{SC27} introduced a selective segmentation method based on
 two-level sets which was shown to be more robust than the Badshah-Chen model.
We also refer to \cite{SC2,MR1} for selective segmentation models which include
different fitting constraints, using coefficient of variation and
the centroid of $\mathcal{M}$ respectively.
None of these models have a restriction on the size of the object or objects to be detected and depending on the initialisation these methods have the potential to detect more or
 fewer objects than the user desired.
To address this and to improve on \cite{SC27}, Rada and Chen \cite{SC28} introduced a
model combining the Badshah-Chen \cite{SC1} model with a constraint on the area of
the objects to be segmented.
The reference area used to constrain the area within the
contour is that of the polygon formed by the markers in $\mathcal{M}$.
Spencer and Chen \cite{SC} introduced a model with the distance fitting
penalty as a standalone term in the energy functional, unbounding it from the edge
detector term of the Gout et al. model.

All of the above selective segmentation models discussed are non-convex and hence
the final result depends on the initialisation. Spencer and Chen \cite{SC}, in the same paper,
reformulated the model they introduced to a convex form using convex relaxation and an exact penalty term as in \cite{SC8}. Their model uses Euclidean distance from the marker set $\mathcal{M}$ as a distance penalty term, however we propose replacing this with the edge-weighted geodesic distance from $\mathcal{M}$ (we call this simply the geodesic distance). This distance increases at edges in the image and is more intuitive for selective segmentation.
The proposed model is given as a convex relaxed model with exact penalty term and we give a general existence and uniqueness proof for the viscosity solution to the PDE given by its Euler-Lagrange equation, which is also applicable to a whole class of PDEs arising in image segmentation. We note that the use of geodesic distance for segmentation has been considered before \cite{bai2009geodesic,protiere2007interactive}, however the models only use geodesic distance as the fitting term within the regulariser, so are liable to make segmentation errors for poor initialisation or complex images. Here we take a different approach, by including geodesic distance as a standalone fitting term, separate from the regulariser, and using intensity fitting terms to ensure robustness.


In this paper we only consider 2D images, however for completion we remark that 3D segmentation models do exist \cite{SC13,SC31} and it is simple to extend the proposed model to 3D.
The contributions of this paper can be summarised as follows:
\begin{itemize}
\item{We incorporate the geodesic distance as a distance penalty term within the variational framework.}
\item{We propose a convex selective segmentation model using this penalty term and demonstrate how it can achieve results which cannot be achieved by other models.}
\item{We improve the geodesic penalty term, focussing on improving robustness to noise and improving segmentation when object edges are blurred.}
\item{We give an existence and uniqueness proof for the viscosity solution for the PDEs associated with a whole class of segmentation models (both global and selective).}
\end{itemize}
We find that the proposed model gives accurate segmentation results for a wide range of parameters and, in particular, when segmenting the same objects from the same modality images, i.e. segmenting lungs from CT scans, the parameters are very similar from one image to the next to obtain accurate results.
Therefore, this model may be used to assist the preparation of large training sets for deep learning studies \cite{Learn17a,Learn17b,DBLP:journals/corr/abs-1710-04043} that concern segmentation of particular objects from images. 

The paper is structured as follows; in \S 2 we review some
global and selective segmentation models.
In \S 3 we discuss the geodesic distance penalty term, propose a new convex model and also address weaknesses in the na\"ive implementation of the geodesic distance term. In \S 4 we discuss the non-standard AOS scheme, introduced in \cite{SC}, which we use to solve the model. In \S 5 we give an existence and uniqueness proof for a general class of PDEs arising in image segmentation, thereby showing that for a given initialisation the solution to our model is unique. In \S 6 we compare the results of the proposed model to other selective segmentation models, show that the proposed model is less parameter dependent than other models and is more robust to user input. Finally, in \S 7 we provide some concluding remarks.


\section{Review of Variational Segmentation Models}

Although we focus on selective segmentation, it is illuminating to introduce some global segmentation models first. Throughout this paper we denote the original image by $z(x,y)$ with image domain $\Omega\subset\mathbb{R}^{2}$.

\subsection{Global Segmentation}

The model of Mumford and Shah \cite{SC22} is one of the most famous and
important variational models in image segmentation.
We will review its two-dimensional piecewise constant variant,
commonly  known as the Chan-Vese model \cite{SC10},
which takes the form
\begin{equation}\label{eqn:cv}
F_{CV}(\Gamma,c_{1},c_{2})
=\mu\cdot length (\Gamma)+\lambda_{1}\int_{\Omega_{1}}
|z(x,y)-c_{1}|^{2}\,\mathrm{d}\Omega +\lambda_{2}\int_{\Omega_{2}}|z(x,y)-c_{2}|^{2}\,\mathrm{d}\Omega
\end{equation}
where
the foreground $\Omega_{1}$ is the subdomain to be segmented,
  the background $\Omega_{2}=\Omega\backslash\Omega_{1}$
  and $\mu,\lambda_{1},\lambda_{2}$ are fixed non-negative parameters.
The values $c_{1}$ and $c_{2}$ are the average intensities of $z(x,y)$
inside   $\Omega_{1}$ and $\Omega_{2}$ respectively.
We use a level set function
\[
\phi(x,y)=
\begin{cases}
    >0,& (x,y)\in\Omega_{1},\\
    0,& (x,y)\in\Gamma,\\
    <0, & otherwise,\\
\end{cases}
\]
to track  $\Gamma = \{ (x,y)\in\Omega\, |\, \phi(x,y)=0\}$ (an idea popularised by Osher and Sethian \cite{OS1})
 and reformulate  (\ref{eqn:cv}) as
 \begin{equation}\label{eqn:cvls}
\begin{split}
F_{CV}(\phi,c_{1},c_{2})=&\mu\int_{\Omega}|\nabla H_{\varepsilon}(\phi)|\,\mathrm{d}\Omega+\lambda_{1}\int_{\Omega}(z(x,y)-c_{1})^{2}H_{\varepsilon}(\phi)\,\mathrm{d}\Omega\\
&
\ \hspace*{2.7cm}+\lambda_{2}\int_{\Omega}(z(x,y)-c_{2})^{2}(1-H_{\varepsilon}(\phi))\,\mathrm{d}\Omega,\\
\end{split}
\end{equation}
with $H_{\varepsilon}(\phi)$ a smoothed Heaviside function such as $H_{\varepsilon}(\phi) = \frac{1}{2} + \frac{1}{\pi}\arctan(\frac{\phi}{\varepsilon})$ for some $\varepsilon$, we set $\varepsilon = 1$ throughout. We solve this in two stages, first
with $\phi$ fixed we minimise $F_{CV}$ with respect to $c_{1}$ and $c_{2}$, obtaining
\begin{equation}\label{eqn:c1c2}
c_{1}=\frac{\int_{\Omega}H_{\varepsilon}(\phi)\cdot z(x,y)\,\mathrm{d}\Omega}{\int_{\Omega}H_{\varepsilon}(\phi) \,\mathrm{d}\Omega},
\hspace{0.5in}
c_{2}=\frac{\int_{\Omega}(1-H_{\varepsilon}(\phi))\cdot z(x,y)\,\mathrm{d}\Omega}{\int_{\Omega}(1-H_{\varepsilon}(\phi)) \,\mathrm{d}\Omega},
\end{equation}
and secondly, with $c_{1}$ and $c_{2}$ fixed we minimise (\ref{eqn:cvls}) with respect to $\phi$. This requires the calculation of the associated Euler-Lagrange equations.
 A drawback of the Chan-Vese energy functional (\ref{eqn:cvls}) is that it is non-convex.
Therefore a minimiser may only be a local minimum and not the global minimum and the final segmentation result is dependent on the initialisation. Chan et al. \cite{SC8} reformulated
(\ref{eqn:cvls}) using an exact penalty term to obtain an equivalent convex model -- we use this same technique in \S\ref{sec:selsegintro} for the Geodesic Model.

\subsection{Selective Segmentation}

\label{sec:selsegintro}

Selective segmentation models make use of user input, i.e. a marker set $\mathcal{M}$ of points near the object or objects to be segmented. Let
\(
\mathcal{M}=\{ (x_{i},y_{i})\in\Omega, 1\le i\le k\}
\)
be such a marker set. The aim of selective segmentation is to design an energy functional where the segmentation contour $\Gamma$ is close to the points of $\mathcal{M}$.

{\bf Early work.} An early model by Caselles et al. \cite{SC5}, commonly known as the Geodesic Active Contour (GAC) model, uses an edge detector function to ensure the contour follows edges, the functional to minimise is given by
\[
\bigintsss_{\Gamma}g(|\nabla z(x,y)|)d\Gamma.
\]
The term $g(|\nabla z(x,y)|)$ is an edge detector, one example is $g(s) = 1/(1+\beta s^{2})$ with $\beta$ a tuning parameter. It is common to smooth the image with a Gaussian filter $G_{\sigma}$ where $\sigma$ is the kernel size, i.e. use $g(|\nabla\left( G_{\sigma}*z(x,y)\right)|)$ as the edge detector. This mitigates the effect of noise in the image, giving a more accurate edge detector. Gout et al. \cite{SC13} built upon the GAC model by incorporating a distance term $\mathcal{D}(x,y)$ into this integral, i.e. the integrand is $\mathcal{D}(x,y)g(|\nabla z|)$. The distance term is a penalty on the distance from $\mathcal{M}$, this model encourages the contour to be near to the set $\mathcal{M}$ whilst also lying on edges. However this model struggles when boundaries between objects and their background
are fuzzy or blurred.
To address this, Badshah and Chen \cite{SC1} introduced a new model which adds the intensity fitting terms from the Chan-Vese model (\ref{eqn:cv}) to the Gout et al. model.
However, their model has poor robustness \cite{SC27}.
To improve on this, Rada and Chen \cite{SC28} introduced a model which adds
an area fitting term into the Badshah-Chen model and is far more robust.

{\bf The Rada-Chen model} \cite{SC28}.
We first briefly introduce this model, defined by
\begin{equation}\label{eqn:rcfunc}
\begin{split}
F_{RC}(\phi,c_{1},c_{2})=&\mu\int_{\Omega}\mathcal{D}(x,y)g(|\nabla z(x,y)|)|\nabla
H_{\varepsilon}(\phi)|\,\mathrm{d}\Omega \\
& + \lambda_{1}\int_{\Omega}(z(x,y)-c_{1})^{2}H_{\varepsilon}(\phi)\,\mathrm{d}\Omega +
\lambda_{2}\int_{\Omega}(z(x,y)-c_{2})^{2}(1-H_{\varepsilon}(\phi))\,\mathrm{d}\Omega\\
& + \gamma\bigg[\left(\int_{\Omega}H_{\varepsilon}(\phi)\,\mathrm{d}\Omega-A_{1}\right)^{2}+
\left(\int_{\Omega}(1-H_{\varepsilon}(\phi))\,\mathrm{d}\Omega-A_{2}\right)^{2}\bigg],
\end{split}
\end{equation}
where $\mu,\lambda_{1},\lambda_{2},\gamma$ are fixed non-negative parameters. There is freedom in choosing the distance term $\mathcal{D}(x,y)$, see \cite{SC28} for some examples. $A_{1}$ is the area of the polygon formed from the points of $\mathcal{M}$ and $A_{2}=|\Omega|-A_{1}$.
The final term of this functional puts a penalty on the area inside a contour being very different to $A_{1}$. One drawback of the Rada-Chen model is that the selective fitting term uses no location information from the marker set $\mathcal{M}$. Therefore the result can be a contour which is separated over the domain into small parts, whose sum area totals the area fitting term.

{\bf Nguyen et al. \cite{SC25}.} This model is based on the GAC model and uses likelihood functions as fitting terms, it has the energy functional
\[
\begin{split}
F_{NG}(\phi)=&\mu\int_{\Omega}g(|\nabla z(x,y)|)
|\nabla H_{\varepsilon}(\phi)|\,\mathrm{d}\Omega \\
&+ \lambda \int_{\Omega} \alpha\left( P_{B}(x,y) - P_{F}(x,y)\right) + \left( 1-\alpha\right) \left(1-2P(x,y)\right) \phi\,\mathrm{d}\Omega
\end{split}
\]
where $P_{B}(x,y)$ and $P_{F}(x,y)$ are the normalised log-likelihoods that the pixel $(x,y)$ is in the foreground and background respectively. $P(x,y)$ is the probability that pixel $(x,y)$ belongs to the foreground, $\alpha\in[0,1]$ and minimisation is constrained, requiring $\phi\in [0,1]$, so $F_{NG}(\phi)$ is convex. This model is good for many examples, see \cite{SC25}, however fails when the boundary of the object to segment is non-smooth or has fine structures. Also, the final result is sometimes sensitive to the marker sets used.

{\bf The Spencer-Chen model \cite{SC}.} The authors introduced the following model
\begin{equation}\label{eqn:scfunc}
\begin{split}
F_{SC}(\phi,c_{1},c_{2})=&\mu\int_{\Omega}g(|\nabla z(x,y)|)
|\nabla H_{\varepsilon}(\phi)|\,\mathrm{d}\Omega
  + \lambda_{1}\int_{\Omega}(z(x,y)-c_{1})^{2}H_{\varepsilon}(\phi)\,\mathrm{d}\Omega\\
  & + \lambda_{2}\int_{\Omega}(z(x,y)-c_{2})^{2}(1-H_{\varepsilon}(\phi))\,\mathrm{d}\Omega
  + \theta\int_{\Omega}\mathcal{D}_{E}(x,y)H_{\varepsilon}(\phi)\,\mathrm{d}\Omega,
\end{split}
\end{equation}
where $\mu,\lambda_{1},\lambda_{2},\theta$ are fixed non-negative parameters.
Note that the regulariser of this model differs from the Rada-Chen model (\ref{eqn:rcfunc}) as
the distance function $\mathcal{D}(x,y)$ has been separated from the edge detector term and
is now a standalone penalty term $\mathcal{D}_{E}(x,y)$. The authors use normalised Euclidean distance $\mathcal{D}_{E}(x,y)$ from the marker set $\mathcal{M}$ as their distance penalty term. We will discuss this later in \S 3 as it is one of the key improvements we make to the Spencer-Chen model, replacing the Euclidean distance term with a geodesic distance term.

{\bf Convex Spencer-Chen model  \cite{SC}.} Spencer and Chen use the ideas of \cite{SC8} to reformulate (\ref{eqn:scfunc}) into a convex minimisation problem.
It can be shown that the Euler-Lagrange equations for $F_{SC}(\phi,c_{1},c_{2})$ have the same stationary solutions as for
\begin{equation}\label{eqn:sc1func}
\begin{split}
F_{SC1}(u,c_{1},c_{2})=&\mu\int_{\Omega}g(|\nabla z(x,y)|)
|\nabla u|\,\mathrm{d}\Omega
  + \int_{\Omega}\left[   \lambda_{1}(z(x,y)-c_{1})^{2}  - \lambda_{2}(z(x,y)-c_{2})^{2}  \right]  u\,\mathrm{d}\Omega \\
  &   + \theta\int_{\Omega}\mathcal{D}_{E}(x,y)u \,\mathrm{d}\Omega,
\end{split}
\end{equation}
with the minimisation constrained to $u\in[0,1]$. This is a constrained convex minimisation which can be reformulated to an unconstrained minimisation using an exact penalty term $\nu(u) := \max\{0, 2|u -\frac{1}{
2}|-1\}$ in the functional, which encourages the minimiser to be in the range $[0,1]$. In \cite{SC} the authors use a smooth approximation $\nu_{\varepsilon}(u)$ to $\nu(u)$ given by
\begin{equation}\label{eqn:nueps}\nu_{\varepsilon}(u) = H_{\varepsilon}\left(\sqrt{(2u-1)^{2}+\varepsilon}-1\right)\left[ \sqrt{(2u-1)^{2}+\varepsilon}-1 \right],
\end{equation}
and perform the unconstrained minimisation of
\begin{equation}\label{eqn:scfunc2}
\begin{split}
F_{SC2}(u,c_{1},c_{2})=&\mu\int_{\Omega}g(|\nabla z(x,y)|)
|\nabla u|\,\mathrm{d}\Omega
  + \int_{\Omega}\left[   \lambda_{1}(z(x,y)-c_{1})^{2}  - \lambda_{2}(z(x,y)-c_{2})^{2}  \right]  u\,\mathrm{d}\Omega \\
  &   + \theta\int_{\Omega}\mathcal{D}_{E}(x,y)u \,\mathrm{d}\Omega + \alpha\int_{\Omega}\nu_{\varepsilon}(u)\,\mathrm{d}\Omega.
\end{split}
\end{equation}
When $
\alpha > \frac{1}{2} \left|\left| \left[   \lambda_{1}(z(x,y)-c_{1})^{2}  - \lambda_{2}(z(x,y)-c_{2})^{2}  \right] + \theta\mathcal{D}_{E}(x,y) \right| \right|_{L^{\infty}}$,
the above functional has the same set of stationary solutions as $F_{SC1}(u,c_{1},c_{2})$. 
It permits us to choose arbitrary $u$ initialisation to obtain the desired selective segmentation result due to its complexity.

{\bf Convex Liu et al. model \cite{Liu2017}.} Recently, a convex model was introduced by Liu et al. which applies a weighting to the data fitting terms, the functional to minimise is given  by
\begin{equation}\label{eqn:liufunc}
\begin{split}
F_{LIU}(u)=&\mu\int_{\Omega}|\nabla u|\,\mathrm{d}\Omega + \mu_{2}\int_{\Omega} |\nabla u|^{2}\,\mathrm{d}\Omega
  + \lambda\int_{\Omega} \omega^{2}(x,y)\left|z - u\right|^{2}\,\mathrm{d}\Omega,
\end{split}
\end{equation}
where $\mu,\mu_{2},\lambda$ are non-negative parameters and $\omega(x,y) = 1-\mathcal{D}(x,y)g(|\nabla z|)$ where $\mathcal{D}(x,y)$ is a distance function from marker set $\mathcal{M}$ (see \cite{Liu2017} for examples).



\section{Proposed Convex Geodesic Selective Model}

\begin{figure}[htb!]
\makebox[\textwidth][c]{
\includegraphics[width=1.3\textwidth]{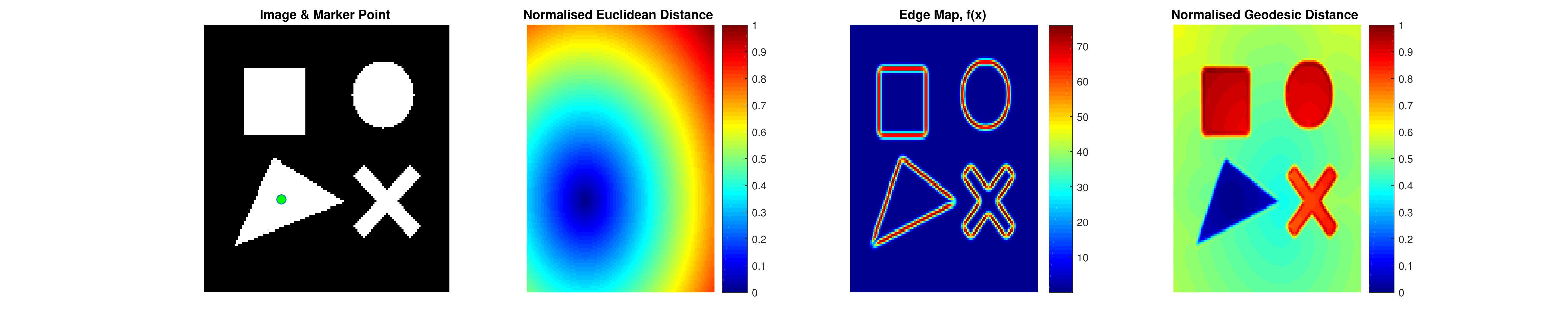}
}\vskip -5mm
\centerline{
(i) \hspace*{3.3cm}
(ii)  \hspace*{3.3cm}
 (iii)  \hspace*{3.3cm}
 (iv) }
\caption{Comparison of distance measures. (i) Simple binary image with marker point; (ii) normalised Euclidean distance from marker point; (iii) edge map function $f(x)$ for the image; (iv) normalised geodesic distance from marker point. \label{fig:eucgeodist}}
\end{figure}

We propose an improved selective model, based on the Spencer-Chen model, which uses geodesic distance from the marker set $\mathcal{M}$ as the distance term, rather than the Euclidean distance. Increasing the distance when edges in the image are encountered gives a more accurate reflection of the true similarity of pixels in an image from the marker set. We propose minimising the convex functional
\begin{equation}\label{eqn:geofunc2}
\begin{split}
F_{CG}(u,c_{1},c_{2})=&\mu\int_{\Omega}g(|\nabla z(x,y)|)
|\nabla u|\,\mathrm{d}\Omega
  + \int_{\Omega}\left[   \lambda_{1}(z(x,y)-c_{1})^{2}  - \lambda_{2}(z(x,y)-c_{2})^{2}  \right]  u\,\mathrm{d}\Omega \\
  &   + \theta\int_{\Omega}\mathcal{D}_{M}(x,y)u \,\mathrm{d}\Omega + \alpha\int_{\Omega}\nu_{\varepsilon}(u)\,\mathrm{d}\Omega,
\end{split}
\end{equation}
where $\mathcal{D}_{M}(x,y)$ is the edge-weighted geodesic distance from the marker set.
In Figure~\ref{fig:eucgeodist}, we compare  the normalised geodesic distance and the Euclidean distance
from the same   marker point (i.e. set $\mathcal{M}$ has
 one point in it); clearly the former gives
 a more intuitively correct distance penalty than the latter. We will refer to this proposed model as the Geodesic Model.

\subsection{Computing the Geodesic Distance Term $\mathcal{D}_{M}(x,y)$}

The geodesic distance from the marker set $\mathcal{M}$ is given by $\mathcal{D}_{M}(x,y) = 0$ for $(x,y)\in\mathcal{M}$ and $\mathcal{D}_{M}(x,y) = \frac{\mathcal{D}_{M}^{0}(x,y)}{||\mathcal{D}_{M}^{0}(x,y)||_{L^{\infty}}}$ for $(x,y)\not\in\mathcal{M}$, where $\mathcal{D}_{M}^{0}(x,y)$ is the solution of the following PDE
\begin{equation}\label{eqn:DG0}
|\nabla\mathcal{D}_{M}^{0}(x,y)| = f(x,y), \qquad \mathcal{D}^{0}_{M}(x_{0},y_{0}) = 0, \,(x_{0},y_{0})\in\mathcal{M}.
\end{equation}
where $f(x,y)$ is defined later on with respect to the image contents.

If $f(x,y)\equiv 1$ (i.e.  $|\nabla\mathcal{D}_{M}^{0}(x,y)|=1$) then the distance penalty $D_{M}(x,y)$ is simply the normalised Euclidean distance $\mathcal{D}_{E}(x,y)$ as used in the Spencer-Chen model (\ref{eqn:scfunc}). We have free rein to design $f(x,y)$ as we wish. Looking at the PDE in (\ref{eqn:DG0}), we see that when $f(x,y)$ small this results in a small gradient in our distance function and it is almost flat. When $f(x,y)$ is large, we have a large gradient in our distance map. In the case of selective image segmentation, we want small gradients in homogeneous areas of the image and large gradients at edges. If we set
\begin{equation}\label{eqn:edgemap}
f(x,y) = \varepsilon_{\mathcal{D}} + \beta_{G}|\nabla z(x,y)|^{2}
\end{equation}
this gives us the desired property that in areas where
$|\nabla z(x,y)|\approx 0$, the distance function increases by some small $\varepsilon_{\mathcal{D}}$; here image $z(x,y)$  is scaled to $[0,1]$. At edges, $|\nabla z(x,y)|$ is large and the geodesic distance increases here. We set  value of $\beta_{G} =1000$ and   $\varepsilon_{\mathcal{D}} = 10^{-3}$ throughout.
In Figure~\ref{fig:eucgeodist},
we see that the geodesic distance plot   gives a low distance penalty on the triangle, which the marker indicates we would like segmented. There is a reasonable penalty on the background, and all other objects in the image have a very high distance penalty (as the geodesic to these points must cross two edges). This contrasts with the Euclidean distance, which gives a low distance penalty to some background pixels and maximum penalty to the pixels furthest away.

\subsection{Comparing Euclidean and Geodesic Distance Terms}
\label{sec:geoadv}

We briefly give some advantages of using the geodesic distance as a penalty term rather than Euclidean distance and a remark on the computational complexity for both distances.
\begin{enumerate}
\item{{\bf Parameter Robustness.} The Geodesic Model is more robust to the choice of the fitting parameter $\theta$, as the penalty on the inside of the shape we want segmented is consistently small. It is only outside the shape where the penalty is large. Whereas with the Euclidean distance term we always have a penalty inside the shape we actually want to segment. This is due to the nature of the Euclidean distance which does not discriminate on intensity -- this penalty can also be quite high if our marker set is small and doesn't cover the whole object.}
\item{{\bf Robust to Marker Set Selection.} The geodesic distance term is far more robust to point selection, for example we can choose just one point inside the object we want to segment and this will give a nearly identical geodesic distance compared to choosing many more points. This is not true of the Euclidean distance term which is very sensitive to point selection and requires markers to be spread in all areas of the object you want to segment (especially at extrema of the object).}
\end{enumerate}

\begin{rmk}[Computational Complexity.] The main concern of using the geodesic penalty term, which we obtain by solving PDE (\ref{eqn:DG0}), would be that it takes a significant amount of time to compute compared to the Euclidean distance. However, using the fast marching algorithm of Sethian \cite{sethian1996fast}, the complexity of computing $\mathcal{D}_{M}(x,y)$ is $\mathcal{O}(N\log(N))$ for an image with $N$ pixels. This is is only marginally more complex than computing the Euclidean distance which has $\mathcal{O}(N)$ complexity \cite{maurer2003linear}.
\end{rmk}

\subsection{Improvements to Geodesic Distance Term}
\begin{figure}[htb!]
\centering
\makebox[\textwidth][c]{
\includegraphics[width=1.2\textwidth, height = 6cm]{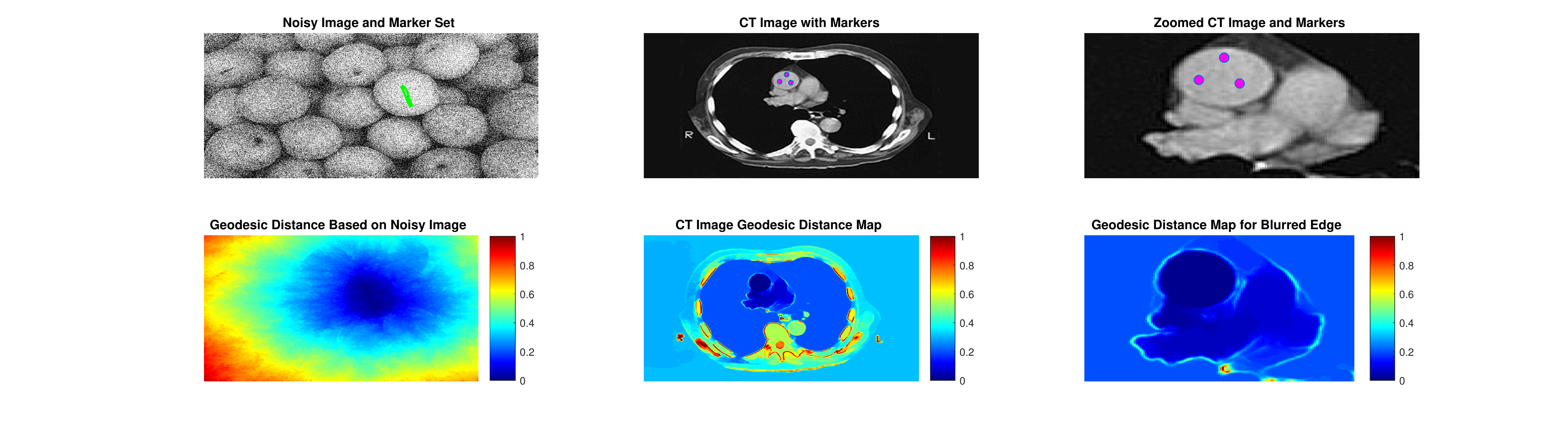}
}%
\caption{Examples of images showing the problems discussed and the resulting geodesic distance maps. Column $1$ shows the lack of robustness to noise, column $2$ shows that outside the patient we have unreasonably low distance penalty, column $3$ shows how the blurred edge under the aorta leads to the distance term being very low throughout the heart.\label{fig:DisadvantagesEPS}}
\end{figure}

We now propose some modifications to the geodesic distance.
Although the geodesic distance presents many advantages for selective image segmentation, we have three key disadvantages of this fitting term, which the Euclidean fitting term does not suffer.
\begin{enumerate}
\item{{\bf Not robust to noise.} The computation of the geodesic distance depends on $|\nabla z(x,y)|^{2}$ in $f(x,y)$ (see (\ref{eqn:DG0})). So, if an image contains a lot of noise, each noisy pixel appears as an edge and we get a misleading distance term.}
\item{{\bf Objects far from $\mathcal{M}$ with low penalty.} As the geodesic distance only uses marker set $\mathcal{M}$ for its initial condition (see (\ref{eqn:DG0})), this can result in objects far from $\mathcal{M}$ having a low distance penalty, which is clearly not desired.}
\item{{\bf Blurred edges.} If we have two objects separated by a blurry edge and we have marker points only in one object, the geodesic distance will be low to the other object, as the edge penalty is weakly enforced for a blurry edge. We would desire low penalty inside the object with markers and a reasonable penalty in the joined object.}
\end{enumerate}
In Figure~\ref{fig:DisadvantagesEPS}, each column shows an example for each of the problems listed above. We now propose solutions to each of these problems.

\begin{figure}[htb!]
\centering
\makebox[\textwidth][c]{
\includegraphics[width=1.2\textwidth, height = 8cm]{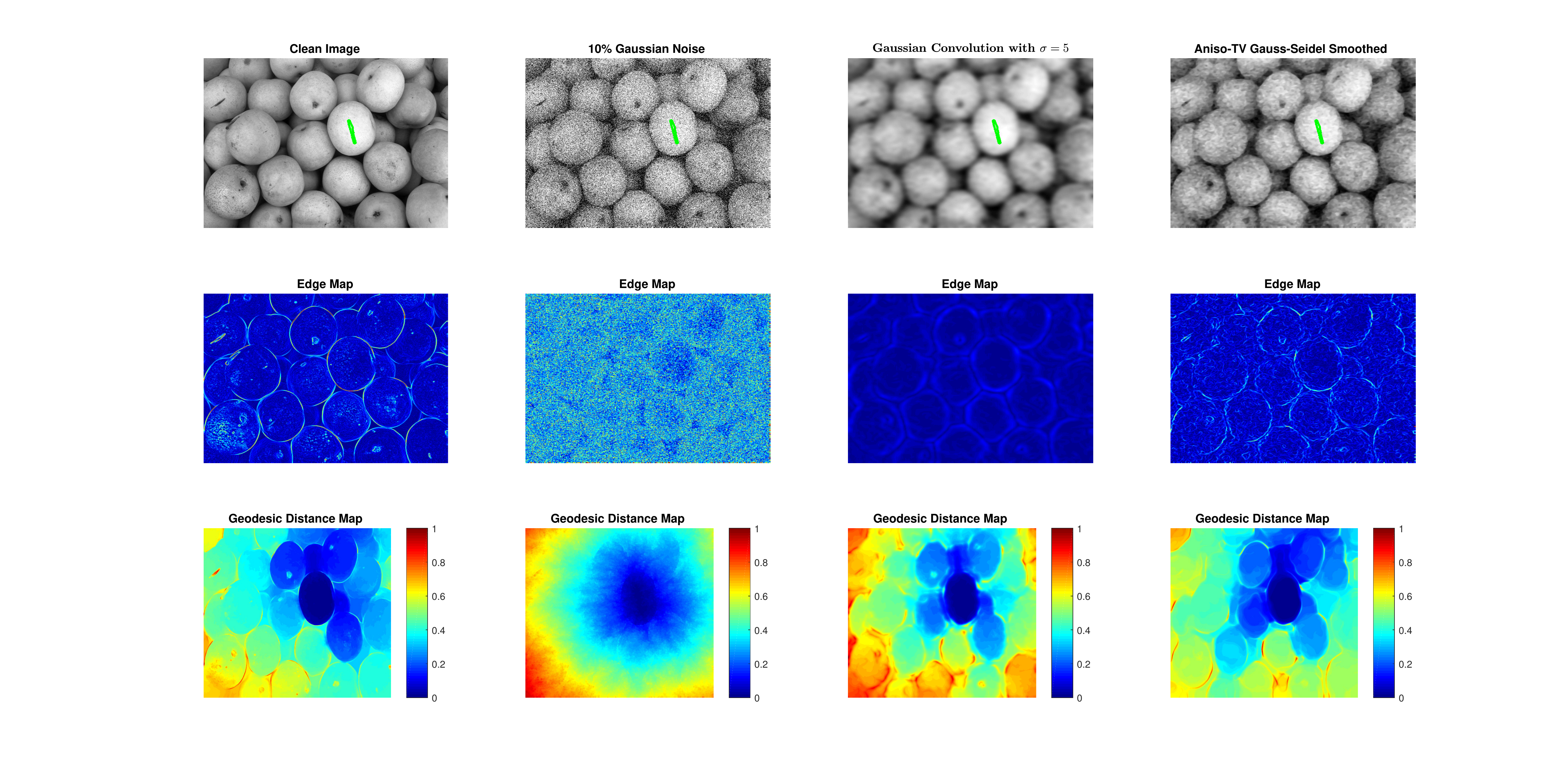}
}%
\caption{The edge maps and geodesic distance maps. (Left to right:) the clean image, the image with 10\% Gaussian noise, the noisy image with Gaussian convolution applied ($\sigma=5$) and for the noisy image with $100$ iterations of anisotropic-TV Gauss-Seidel smoothing. The set $\mathcal{M}$ is shown on the top row, it is the same for each image.\label{fig:prob1sol}}
\end{figure}

\begin{figure}[htb!]
\centering
\makebox[\textwidth][c]{
\includegraphics[width=1.3\textwidth, height = 8cm]{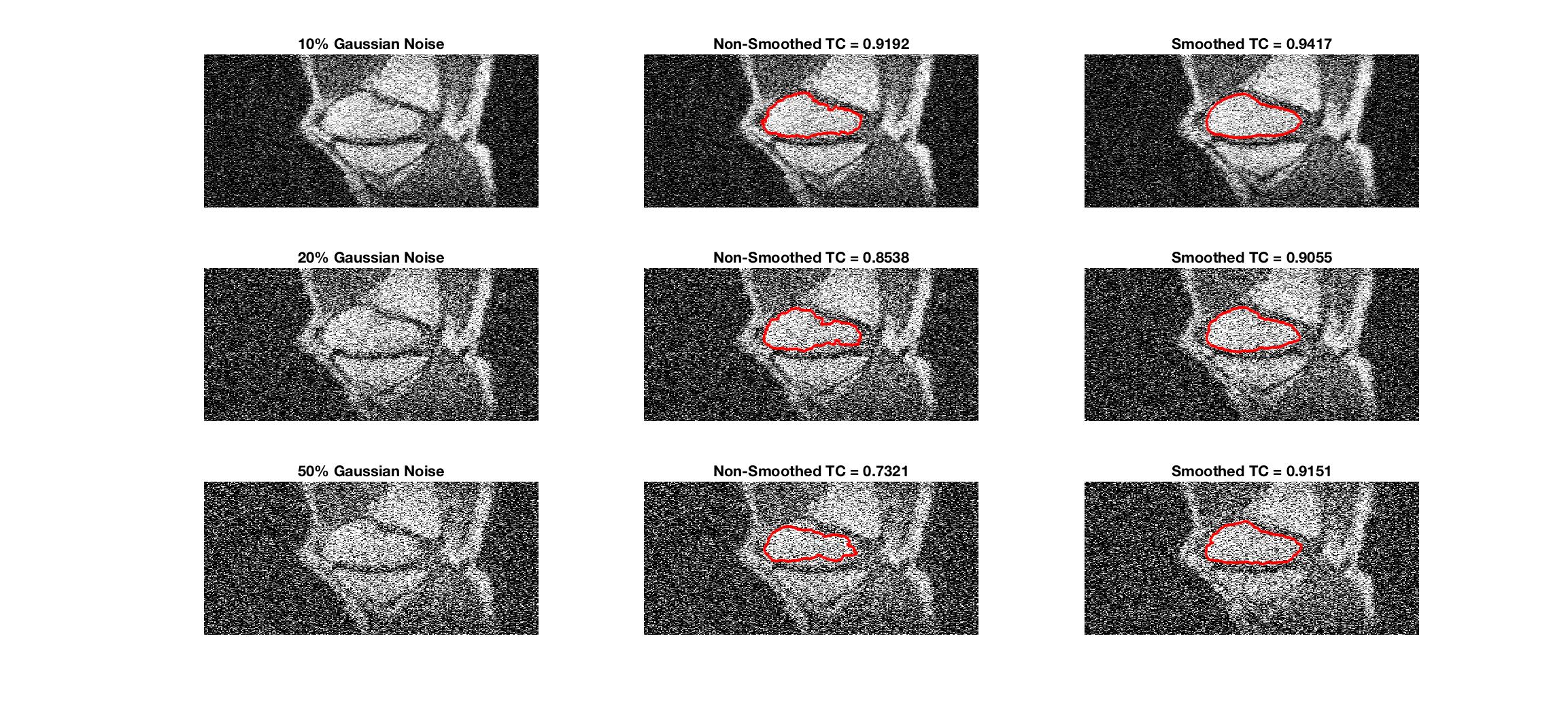}
}%
 \vspace{-0.2in}
\caption{Segmentation results and Tanimoto Coefficients (see \S\ref{sec:TC}) for images with 10\%, 20\% and 30\% Gaussian Noise with and without smoothing, $\lambda_{1} = \lambda_{2} = 5$, $\theta = 3$. \label{fig:SmoothedImages}}
\end{figure}

{\underline{Problem 1}: {\bf Noise Robustness.}} A na\"ive solution to the problem of noisy images would be to apply a Gaussian blur to $z(x,y)$ to remove the effect of the noise, so we change $f(x,y)$ to
\begin{equation}\label{eqn:edgemap}
\tilde{f}(x,y) = \varepsilon_{\mathcal{D}} + \beta_{G}|\nabla G_{\sigma} * z(x,y)|^{2}
\end{equation}
where $G_{\sigma}$ is a Gaussian convolution with standard deviation $\sigma$. However, the effect of Gaussian convolution is that it also blurs edges in the image. This then gives us the same issues described in Problem 3. We see in Figure~\ref{fig:prob1sol} column 3, that the Gaussian convolution reduces the sharpness of edges and this results in the geodesic distance being very similar in adjacent objects -- therefore we see more pixels with high geodesic distance.
Our alternative to Gaussian blur is to consider anisotropic TV denoising. We refer the reader to \cite{catte1992image,perona1990scale} for information on the model, here we just give the PDE which results from its minimisation:
\begin{equation}\label{eqn:anisoTV}
\begin{gathered}
\begin{aligned}
\tilde{\mu}\nabla\cdot \bigg( g(|\nabla z(x,y)|) \frac{\nabla u}{|\nabla u|_{\varepsilon_{2}}} \bigg) + \iota (z(x,y) - u) = 0,
\end{aligned}
\end{gathered}
\end{equation}
where $\tilde{\mu},\iota$ are non-negative parameters (we fix throughout $\tilde{\mu} = 10^{-3}, \iota = 5\times 10^{-4}$).
It is proposed to apply a relatively small number of cheap fixed point Gauss-Seidel iterations (between 100 and 200) to the discretised PDE. We cycle through all pixels $(i,j)$ and update $u_{i,j}$ as follows
\begin{equation}\label{eqn:anisoTVGSLEX}
\begin{gathered}
\begin{aligned}
u_{i,j} = \frac{A_{i,j}u_{i+1,j} + B_{i,j}u_{i-1,j} + C_{i,j}u_{i,j+1} + D_{i,j}u_{i,j-1}}{A_{i,j}+B_{i,j}+C_{i,j}+D_{i,j} + \iota}
\end{aligned}
\end{gathered}
\end{equation}
where $A_{i,j} = \frac{\tilde{\mu}}{h_{x}^{2}}g(|\nabla z(x,y)|)_{i+1/2,j}$, $B_{i,j} = \frac{\tilde{\mu}}{h_{x}^{2}}g(|\nabla z(x,y)|)_{i-1/2,j}$, $C_{i,j} = \frac{\tilde{\mu}}{h_{y}^{2}}g(|\nabla z(x,y)|)_{i,j+1/2}$ and $D_{i,j} = \frac{\tilde{\mu}}{h_{y}^{2}}g(|\nabla z(x,y)|)_{i,j-1/2}$. We update all pixels once per iteration and solve the PDE in (\ref{eqn:DG0}) with $f(x,y)$ replaced by
\begin{equation}\label{eqn:fhat}
{f}_{1}(x,y) = \varepsilon_{\mathcal{D}} + \beta_{G}|\nabla S^{k}(z(x,y))|^{2}
\end{equation}
where $S$ represents the Gauss-Seidel iterative scheme and $k$ is the number of iterations performed (we choose $k=100$ in our tests). In the final column of Figure~\ref{fig:prob1sol} we see that the geodesic distance map more closely resembles that of the clean image than the Gaussian blurred map in column 3 and in Figure~\ref{fig:SmoothedImages} we see that the segmentation results are qualitatively and quantitatively better using the anisotropic smoothing technique.

{\underline{Problem 2}: {\bf Objects far from $\mathcal{M}$ with low penalty.}}

\begin{figure}[htb!]
\centering
\includegraphics[width=\textwidth, height = 5cm]{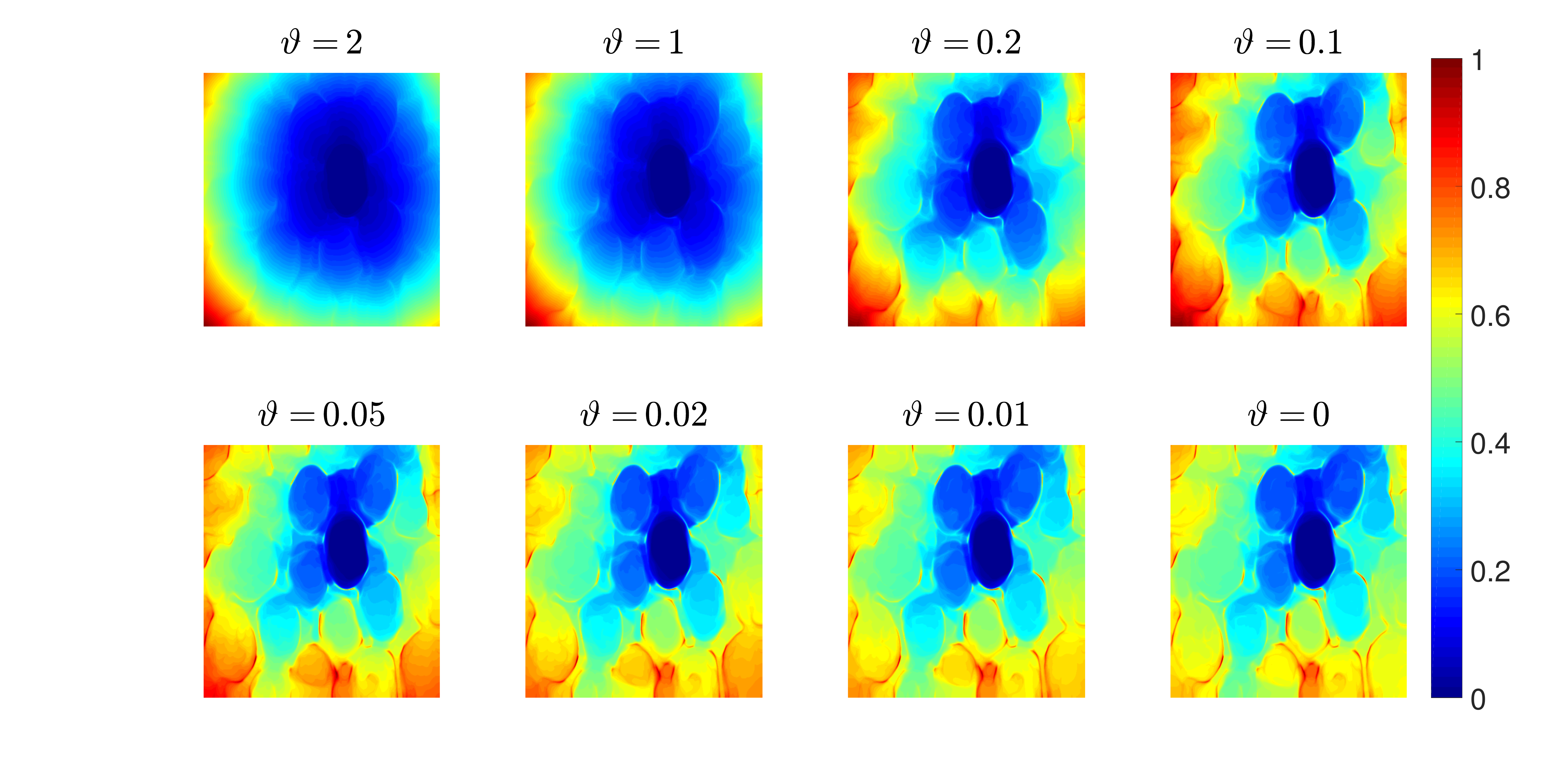}
\caption{Displayed is $\mathcal{D}_{M}(x,y)$ using ${f}_{2}(x,y)$ for various $\vartheta$ values. The marker set is the same as that used in Figure~\ref{fig:prob1sol}. \label{fig:prob2sol}}
\end{figure}

In Figure~\ref{fig:DisadvantagesEPS} column 2 we see that the geodesic distance to the outside of the patient is lower than to their ribs. This is due to the fact that the region outside the body is homogeneous and there is almost zero distance penalty in this region. Similarly for Figure~\ref{fig:prob1sol} column $4$, the distances from the marker set to many surrounding objects is low, even though their Euclidean distance from the marker set is high. We wish to have the Euclidean distance $\mathcal{D}_{E}(x,y)$ incorporated somehow. Our solution is to modify the term ${f}_{1}(x,y)$ from (\ref{eqn:fhat}) to
\begin{equation}\label{eqn:f2_final}
{f}_{2}(x,y) = \varepsilon_{\mathcal{D}} + \beta_{G}|\nabla S^{k}(z(x,y))|^{2} + \vartheta\mathcal{D}_{E}(x,y).
\end{equation}
In Figure~\ref{fig:prob2sol} the effect of this is clear, as $\vartheta$ increases, the distance function resembles the Euclidean distance more. We use value $\vartheta = 10^{-1}$ in all experiments as it adds a reasonable penalty to pixels far from the marker set.

{\underline{Problem 3}: {\bf Blurred edges.}}

If there are blurred edges between objects in an image, the geodesic distance will not increase significantly at this edge. Therefore the final segmentation result is liable to include unwanted objects. We look to address this problem through the use of anti-markers. These are markers which indicate objects that we do not want to segment, i.e. the opposite of marker points, we denote the set of anti-marker points by $\mathcal{AM}$.
\begin{figure}[htb!]
\centering
\makebox[\textwidth][c]{
\includegraphics[width=1.3\textwidth, height = 4cm]{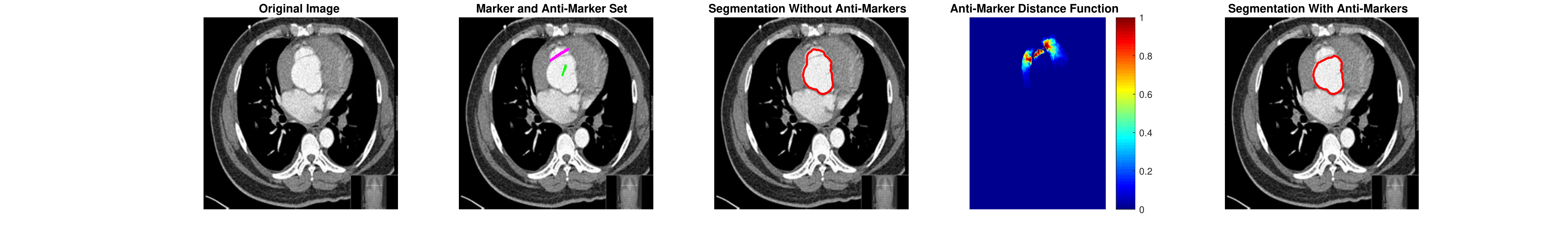}
}%
\caption{(Left to right:) original image, $\mathcal{M}$ (green) and $\mathcal{AM}$ (pink), segmentation result just using marker set, $\mathcal{D}_{AM}(x,y)$ using anti-markers, segmentation result using anti-markers. For these $\mu = 1, \lambda_{1}=\lambda_{2}=5, \theta = 25$. \label{fig:AntiMarkerDist}}
\end{figure}
We propose to use a geodesic distance map from the set $\mathcal{AM}$ denoted by $\mathcal{D}_{AM}(x,y)$ which penalises pixels near to the set $\mathcal{AM}$ and doesn't add any penalty to those far away. We could na\"ively choose $\mathcal{D}_{AM}(x,y) = 1 - \tilde{\mathcal{D}}_{GAM}(x,y)$ where $\tilde{\mathcal{D}}_{GAM}(x,y)$ is the normalised geodesic distance from $\mathcal{AM}$. However this puts a large penalty on those pixels inside the object we actually want to segment (as $\tilde{\mathcal{D}}_{GAM}(x,y)$ to those pixels is small). To avoid this problem, we propose the following anti-marker distance term
\[
\mathcal{D}_{AM}(x,y) = \frac{\exp\left(-\tilde{\alpha}\tilde{\mathcal{D}}_{GAM}(x,y)\right)-\exp\left(-\tilde{\alpha}\right)}{1-\exp\left(-\tilde{\alpha}\right)}
\]
where $\tilde{\alpha}$ is a tuning parameter. We choose $\tilde{\alpha} = 200$ throughout. This distance term ensures rapid decay of the penalty away from the set $\mathcal{AM}$ but still enforces high penalty around the anti-marker set itself. See Figure~\ref{fig:AntiMarkerDist} where a segmentation result with and without anti-markers is shown. As $\mathcal{D}_{AM}(x,y)$ decays rapidly from $\mathcal{AM}$, we do require that the anti-marker set be close to the blurred edge and away from the object we desire to segment.

\subsection{The new model and its Euler-Lagrange equation}

{\bf The Proposed Geodesic Model.} Putting the above $3$ ingredients together, we propose the model
\vspace{-0.2in}
\begin{equation}
\resizebox{0.92\textwidth}{!}{$
\begin{split}
\min_{u,c_{1},c_{2}} \Big\{ F_{GEO}(u,c_{1},c_{2})=& \int_{\Omega}\left[   \lambda_{1}(z(x,y)-c_{1})^{2}  - \lambda_{2}(z(x,y)-c_{2})^{2}  \right]  u\,\mathrm{d}\Omega \\
  & + \mu\int_{\Omega}g(|\nabla z(x,y)|)
|\nabla u|\,\mathrm{d}\Omega
    + \theta\int_{\Omega}\mathcal{D}_{G}(x,y)u \,\mathrm{d}\Omega + \alpha\int_{\Omega}\nu_{\varepsilon}(u)\,\mathrm{d}\Omega \Big\},
\end{split}$}\label{eqn:geofunc2F}
\end{equation}
where $\mathcal{D}_{G}(x,y) = \left(\mathcal{D}_{M}(x,y) + \mathcal{D}_{AM}(x,y)\right)/2$ and $\mathcal{D}_{M}(x,y)$ is the geodesic distance from the marker set $\mathcal{M}$. We compute $\mathcal{D}_{M}(x,y)$ using (\ref{eqn:DG0}) where $f(x,y) = f_{2}(x,y)$ defined in (\ref{eqn:f2_final}). Using Calculus of Variations, solving (\ref{eqn:geofunc2F}) with respect to $c_{1}, \ c_{2}$, with $u$ fixed, leads to
\begin{equation}\label{eqn:geoc1c2}
c_{1}(u)=\frac{\int_{\Omega}u\cdot z(x,y)\,\mathrm{d}\Omega}{\int_{\Omega}u \,\mathrm{d}\Omega},
\hspace{0.5in}
c_{2}(u)=\frac{\int_{\Omega}(1-u)\cdot z(x,y)\,\mathrm{d}\Omega}{\int_{\Omega}(1-u) \,\mathrm{d}\Omega},
\end{equation}
and the minimisation with respect to $u$ (with $c_{1}$ and $c_{2}$ fixed) gives the PDE
\begin{equation} 
\begin{gathered}
\begin{aligned}
\mu\nabla\cdot \left( g(|\nabla z(x,y)|) \frac{\nabla u}{|\nabla u|_{\varepsilon_{2}}} \right) - \Big[ \lambda_{1} (z(x,y)-c_{1})^{2} &- \lambda_{2}(z(x,y) - c_{2})^{2}\Big]
\\
&- \theta\mathcal{D}_{G}(x,y) - \alpha\nu'_{\varepsilon}(u)= 0
\end{aligned}
\end{gathered}\label{eqn:geoel2}
\end{equation}
in $\Omega$, where we replace $|\nabla u|$ with $|\nabla u|_{\varepsilon_{2}} = \sqrt{u_{x}^{2}+u_{y}^{2}+\varepsilon_{2}}$ to avoid zero denominator; we choose $\varepsilon_{2}=10^{-6}$ throughout. We also have Neumann boundary conditions $\frac{\partial u}{\partial {\bm n}} =  0$ on $\partial\Omega$
where ${\bm n}$ is the outward unit normal vector.

Next we discuss a numerical scheme for solving this PDE \eqref{eqn:geoel2}.
However it should be remarked that updating $c_{1}(u), c_{2}(u)$ should be done as soon as $u$ is updated; practically $c_1, c_2$ converge very quickly since
the object intensity $c_1$ does not change much. 


\section{An additive operator splitting algorithm}
\label{sec:AOS}

\begin{figure}[htb!]
\centering
\makebox[\textwidth][c]{
  \begin{subfigure}[b]{0.3\textwidth}
\includegraphics[width=\textwidth,height = 4cm]{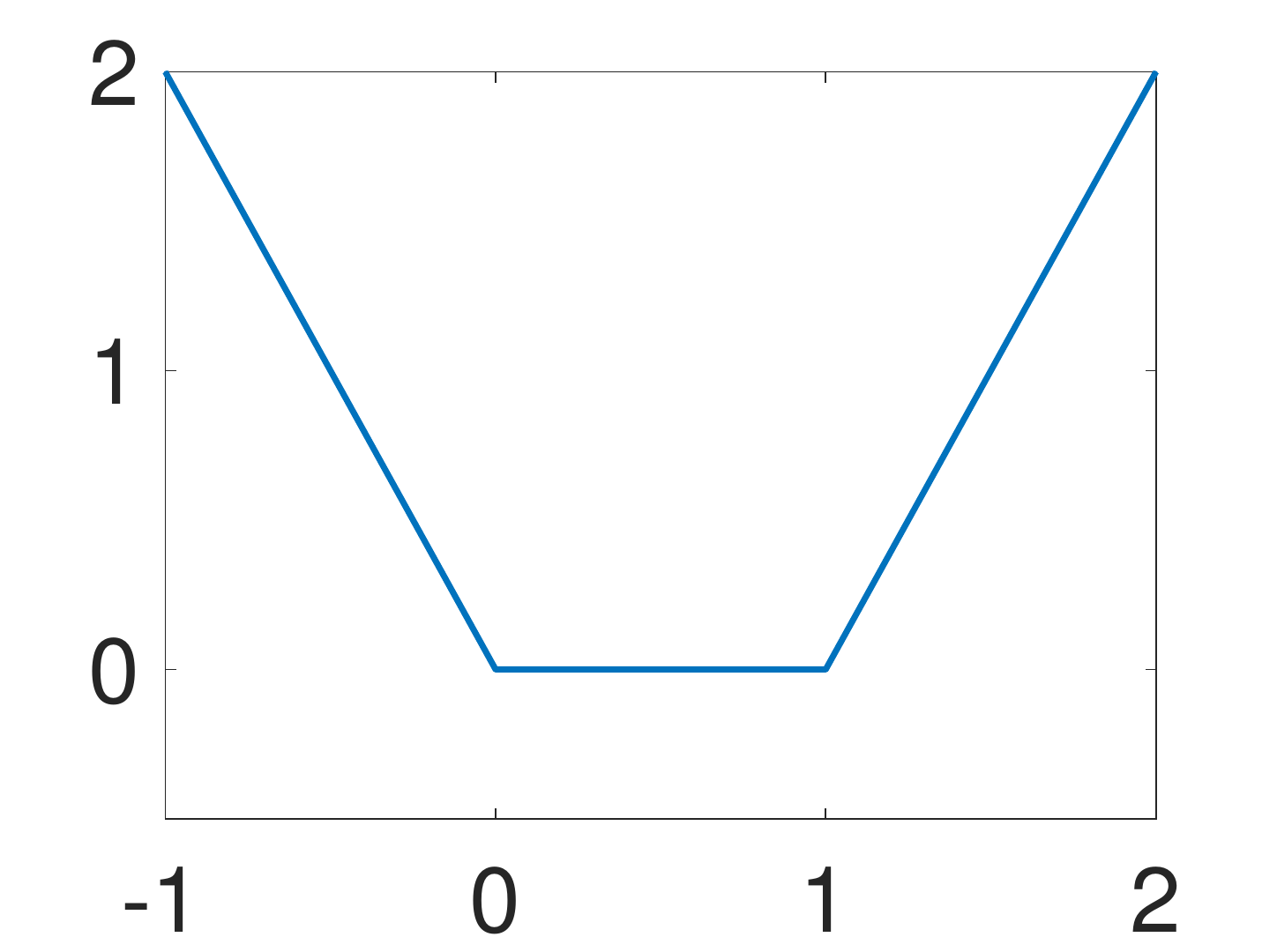}
        \caption{$\nu(u)$.}
    \end{subfigure}
      \begin{subfigure}[b]{0.3\textwidth}
\includegraphics[width=\textwidth,height = 4cm]{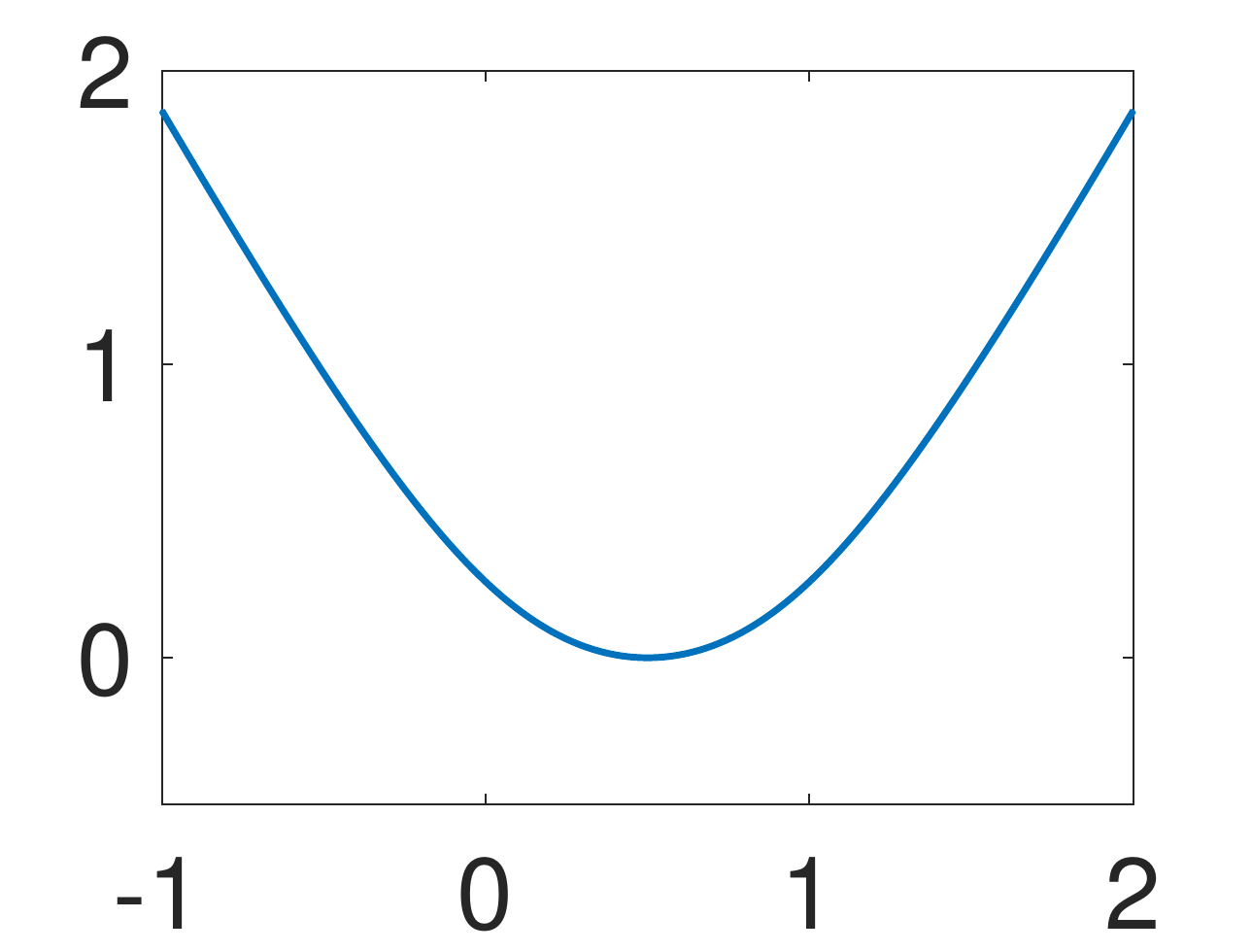}
        \caption{$\nu_{\varepsilon}(u)$ for $\varepsilon=1$.}
    \end{subfigure}
          \begin{subfigure}[b]{0.3\textwidth}
\includegraphics[width=\textwidth,height = 4cm]{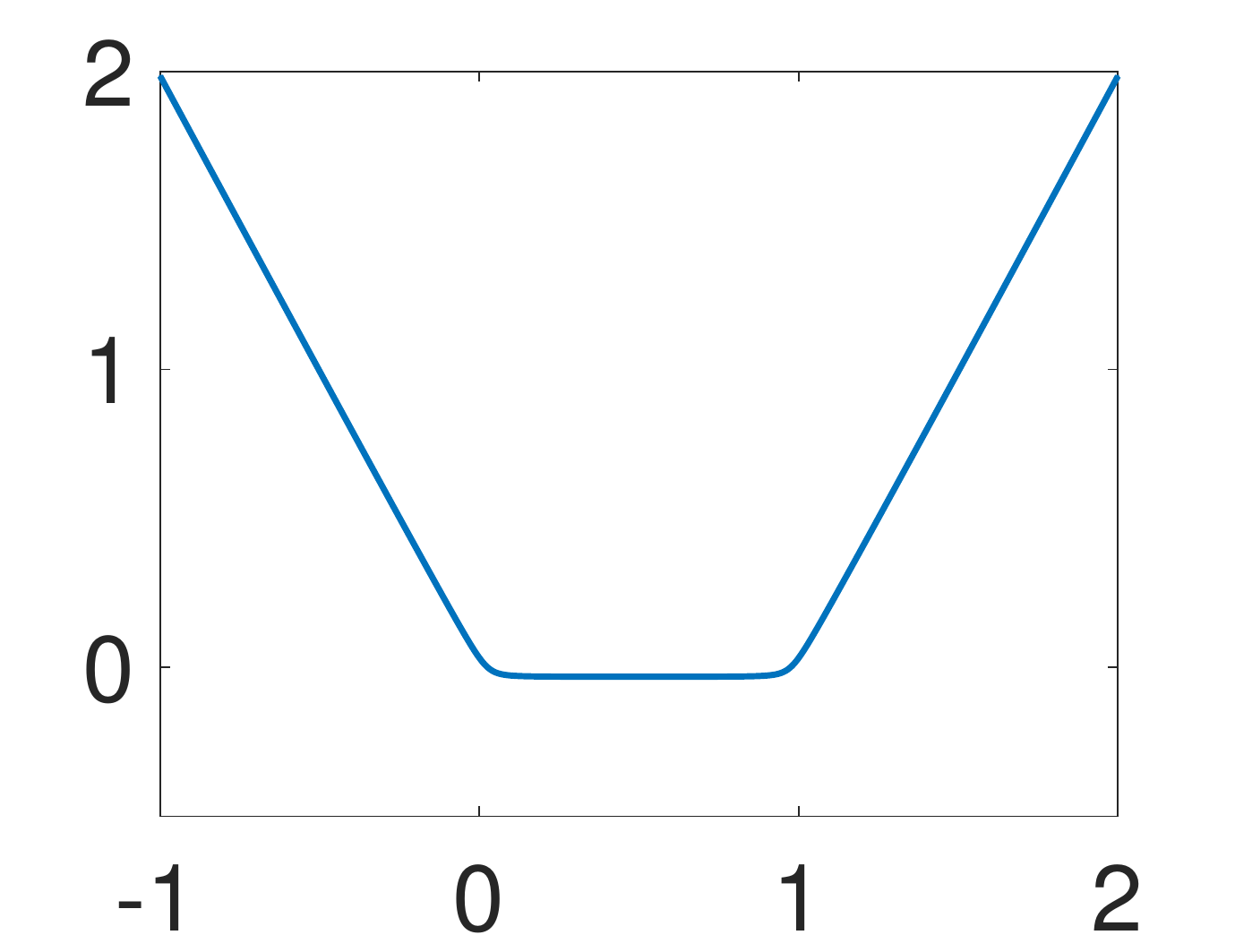}
        \caption{$\nu_{\varepsilon}(u)$ for $\varepsilon=0.1$.}
    \end{subfigure}    }%
    \caption{(a) The exact penalty function $\nu(u)$ and (b,c) $\nu_{\varepsilon}(u)$ for different $\varepsilon$ values. \label{fig:vepsplots}}
\end{figure}

\begin{figure}[htb!]
\centering
\makebox[\textwidth][c]{
  \begin{subfigure}[b]{0.3\textwidth}
\includegraphics[width=\textwidth,height = 4cm]{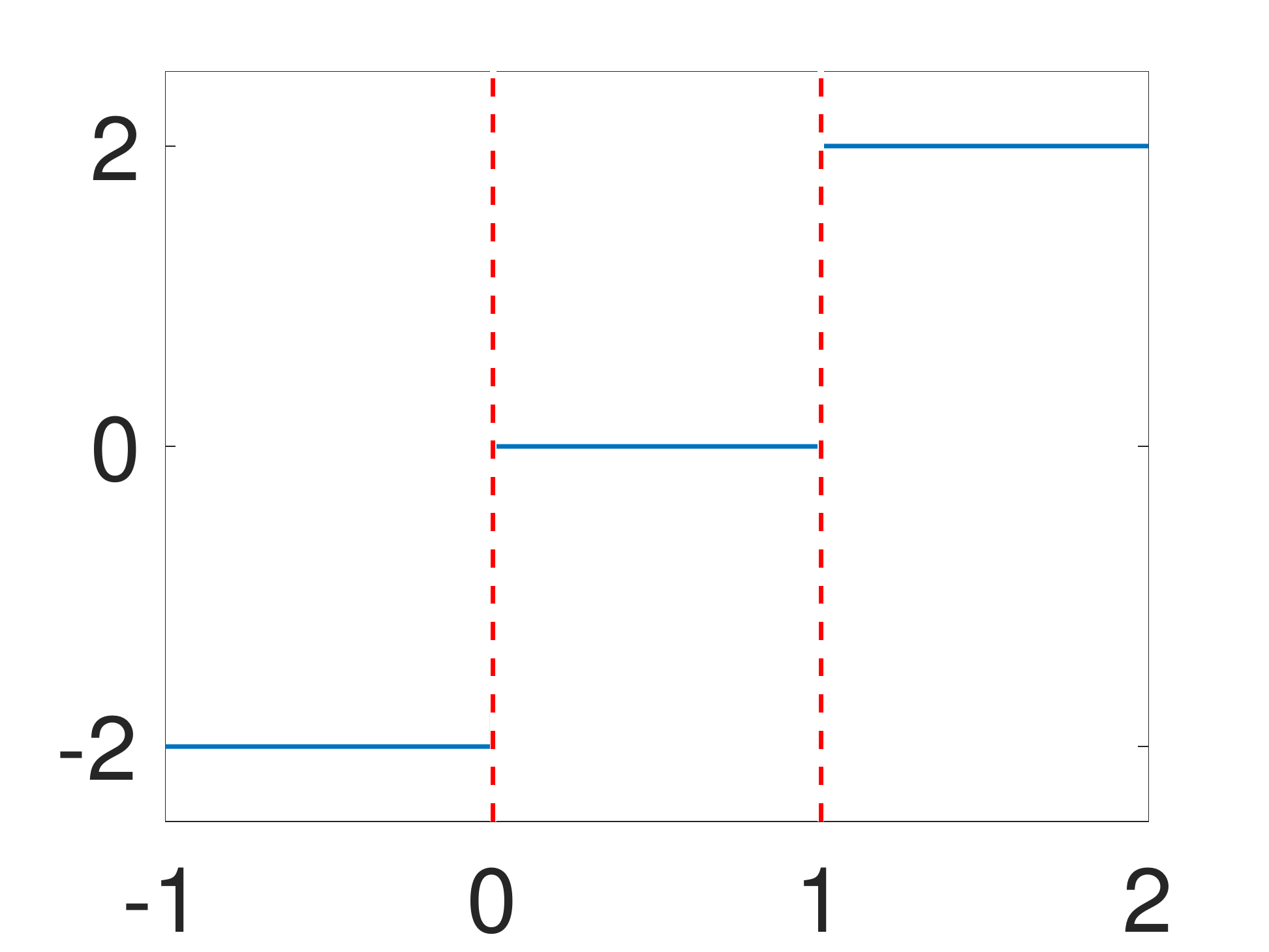}
        \caption{$\nu'(u)$.}
    \end{subfigure}
      \begin{subfigure}[b]{0.3\textwidth}
\includegraphics[width=\textwidth,height = 4cm]{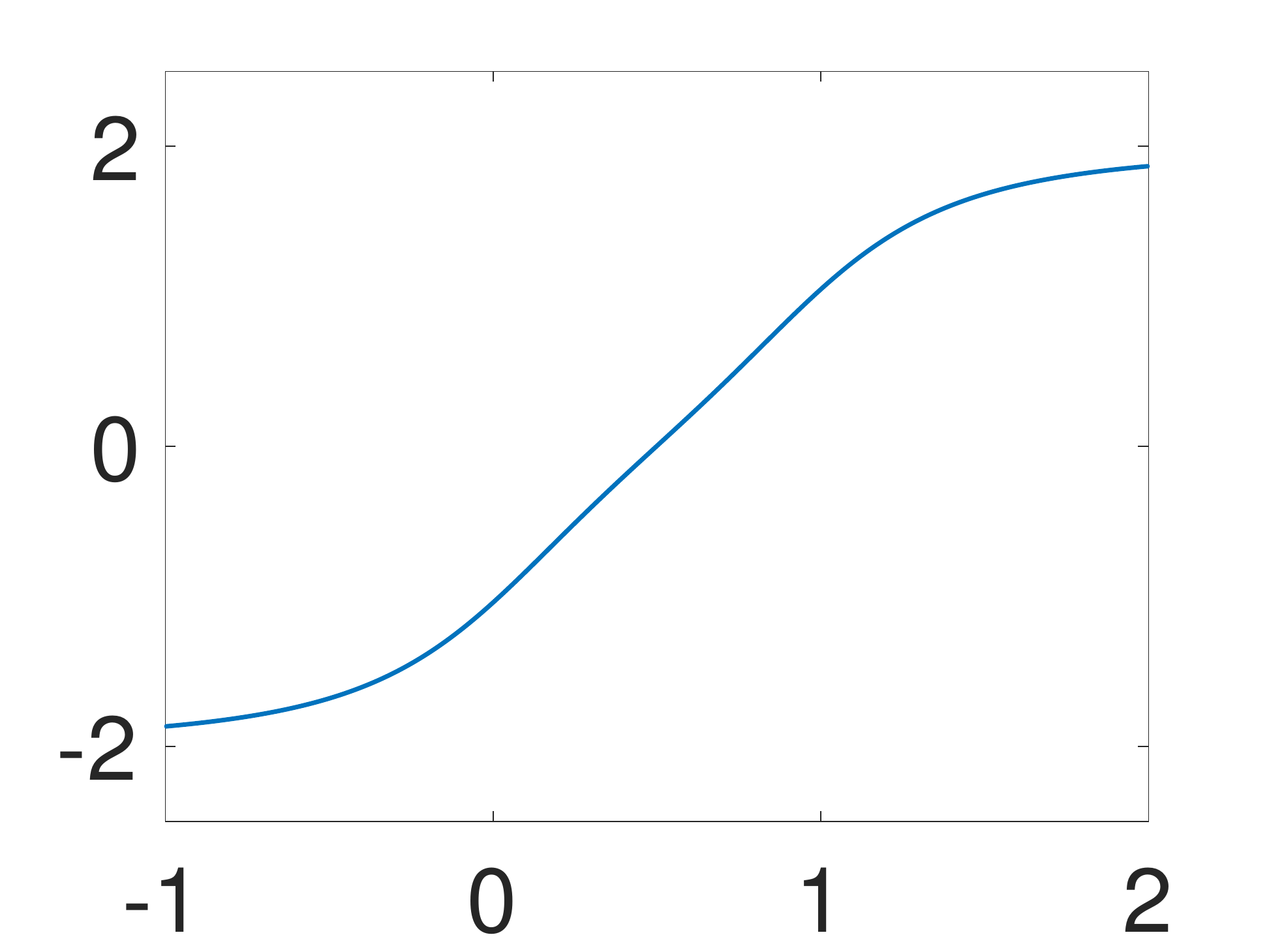}
        \caption{$\nu'_{\varepsilon}(u)$ for $\varepsilon=1$.}
    \end{subfigure}
          \begin{subfigure}[b]{0.3\textwidth}
\includegraphics[width=\textwidth,height = 4cm]{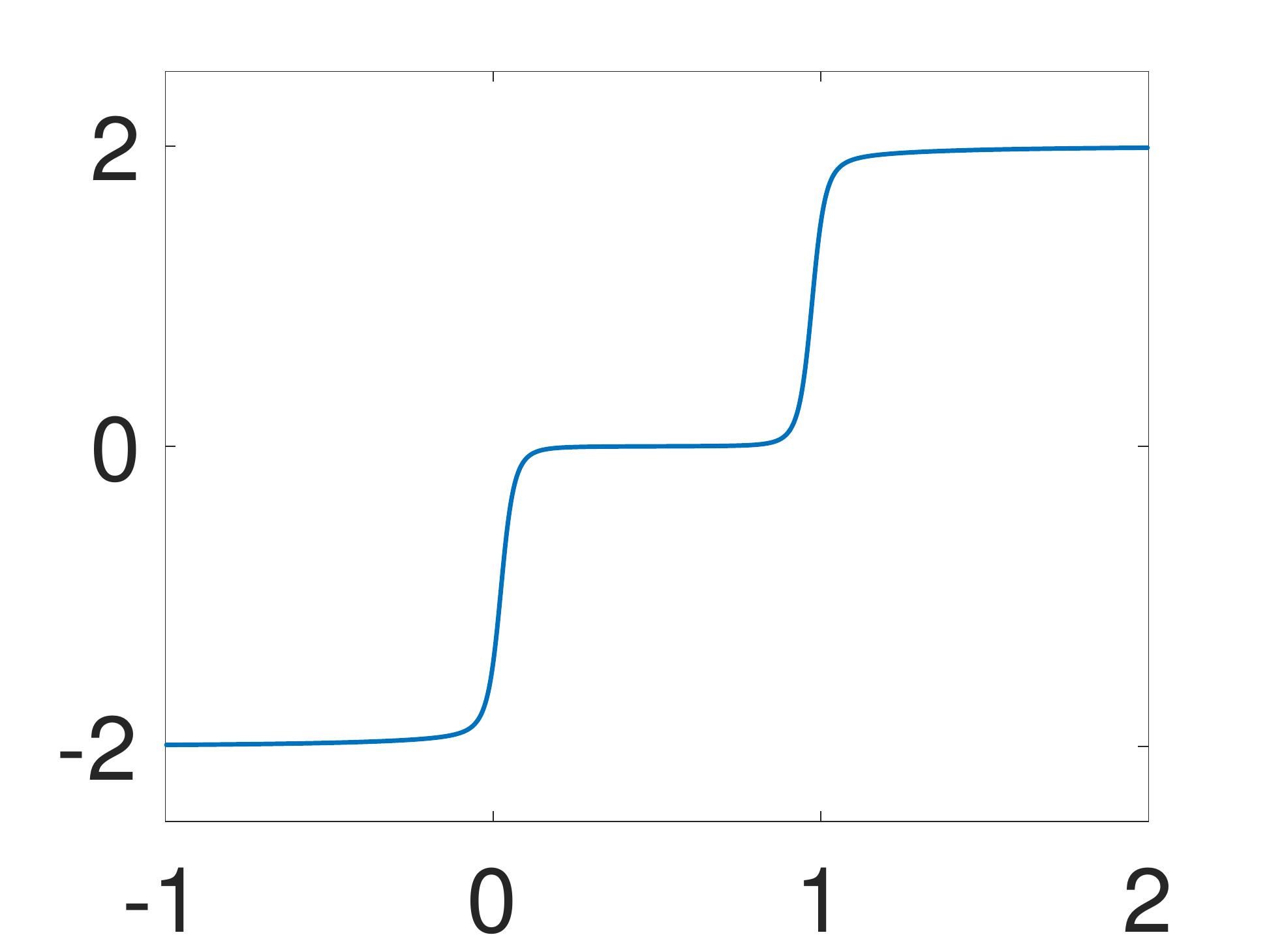}
        \caption{$\nu'_{\varepsilon}(u)$ for $\varepsilon=0.1$.}
    \end{subfigure}    }%
        \caption{(a) $\nu'(u)$ (discontinuities shown in red) and (b,c) $\nu'_{\varepsilon}(u)$ for different $\varepsilon$ values. \label{fig:vepsprimeplots}}
\end{figure}

Additive Operator Splitting (AOS) is a widely used method \cite{Gord74,LU199125,BC36} as seen from more recent works
\cite{BC1, BC2, SC1,SC2,SC28, SC}  on the diffusion type equation such as
\begin{equation}\label{eqn:AOSeqn}
\frac{\partial u}{\partial t} = \mu\nabla\cdot (G(u)\nabla u) - f.
\end{equation}
AOS allows us to split the two dimensional problem into two one-dimensional problems, which we solve and then combine. Each one dimensional problem gives rise to a tridiagonal system of equations which can be solved efficiently, hence AOS is a very efficient method for solving diffusion-like equations. AOS is a semi-implicit method and permits far larger time-steps than the corresponding explicit schemes would. Hence AOS is more stable than an explicit method \cite{BC36}. We rewrite the above equation as
\[
\frac{\partial u}{\partial t} = \mu \bigg( \partial_{x}\left(G(u)\partial_{x} u\right)+\partial_{y}\left(G(u)\partial_{y} u\right)\bigg) - f.
\]
and after discretisation, we can rewrite this as \cite{BC36}
\[
u^{k+1} = \frac{1}{2}\sum_{\ell=1}^{2}\bigg(I-2\tau\mu A_{\ell}(u^{k})\bigg)^{-1}\left(u^{k}+\tau f\right)
\]
where $\tau$ is the time-step, $A_{1}(u) = \partial_{x}(G(u)\partial_{x})$ and $A_{2}(u) = \partial_{y}(G(u)\partial_{y})$. For notational convenience we write $G = G(u)$. The matrix $A_{1}(u)$ can be obtained as follows
\begin{equation*}
\resizebox{\textwidth}{!}{$
\begin{gathered}
\begin{aligned}
\bigg( A_{1}(u^{k})u^{k+1}\bigg)_{i,j} = \bigg( \partial_{x}\Big(G \partial_{x} u^{k+1}\Big)\bigg)_{i,j}
=& \Bigg( \frac{G_{i+\frac{1}{2},j}}{h_{x}^{2}} \Bigg) u^{k+1}_{i+1,j} + \Bigg( \frac{G_{i-\frac{1}{2},j}}{h_{x}^{2}} \Bigg) u^{k+1}_{i-1,j}
- \Bigg( \frac{G_{i+\frac{1}{2},j}+G_{i-\frac{1}{2},j}}{h_{x}^{2}}  \Bigg) u^{k+1}_{i,j}
\end{aligned}
\end{gathered}
$}
\end{equation*}
and similarly to \cite{SC28,SC}, for the half points in $G$ we take the average of the surrounding pixels, e.g. $G_{i+\frac{1}{2},j} = \frac{G_{i+1,j} +G_{i,j} }{2}$.
Therefore we must solve two tridiagonal systems to obtain $u^{k+1}$, the Thomas algorithm allows us to solve each of these efficiently \cite{BC36}. The AOS method described here assumes $f$ does not depend on $u$, however in our case it depends on $\nu'_{\varepsilon}(u)$ (see (\ref{eqn:geoel2})) which has jumps around 0 and 1, so the algorithm has stability issues. This was noted in \cite{SC} and the authors adapted the formulation of (\ref{eqn:geoel2}) to offset the changes in $f$. Here we repeat their arguments for adapting AOS when the exact penalty term $\nu'_{\varepsilon}(u)$ is present (we refer to Figures~\ref{fig:vepsplots} and \ref{fig:vepsprimeplots} for plots of the penalty function and its derivative, respectively).

The main consideration is to extract a linear part out of the nonlinearity in $f=f(u)$.
If we evaluate the Taylor expansion of $\nu'_{\varepsilon}(u)$ around $u=0$ and $u=1$ and group the terms into the constant and linear components in $u$, we can respectively write $\nu'_{\varepsilon}(u) = a_{0}(\varepsilon) + b_{0}(\varepsilon) u + \mathcal{O}(u^{2})$ and $\nu'_{\varepsilon}(u) = a_{1}(\varepsilon)  + b_{1}(\varepsilon) u + \mathcal{O}(u^{2})$. We actually find that $b_{0}(\varepsilon) =b_{1}(\varepsilon) $ and denote the linear term as $b$ from now on. Therefore, for a change in $u$ of $\delta u$ around $u=0$ and $u=1$, we can approximate the change in $\nu'_{\varepsilon}(u)$ by $b\cdot\delta u$.
To focus on the jumps, 
define the interval in which $\nu'_{\varepsilon}(u)$ jumps as
\[
I_{\zeta} := [0-\zeta,0+\zeta]\cup [1-\zeta,1+\zeta]
\]
and refine the linear function by
\[
\tilde{b}_{i,j}^{k} =
\begin{cases}
b, & u_{i,j}^{k} \in I_{\zeta}\\
0, & else.
\end{cases}
\]
Using these we can now offset the change in $\nu'_{\varepsilon}(u^{k})$ by changing the formulation (\ref{eqn:AOSeqn}) to
\[
\frac{\partial u}{\partial t} = \mu\nabla\cdot (G(u)\nabla u)
 - \alpha\tilde{b}^{k} u    +  \big[ \alpha\tilde{b}^{k}u  - f\big]
\]
or in AOS form
\(
u^{k+1} =u^{k}+\tau\mu\nabla\cdot (G(u^{k})\nabla u^{k+1}) 
- \tau\alpha\tilde{b}^{k} u^{k+1}
+ \big[ \tau\alpha\tilde{b}^{k}u^{k} - f^{k} \big]
\)
which, following the derivation in \cite{SC}, can be reformulated as
\[
u^{k+1} = \frac{1}{2}\sum_{\ell=1}^{2}{\underbrace{\bigg(I+\tilde{B}^{k}-2\tau\mu A_{\ell}\left(u^{k}\right)\bigg)}_{Q_{1}}}^{-1}\left(\left(I+\tilde{B}^{k}\right) u^{k}+\tau f^{k}\right)
\]
where $\tilde{B}^{k} = \text{diag}(\tau\alpha\tilde{b}^{k})$. We note that $Q_{1}$ is invertible as it is strictly diagonally dominant. This scheme improves on (\ref{eqn:AOSeqn}) as now, changes in $f^{k}$ are damped. However, it is found in \cite{SC} that although this scheme does satisfy most of the discrete scale space conditions of Weickert \cite{BC36} (which guarantee convergence of the scheme), it does not satisfy all of them. In particular the matrix $Q_{1}$ doesn't have unit row sum and is not symmetrical. The authors adapt the scheme above to the equivalent
\begin{equation}\label{eqn:AOSgeo}
u^{k+1} = \frac{1}{2}\sum_{\ell=1}^{2}{\underbrace{\bigg(I-2\tau\mu \left(I+\tilde{B}^{k}\right)^{-1}A_{\ell}\left(u^{k}\right)\bigg)}_{Q_{2}}}^{-1}\left(u^{k}+\tau \left(I+\tilde{B}^{k}\right)^{-1}f^{k}\right),
\end{equation}
 where the matrix $Q_{2}$ does have unit row sum, however the matrix is not always symmetrical. We can guarantee convergence for $\zeta = 0.5$ (in which case $Q_{2}$ must be symmetrical) but we desire to use a small $\zeta$ to give a small interval $I_{\zeta}$.
We find experimentally that convergence is achieved for any small value of $\zeta$, this is due to the fact that at convergence the solution $u$ is almost binary \cite{SC8}. Therefore, although initially $Q_{2}$ is asymmetrical at some pixels, at convergence all pixels have values which fall within $I_{\zeta}$ and $I+\tilde{B}^{k}$ is a matrix with all diagonal entries $1+\tau\alpha b$. Therefore we find that at convergence $Q_{2}$ is symmetrical and the discrete scale space conditions are all satisfied. In all of our tests we fix $\zeta  = 0.01$.
\begin{megaalgorithm}[htb!]
\caption{Solution of the Geodesic Model}
\label{alg:AOS}
\begin{algorithmic}
\STATE Set $\mu,\lambda,\theta$. Compute $g(|\nabla z(x,y)|) = \frac{1}{1+\beta_{G}|\nabla z(x,y)|^{2}}$ and $\mathcal{D}_{G}(x,y) = \frac{\mathcal{D}_{G}^{0}(x,y)}{||\mathcal{D}_{G}^{0}(x,y)||_{L^{\infty}}}$,
\STATE  with $\mathcal{D}_{G}^{0}(x,y)$ the solution of (\ref{eqn:DG0}). Initialise $u^{(0)}$ arbitrarily.
\FOR{$iter=1$ to $max\_iterations$}
\STATE Calculate $c_{1}$ and $c_{2}$ using (\ref{eqn:geoc1c2}).
\STATE Calculate $r =  \lambda_{1} (z-c_{1})^{2} - \lambda_{2}(z - c_{2})^{2} + \theta\mathcal{D}_{G}$.
\STATE Set $\alpha = || r ||_{L^{\infty}}$.
\STATE Calculate $f^{k} =  r + \alpha\nu'_{\varepsilon}(u^{k})$.
\STATE Update $u^{k}$ to $u^{k+1}$ using the AOS scheme (\ref{eqn:AOSgeo}).
\ENDFOR
\STATE $u^{*} \leftarrow u^{k}$.
\end{algorithmic}
\end{megaalgorithm}


\section{Existence and Uniqueness of the Viscosity Solution}

In this section we use the viscosity solution framework and the work of Ishii and Sato \cite{ishiisato} to prove that, for a class of PDEs in image segmentation, the solution exists and is unique. In particular, we prove the existence and uniqueness of the viscosity solution for the PDE which is determined by the Euler-Lagrange equation for the Geodesic Model. Throughout, we will assume $\Omega$ is a bounded domain with $C^{1}$ boundary.

From the work of \cite{crandall,ishiisato}, we have the following Theorem for analysing the solution of a partial differential equation of the form $F(\bm{x},u,Du,D^{2}u)=0$ where $F: \mathbb{R}^{n}\times\mathbb{R}\times\mathbb{R}^{n}\times\mathscr{M}^{n}\rightarrow\mathbb{R}$, $\mathscr{M}^{n}$ is the set of $n\times n$ symmetric matrices, $Du$ is the gradient of $u$ and $D^{2}u$ is the Hessian of $u$. For simplicity, and in a slight abuse of notation, we use $x := \bm{x}$ for the vector of a general point in $\mathbb{R}^{n}$.

\begin{thm}[Theorem 3.1 \cite{ishiisato}] Assume that the following conditions (C1)--(C2) and (I1)--(I7) hold. Then for each $u_{0}\in C(\overline{\Omega})$ there is a unique viscosity solution $u\in C([0,T)\times\overline{\Omega})$ of (\ref{eqn:parabpde}) and (\ref{eqn:parabpdeBC}) satisfying (\ref{eqn:parabpdeg}). \label{thm:satothm}
\vspace{-0.3in}
\end{thm}
\begin{equation}\label{eqn:parabpde}
\frac{\partial u}{\partial t} + F(t,x,u,Du,D^{2}u) = 0 \hspace{0.2in} \text{in}  \hspace{0.2in} Q=(0,T)\times\Omega,
\end{equation}
\vspace{-0.0875in}
\begin{equation}\label{eqn:parabpdeBC}
B(x,Du) = 0\hspace{0.2in} \text{in} \hspace{0.2in}S=(0,T)\times\partial\Omega,
\end{equation}
\vspace{-0.0875in}
\begin{equation}\label{eqn:parabpdeg}
u(0,x) = u_{0}(x) \hspace{0.2in}\text{for} \hspace{0.2in}x\in\overline{\Omega}.
\end{equation}

{\bf Conditions (C1)--(C2).}
\begin{enumerate}
\item[(C1)]{$F(t,x,u,p,X) \le F(t,x,v,p,X)$ for $u\le v$.}
\item[(C2)]{$F(t,x,u,p,X) \le F(t,x,u,p,Y)$ for $X,Y\in\mathscr{M}^{n}$ and $Y\le X$.}
\end{enumerate}

{\bf Conditions (I1)--(I7).} Assume $\Omega$ is a bounded domain in $\mathbb{R}^{n}$ with $C^{1}$ boundary.
\begin{enumerate}
\item[(I1)]{$F \in C \left( [0,T] \times \overline{\Omega} \times \mathbb{R} \times \left( \mathbb{R}^{n}\backslash\{ 0 \} \right) \times \mathscr{M}^{n} \right)$.}
\item[(I2)]{There exists a constant $\gamma\in\mathbb{R}$ such that for each $(t,x,p,X)\in[0,T]\times\overline{\Omega}\times\left( \mathbb{R}^{n}\backslash\{ 0\} \right) \times\mathscr{M}^{n}$ the function $u\mapsto F(t,x,u,p,X) - \gamma u$ is non-decreasing on $\mathbb{R}$.}
\item[(I3)]{$F$ is continuous at $(t,x,u,0,0)$ for any $(t,x,u)\in[0,T]\times\overline{\Omega}\times\mathbb{R}$ in the sense that
\[
-\infty < F_{*}(t,x,u,0,0) = F^{*}(t,x,u,0,0)<\infty
\]
holds. Here $F^{*}$ and $F_{*}$ denote, respectively, the upper and lower semi-continuous envelopes of $F$, which are defined on $[0,T]\times\overline{\Omega}\times\mathbb{R}\times\mathbb{R}^{n}\times\mathscr{M}^{n}$.}
\item[(I4)]{$B\in C\left(\mathbb{R}^{n} \times\mathbb{R}^{n} \right)\cap C^{1,1}\left( \mathbb{R}^{n} \times \left(\mathbb{R}^{n}\backslash\{ 0 \}\right)  \right)$, where $C^{1,1}$ is the H{\"o}lder functional space.}
\item[(I5)]{For each $x\in\mathbb{R}^{n}$ the function $p\mapsto B(x,p)$ is positively homogeneous of degree one in $p$, i.e. $B(x,\lambda p) = \lambda B(x,p)$ for all $\lambda\ge 0$ and $p\in\mathbb{R}^{n}\backslash\{ 0\}$.}
\item[(I6)]{There exists a positive constant $\Theta$ such that $\langle {\bm{n}}(x), D_{p}B(x,p)\rangle\ge\Theta$ for all $x\in\partial\Omega$ and $p\in\mathbb{R}^{n}\backslash\{ 0\}$. Here ${\bm{n}}(x)$ denotes the unit outward normal vector of $\Omega$ at $x\in\partial\Omega$.}
\item[(I7)]{For each $R>0$ there exists a non-decreasing continuous function $\omega_{R}:[0,\infty)\rightarrow[0,\infty)$ satisfying $\omega_{R}(0)=0$ such that if $X,Y\in\mathscr{M}^{n}$ and $\mu_{1},\mu_{2}\in[0,\infty)$ satisfy
\begin{equation}\label{eqn:assumpI3}\begin{bmatrix}
X & 0 \\
0 & Y \\
\end{bmatrix} \le\mu_{1}
\begin{bmatrix}
I & -I \\
-I & I \\
\end{bmatrix} +\mu_{2}
\begin{bmatrix}
I &0\\
0& I \\
\end{bmatrix}\end{equation}
then
\begin{equation*}
\begin{gathered}
\begin{aligned}
F(t,x,u,p,X) - F(t,y,u,q,-Y)  \ge &-\omega_{R}\Big(  \mu_{1}\left(  |x-y|^{2}+\rho(p,q)^{2} \right) + \mu_{2} + |p-q| \\
&+ |x-y|\left(\max(|p |, |q|) +1 \right)\Big)\\
\end{aligned}
\end{gathered}
\end{equation*}
for all $t\in[0,T], x,y\in\overline{\Omega}, u\in\mathbb{R}$, with $|u|\le R$, $p,q\in\mathbb{R}^{n}\backslash\{ 0\}$ and $\rho(p,q) = \min\left( \frac{|p-q|}{\min(|p|,|q|)},1 \right)$.}
\end{enumerate}

\subsection{Existence and uniqueness for the Geodesic Model}\label{sec:proofforgeo}

We now prove that there exists a unique solution for the PDE (\ref{eqn:geoel2}) resulting from the minimisation  of the functional for the Geodesic Model (\ref{eqn:geofunc2F}). 

\begin{rmk}
It is important to note that although the values of $c_{1}$ and $c_{2}$ depend on $u$, they are fixed when we solve the PDE for $u$ and therefore the probem is a local one and Theorem~\ref{thm:satothm} can be applied. Once we update $c_{1}$ and $c_{2}$, using the updated $u$, then we fix them again and apply Theorem~\ref{thm:satothm}. In practice, as we near convergence, we find $c_{1}$ and $c_{2}$ stabilise so we typically stop updating $c_{1}$ and $c_{2}$ once the change in both values is below a tolerance.
\end{rmk}

To apply the above theorem to the proposed model (\ref{eqn:geoel2}),
the key step will be to verify the nine conditions.
First, we multiply (\ref{eqn:geoel2}) by the factor $|\nabla u|_{\varepsilon_{2}}$, obtaining the nonlinear PDE
\begin{equation}\label{eqn:old11}
\begin{aligned}
-\mu|\nabla u|_{\varepsilon_{2}}\nabla\cdot \bigg( G(x,\nabla z) \frac{\nabla u}{|\nabla u|_{\varepsilon_{2}}} \bigg) + |\nabla u|_{\varepsilon_{2}}\bigg[ \lambda_{1} (z(x,y)-c_{1})^{2} -\lambda_{2} (z(x,y) - c_{2})^{2}\bigg.&\\
 \bigg. + \theta\mathcal{D}_{G}(x,y) + \alpha\nu'_{\varepsilon}(u) &\bigg]= 0
\end{aligned}
\end{equation}
where $G(x,\nabla z) = g(|\nabla z(x,y)|)$. We can rewrite this as
\begin{equation}\label{eqn:Fgeo}
F(x,u,p,X) = -\mu\, \text{trace}\left(  A(x,p)X \right)  - \mu\langle \nabla G(x,\nabla z), p \rangle + |p|k(u) +  |p|f(x)=0
\end{equation}
where $f(x) =  \lambda_{1}(z(x)-c_{1})^{2} - \lambda_{2}(z(x) - c_{2})^{2}$, $k(u)=\alpha\nu'_{\varepsilon}(u)$, $p  = (p_{1},p_{2}) = |\nabla u|_{\varepsilon_{2}}$, $X$ is the Hessian of $u$ and
\begin{equation}\label{eqn:defA}
A(x,p) =
\begin{bmatrix}
G(x,\nabla z)\frac{p_{2}^{2}}{|p|^{2}} & - G(x,\nabla z)\frac{p_{1}p_{2}}{|p|^{2}} \\
- G(x,\nabla z)\frac{p_{1}p_{2}}{|p|^{2}} & G(x,\nabla z)\frac{p_{1}^{2}}{|p|^{2}} \\
\end{bmatrix}
\end{equation}

\begin{thm}[Theory for the Geodesic Model]\label{thm:no2}
The parabolic PDE
$\frac{\partial u}{\partial t} + F(t,x,u,Du,D^{2}u) = 0 $ with $u_{0} = u(0,x)\in C(\overline{\Omega})$, $F$ as defined in (\ref{eqn:Fgeo}) and Neumann boundary conditions
has a unique solution $u=u(t,x)$ in $C([0,T)\times\overline{\Omega})$.
\end{thm}
{\bf Proof}. By Theorem \ref{thm:satothm}, it remains to verify that $F$ satisfies (C1)--(C2) and (I1)--(I7).  We will show that each of the conditions is satisfied. Most are simple to show, the exception being (I7) which is non-trivial.

{\underline{\emph{(C1)}:}} Equation (\ref{eqn:Fgeo}) only has dependence on $u$ in the term $k(u)$, we therefore have a restriction on the choice of $k$, requiring $k(v)\ge k(u)$ for $v\ge u$. This is satisfied for $k(u)=\alpha\nu'_{\varepsilon}(u)$ with $\nu'_{\varepsilon}(u)$ defined as in (\ref{eqn:nueps}).

{\underline{\emph{(C2)}:}} We find for arbitrary $s=(s_{1},s_{2})\in\mathbb{R}^{2}$ that $s^{T} A(x,p) s \ge 0$ and so $ A(x,p)\ge 0$. It follows that $-\text{trace}(A(x,p) X)  \le -\text{trace}(A(x,p) Y) $, therefore this condition is satisfied.

{\underline{\emph{(I1)}:}} $A(x,p)$ is only singular at $p=0$, however it is continuous elsewhere and
satisfies this condition.

{\underline{\emph{(I2)}:}} In $F$ the only term which depends on $u$ is $k(u)=\alpha\nu'_{\varepsilon}(u)$. With $\nu'_{\varepsilon}(u)$ defined as in (\ref{eqn:nueps}), in the limit $\varepsilon\rightarrow 0$ this function is a step function from $-2$ on $(\infty,0)$, $0$ on $[0,1]$ and $2$ on $(0,\infty)$. So we can choose any constant $\varepsilon < -2$. With $\varepsilon\neq 0$ there is smoothing at the end of the intervals, however there is still a lower bound on $L$ for $\nu'_{\varepsilon}(u)$ and we can choose any constant $\gamma < L$.

{\underline{\emph{(I3)}:}} $F$ is continuous at $(x, 0, 0)$ for any $x \in \Omega$ because $F^{*}(x, 0, 0) = F_{*}
(x, 0, 0) = 0.$ Hence this condition is satisfied.

{\underline{\emph{(I4)}:}} The Euler-Lagrange equations give Neumann boundary conditions
\[B(x,\nabla u) = \frac{\partial u}{\partial {\bm n}} ={\bm{n}}\cdot\nabla~u = \langle {\bm n} , \nabla u \rangle =  0
\]
 on $\partial\Omega$,
where ${\bm n}$ is the outward unit normal vector, and we see that $B(x,\nabla u)\in C^{1,1}\left(\mathbb{R}^{n}\times\mathbb{R}^{n}\backslash\{ 0 \}\right)$ and therefore this condition is satisfied.

{\underline{\emph{(I5)}:}} By the definition above, $B(x,\lambda\nabla u) = \langle {\bm n} , \lambda \nabla u \rangle =\lambda \langle {\bm n} ,  \nabla u \rangle = \lambda B(x,\nabla u) $. So this condition is satisfied.

{\underline{\emph{(I6)}:}} As before we can use the definition, $\langle \bm{n}(x), D_{p}B(x,p)\rangle = \langle \bm{n}(x), \bm{n}(x) \rangle = |\bm{n}(x)|^{2}$. So we can choose $\Theta = 1$ and the condition is satisfied.

{\underline{\emph{(I7)}:}} This is the most involved condition to prove and uses many other results. For clarity of the overall paper, we postpone the proof to Appendix A.
\hfill $\Box$

\subsection{Generalisation to other related models}

Theorems~\ref{thm:satothm} and \ref{thm:no2}
can be generalised to a few other models.
This amounts to writing each model as
a PDE of the form (\ref{eqn:Fgeo}) where $k(u)$ is monotone and $f(x), k(u)$ are bounded.
This is summarised in the following Corollary:
\begin{cor}\label{cor:uniquecorollary}
Assume that $c_1$ and  $c_2$ are fixed, with the terms
$f(x)$ and $k(u)$ respectively defined as follows for a few related models:
\begin{itemize}
\item{{\bf Chan-Vese \cite{SC10}}: $f(x) = f_{CV}(x):=\lambda_{1} (z(x)-c_{1})^{2} - \lambda_{2}(z(x)-c_{2})^{2}$, $k(u) = 0$}.
\item{{\bf Chan-Vese (Convex) \cite{SC8}}: $f(x) = f_{CV}(x)$, $k(u) = \alpha\nu'_{\varepsilon}(u)$}.
\item{{\bf Geodesic Active Contours \cite{SC5} and Gout et al. \cite{SC13}}: $f(x) = 0$, $k(u) = 0$}.
\item{{\bf Nguyen et al. \cite{SC25}}: $f(x) = \alpha\left( P_{B}(x,y) - P_{F}(x,y)\right) + \left( 1-\alpha\right) \left(1-2P(x,y)\right)$, $k(u) = 0$}.
\item{{\bf Spencer-Chen (Convex) \cite{SC}}: $f(x) = f_{CV}(x) + \theta\mathcal{D}_{E}(x)$, $k(u) = \alpha\nu'_{\varepsilon}(u)$}.
\end{itemize}
Then if we define a PDE of the general form
\[
-\mu \nabla\cdot\left( G(x)\frac{\nabla u}{|\nabla u|_{\varepsilon_{2}}} \right) + k(u)  + f(x)= 0
\]
with
\begin{enumerate}
\item[(i)]{Neumann boundary conditions $\frac{\partial u}{\partial {\bm n}} = 0$ (${\bm n}$ the outward normal unit vector)}
\item[(ii)]{$k(u)$ satisfies $k(u)\ge k(v)$ if $u\ge v$}
\item[(iii)]{$k(u)$ and $f(x)$ are bounded; and}
\item[(iv)]{$G(x) = Id$ or $G(x) = f(|\nabla z(x)|) = \frac{1}{1+|\nabla z(x)|^{2}}$,}
\end{enumerate}
we have a unique solution $u\in C([0,T)\times\overline{\Omega})$ for a given initialisation.
Consequently we conclude that all above models admit a unique solution.
\end{cor}

{\bf Proof}. The conditions (i)--(iv) are hold for all of these models. All of these models require Neumann boundary conditions and use the permitted G(x). The monotonicity of $\nu'_{\varepsilon}(u)$ is discussed in the proof of (C1) for Theorem~\ref{thm:no2} and the boundedness of $f(x)$ and $k(u)$ is clear in all cases.
\hfill $\Box$


\begin{rmk}
Theorem~\ref{thm:no2} and Corollary~\ref{cor:uniquecorollary} also generalise to cases where $G(x) = \frac{1}{1+\beta|\nabla z|^{2}}$ and to $G(x) = \mathcal{D}(x) g(|\nabla z|)$ where $\mathcal{D}(x)$ is a distance function such as in \cite{RC15,gout2005segmentation,RC16,SC}. The proof is very similar to that shown in \S\ref{sec:proofforgeo}, relying on Lipschitz continuity of the function $G(x)$.
\end{rmk}

\begin{rmk}
We cannot apply the classical viscosity solution framework to the Rada-Chen model \cite{SC28} as this is a non-local problem with $k(u) = 2\nu \left( \int_{\Omega}H_{\varepsilon}(u)\,\mathrm{d}\Omega-A_{1} \right)$.
\end{rmk}



\section{Numerical Results}
\label{sec:TC}

In this section we will demonstrate the advantages of the Geodesic Model for selective image segmentation over related and previous models. Specifically we shall compare
\begin{itemize}
\item{\bf M1 } --- the Nguyen et al. (2012) model \cite{SC25};
\item{\bf M2 } --- the Rada-Chen (2013) model \cite{SC28};
\item{\bf M3 } --- the convex Spencer-Chen  (2015) model \cite{SC};
\item{\bf M4 } --- the convex Liu et al. (2017) model \cite{Liu2017};
\item{\bf M5 } --- the reformulated Rada-Chen model with geodesic distance penalty (see Remark~\ref{rmk:M4});
\item{\bf M6 } --- the reformulated Liu et al. model with geodesic distance penalty (see Remark~\ref{rmk:M4});
\item{\bf M7 } --- the proposed convex Geodesic Model (Algorithm \ref{alg:AOS}).
\end{itemize}

\begin{rmk}[A note on {\bf{M5}} and {\bf{M6}}] \label{rmk:M4}\ We include {\bf M5 -- M6} to test how the geodesic distance penalty term can improve {\bf M2} \cite{SC28} and {\bf M4} \cite{Liu2017}. These were obtained as follows:
\begin{itemize}
\item{we extend {\bf{M2}} to {\bf{M5}} simply by including the geodesic distance function $\mathcal{D}_{G}(x,u)$ in the functional.}
\item{we extend {\bf{M4}} to {\bf{M6}} with a minor reformulation to include data fitting terms.
Specifically, the model {\bf{M6}} is
\begin{equation}\label{eqn:convliucvfunc}
\begin{split}
\min_{u,c_{1},c_{2}}\Big\{F_{CV\omega}(u,c_{1},c_{2})=& \int_{\Omega}\omega^{2}(x,y)\left[   \lambda_{1}(z(x,y)-c_{1})^{2}  - \lambda_{2}(z(x,y)-c_{2})^{2}  \right]  u\,\mathrm{d}\Omega \\
  +&  \mu\int_{\Omega}g(|\nabla z|))|\nabla u|\,\mathrm{d}\Omega + \theta\int_{\Omega}\mathcal{D}_{G}(x,y)u \,\mathrm{d}\Omega + \alpha\int_{\Omega}\nu_{\varepsilon}(u)\,\mathrm{d}\Omega\Big\}
\end{split}
\end{equation}
for $\mu,\lambda_{1},\lambda_{2}$ non-negative fixed parameters, $\alpha$ and $\nu_{\varepsilon}(u)$ as defined in
(\ref{eqn:nueps}) and $\omega$ as defined for the convex Liu et al. model. This is a convex model and is the same as the proposed Geodesic Model {\bf{M7}} but with weighted intensity fitting terms.}
\end{itemize}
\end{rmk}

\begin{figure}[htb!]
\centering
\makebox[\textwidth][c]{
\includegraphics[width=0.25\textwidth,height = 3cm]{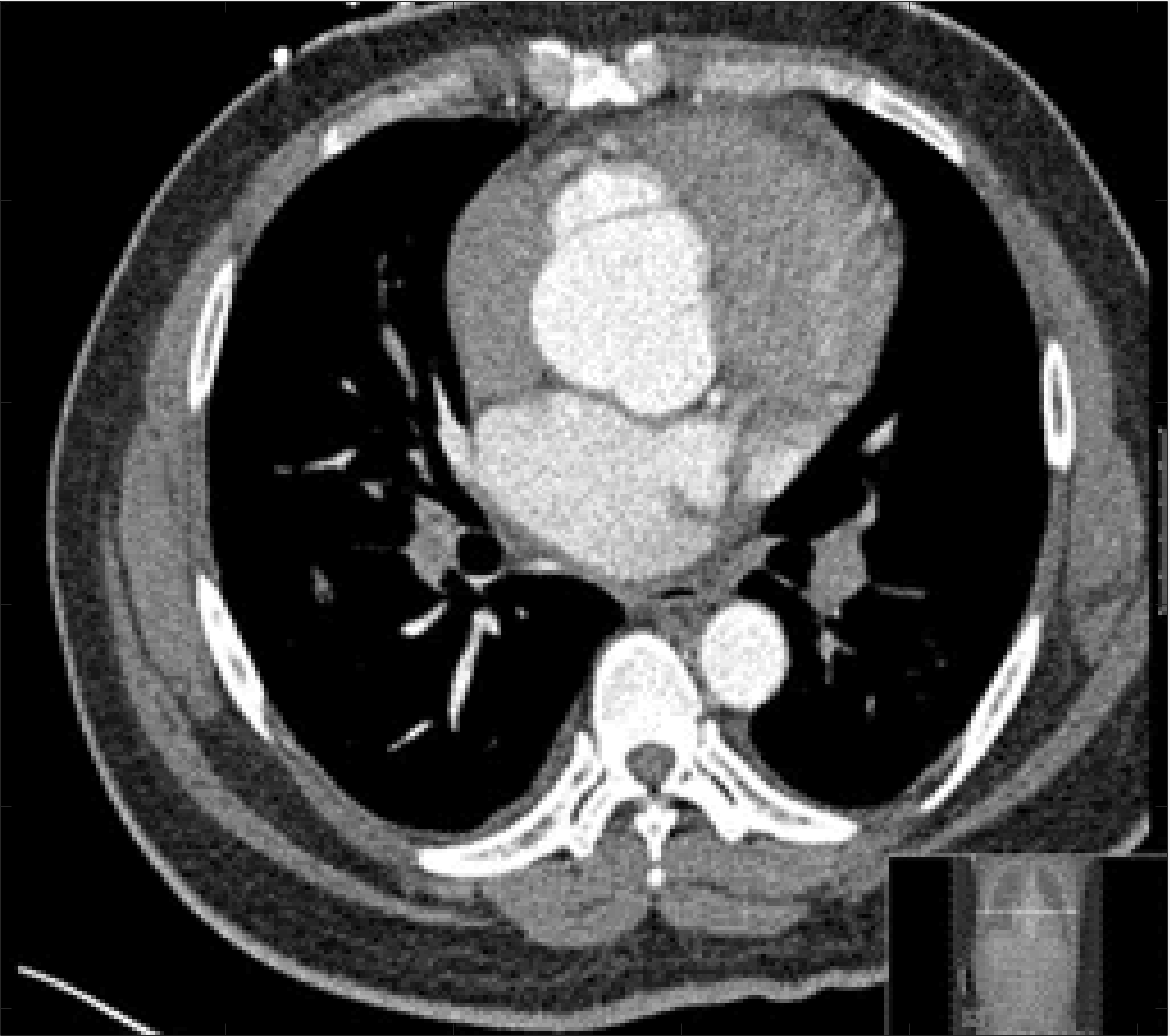}
\quad
\includegraphics[width=0.25\textwidth,height = 3cm]{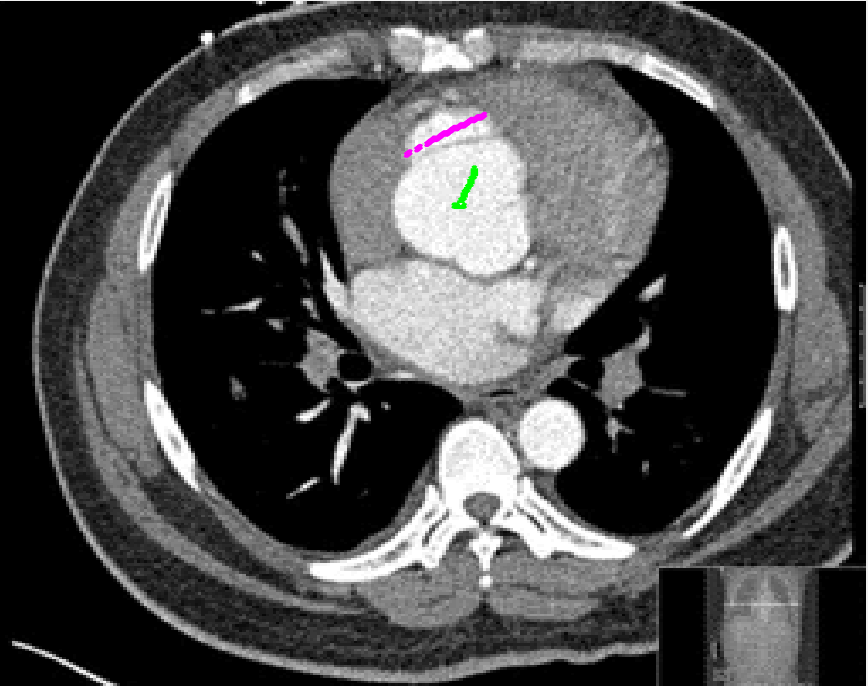}
\quad
\includegraphics[width=0.25\textwidth,height = 3cm]{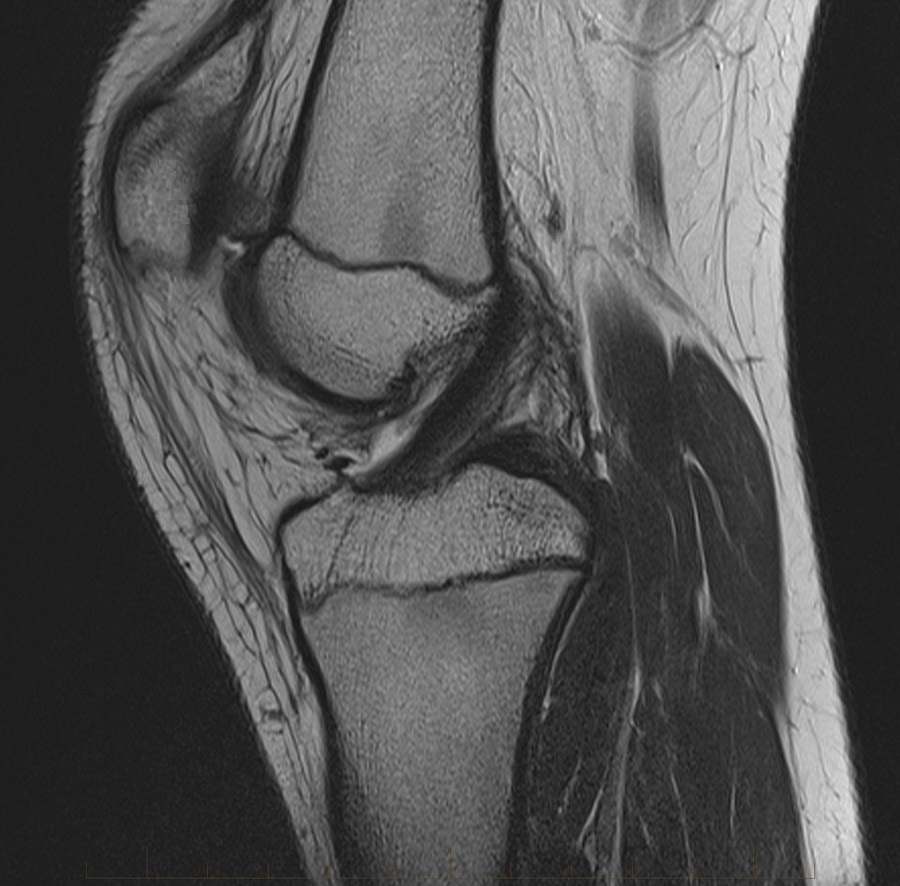}
\quad
\includegraphics[width=0.25\textwidth,height = 3cm]{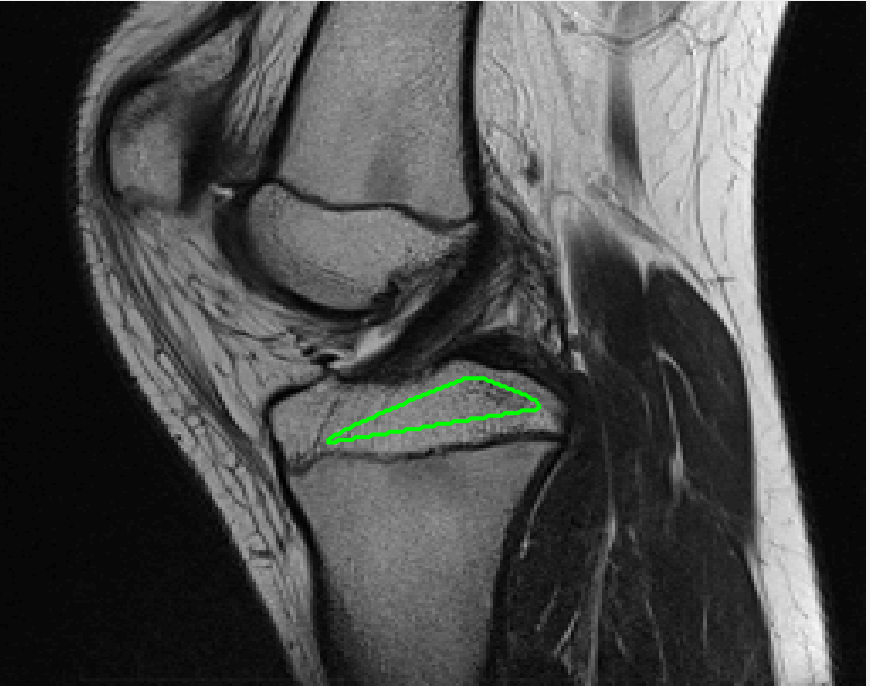}
}\\
\centerline{
(i)\hspace*{0.25\textwidth}(ii)\hspace*{0.25\textwidth}(iii)\hspace*{0.25\textwidth}(iv)}
\caption{Test 1 setting:\ (i) Image 1; \
(ii) Image 1 with marker and anti-marker set shown in green and pink respectively; \ (iii) Test Image 2;\
(iv) Image 2 with marker set shown.\label{fig:testimages}}
\end{figure}

Four sets of test results are shown below.
In Test 1 we compare models {\bf M1 -- M6}  to the proposed model {\bf M7} for two images which are hard to segment. The first is a CT scan from which we would like to segment the lower portion of the heart, the second is an MRI scan of a knee and we would like to segment the top of the Tibia. See Figure~\ref{fig:testimages} for the test images and the marker sets used in the experiments.
In Test 2 we will review the sensitivity of the proposed model to the main parameters.
In Test 3 we will give several results achieved by the model using marker and anti-marker sets. In Test 4 we show the initialisation independence and marker independence of the Geodesic Model on real images.

For {\bf M7}, we denote by $\tilde{u}$   the thresholded $u > \tilde{\gamma}$ at some value $\tilde{\gamma}\in(0,1)$ to define the segmented region. Although the threshold can be chosen arbitrarily in $(0,1)$ from the work by \cite[Thm 1]{SC8} and \cite{SC},
we usually take $\tilde{\gamma}=0.5$.

{\underline{\it Quantitative Comparisons.}}
To measure   the quality of a segmentation, we use the Tanimoto Coefficient (TC) (or Jaccard Coefficient \cite{NPH:NPH37}) defined by
\[
TC(\tilde{u},GT) = \frac{|\tilde{u}\cap GT|}{|\tilde{u}\cup GT|}
\]
where 
 GT is the `ground truth' segmentation and
$\tilde{u}$ is the result from a particular model.  
 This measure takes value one for a segmentation which coincides perfectly with the ground truth and reduces to zero as the quality of the segmentation gets worse. In the other tests, where a ground truth is not available, we use visual plots.

{\underline{\it Parameter Choices and Implementation.}} We set $\mu = 1$, $\tau = 10^{-2}$ and vary $\lambda = \lambda_{1} = \lambda_{2}$ and $\theta$. Following \cite{SC8} we let $\alpha = ||\lambda_{1}(z-c_{1})^{2}-\lambda_{2}(z-c_{2})^{2}+\theta\mathcal{D}_{G}(x,y)||_{L^{\infty}}$. To implement the marker points in MATLAB we use \texttt{roipoly} for choosing a small number of points by clicking and also \texttt{freedraw} which allows the user to draw a path of marker points. The stopping criteria used is the dynamic residual falling below a given threshold, i.e. once $||u^{k+1}-u^{k}||/||u^{k}|| < tol$ the iterations stop (we use $tol = 10^{-6}$ in the tests shown).

{\bf Test 1 -- Comparison of models {\bf M1 -- M7}}.

\begin{figure}[h!]
\centering
\makebox[\textwidth][c]{
  \begin{subfigure}[b]{1.1\textwidth}
        \includegraphics[width=0.33\textwidth,height = 4cm]{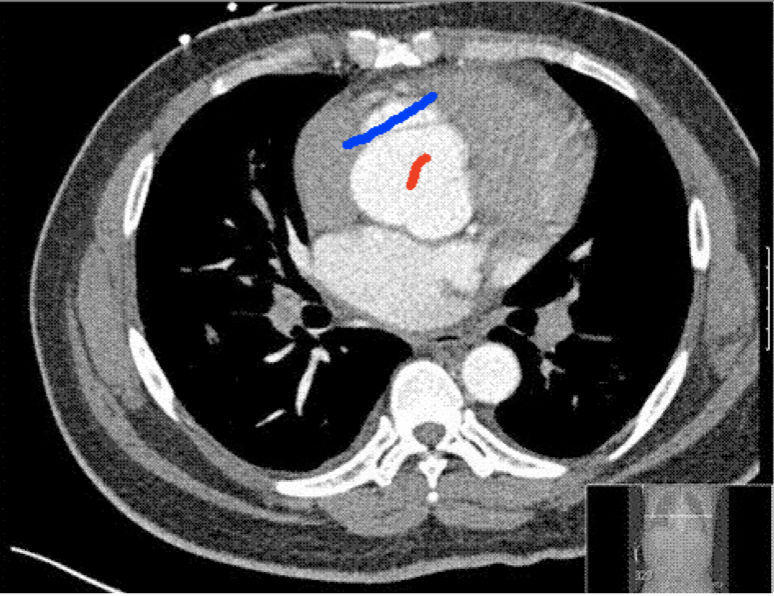}\quad
        \includegraphics[width=0.33\textwidth,height = 4cm]{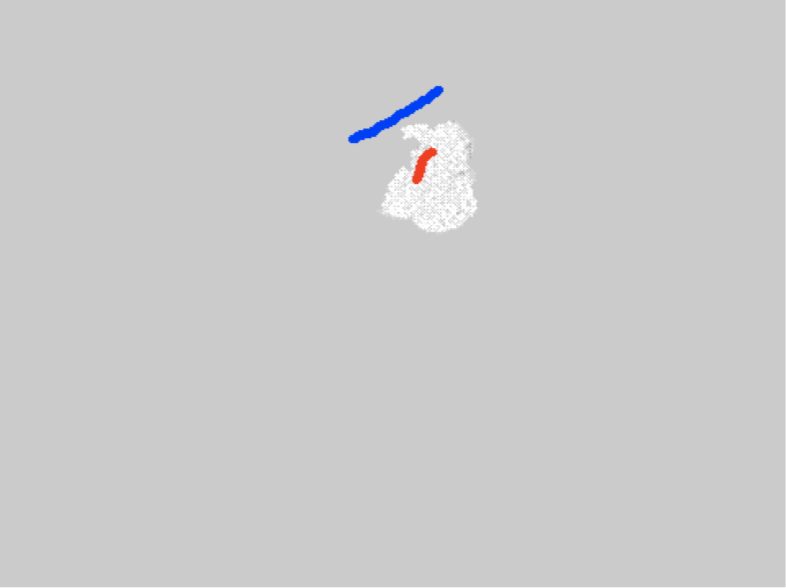}\quad
        \includegraphics[width=0.33\textwidth,height = 4cm]{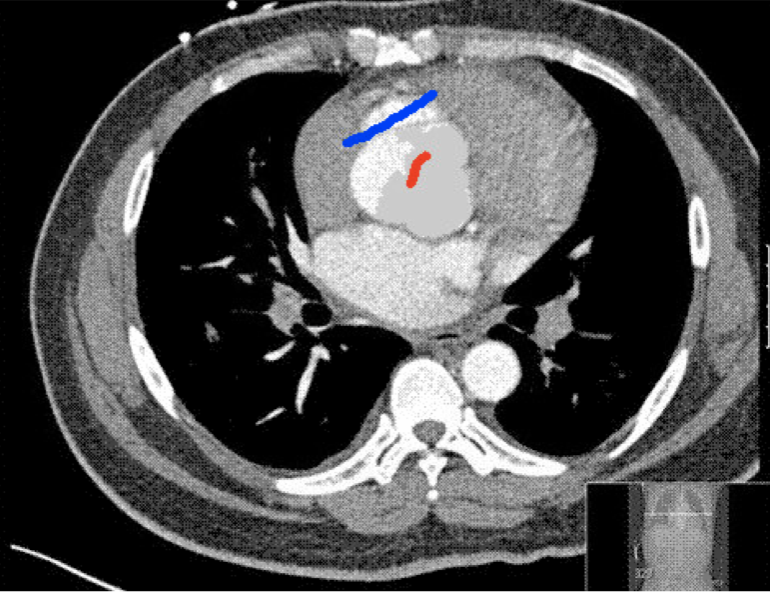}
        \caption{{\bf M1} (Left to right:) Test Image 1 with markers (red) and anti-markers (blue), foreground segmentation and background segmentation (we used published software, no parameter choice required).}
        \end{subfigure}
    }%

\vspace{0.2in}

\centering
\makebox[\textwidth][c]{
  \begin{subfigure}[b]{0.28\textwidth}
        \includegraphics[width=\textwidth]{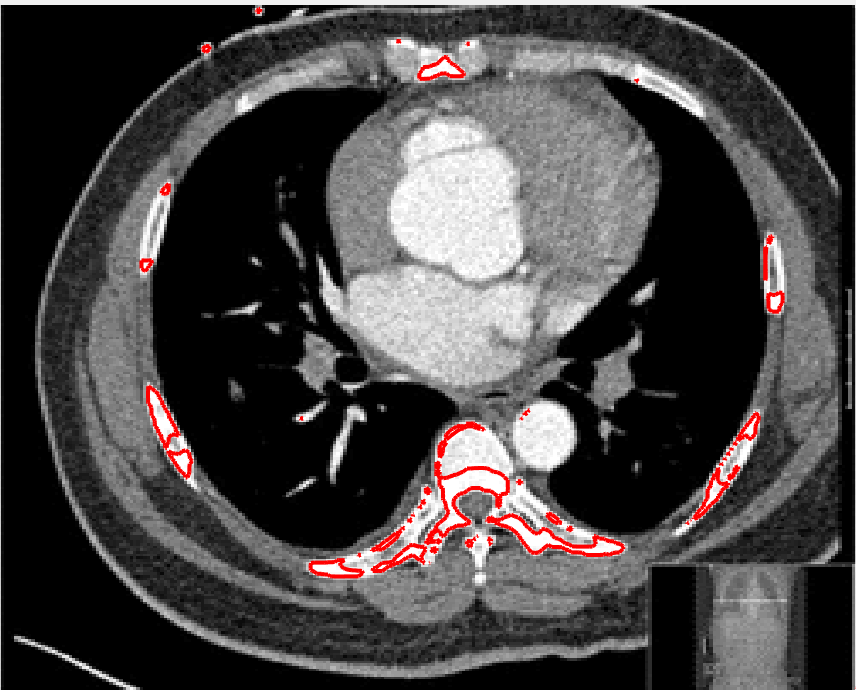}
        \caption{{\bf M2} $\lambda=1$, $\gamma=10$.}
    \end{subfigure}\quad%
      \begin{subfigure}[b]{0.28\textwidth}
        \includegraphics[width=\textwidth]{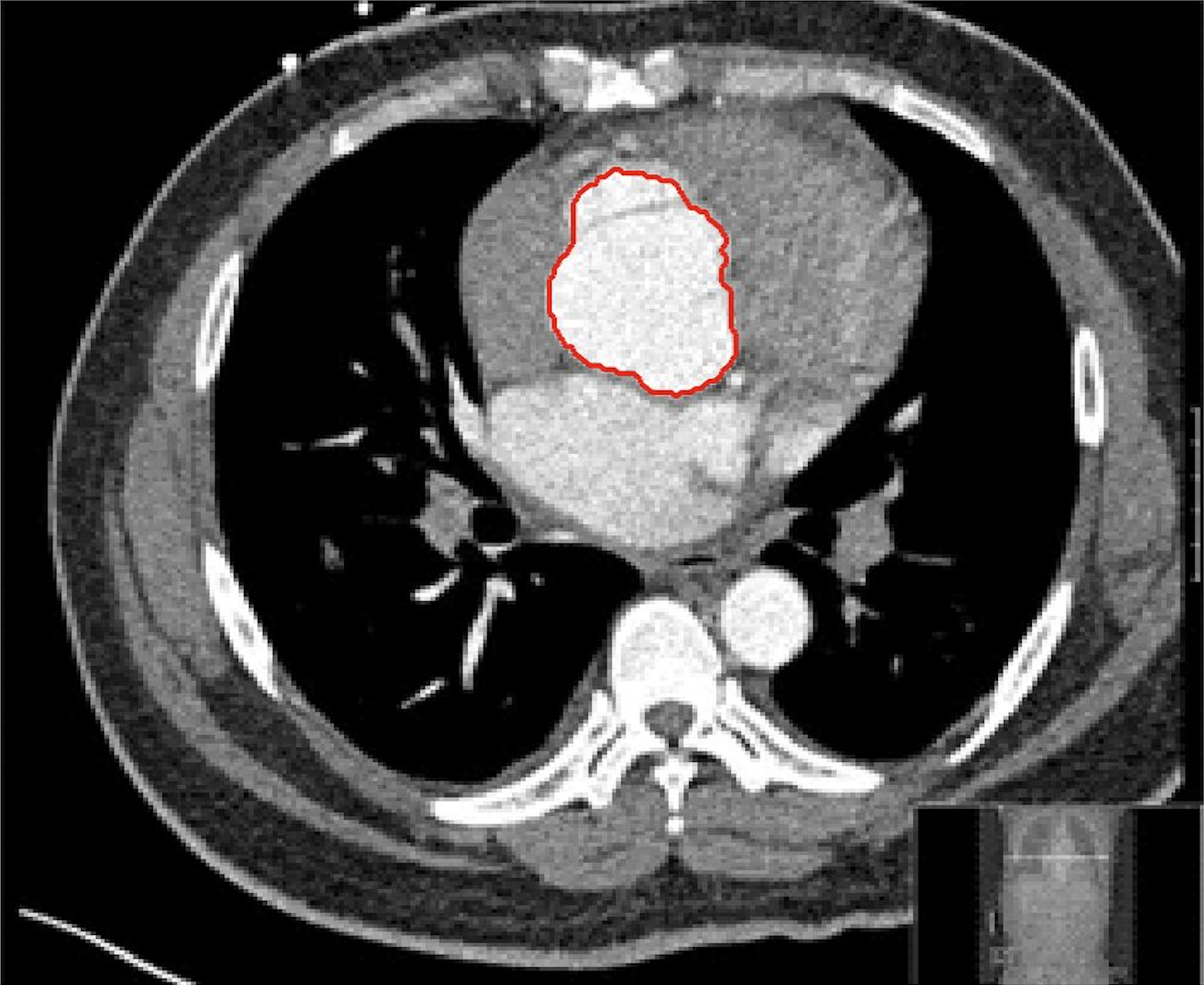}
        \caption{{\bf M3} $\lambda = 5$, $\theta = 3$.}
    \end{subfigure}\quad%
      \begin{subfigure}[b]{0.28\textwidth}
        \includegraphics[width=\textwidth]{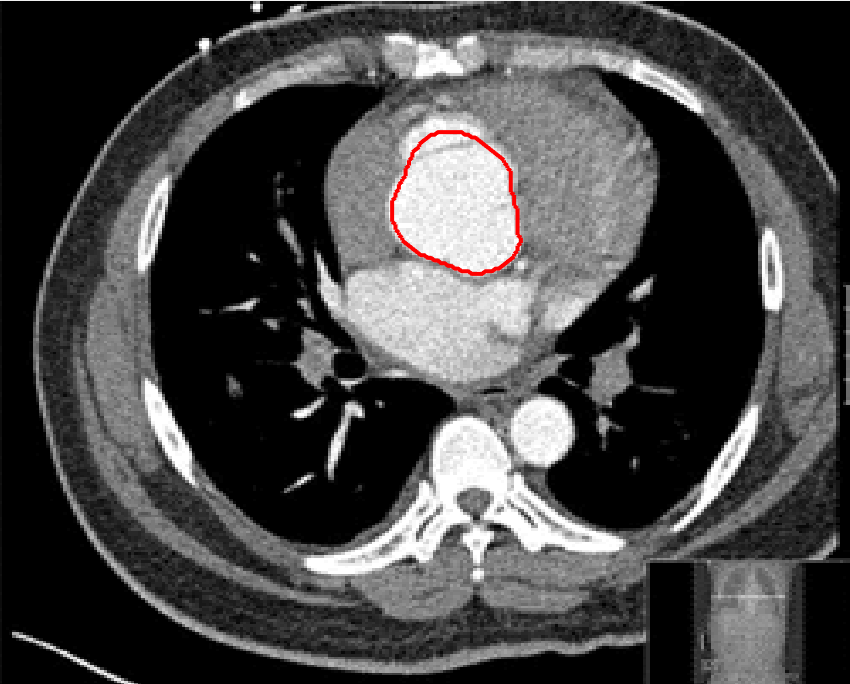}
        \caption{{\bf M4} $\lambda = 1/4$.}
    \end{subfigure}
    }\\
 \makebox[\textwidth][c]{
      \begin{subfigure}[b]{0.28\textwidth}
        \includegraphics[width=\textwidth]{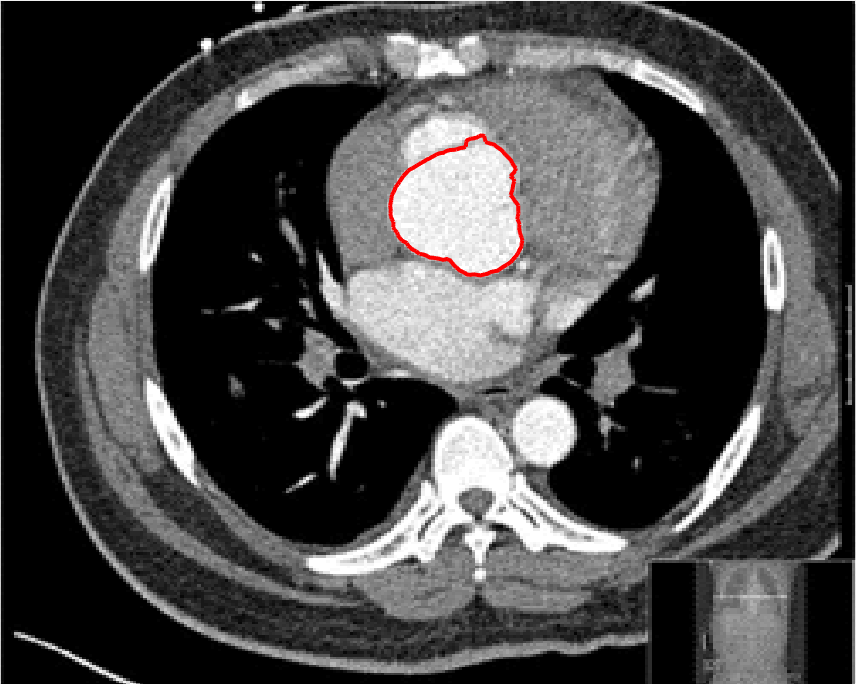}
        \caption{{\bf M5} $\lambda = 5, \gamma = 3, \theta = \frac{1}{10}$.}
    \end{subfigure}\quad%
      \begin{subfigure}[b]{0.28\textwidth}
        \includegraphics[width=\textwidth]{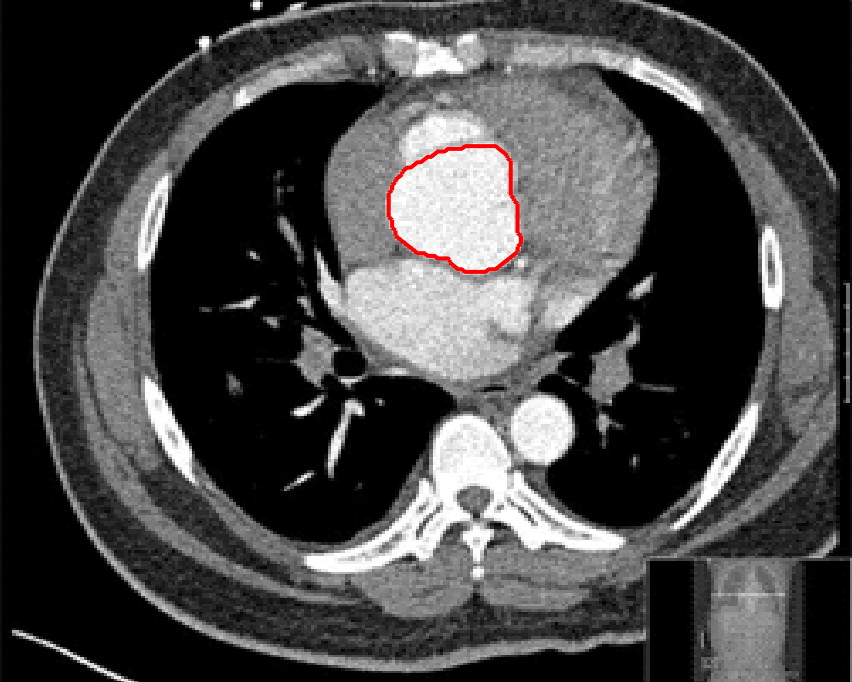}
         \caption{{\bf M6} $\lambda = 15, \theta = 3$.}
    \end{subfigure}\quad%
          \begin{subfigure}[b]{0.28\textwidth}
        \includegraphics[width=\textwidth]{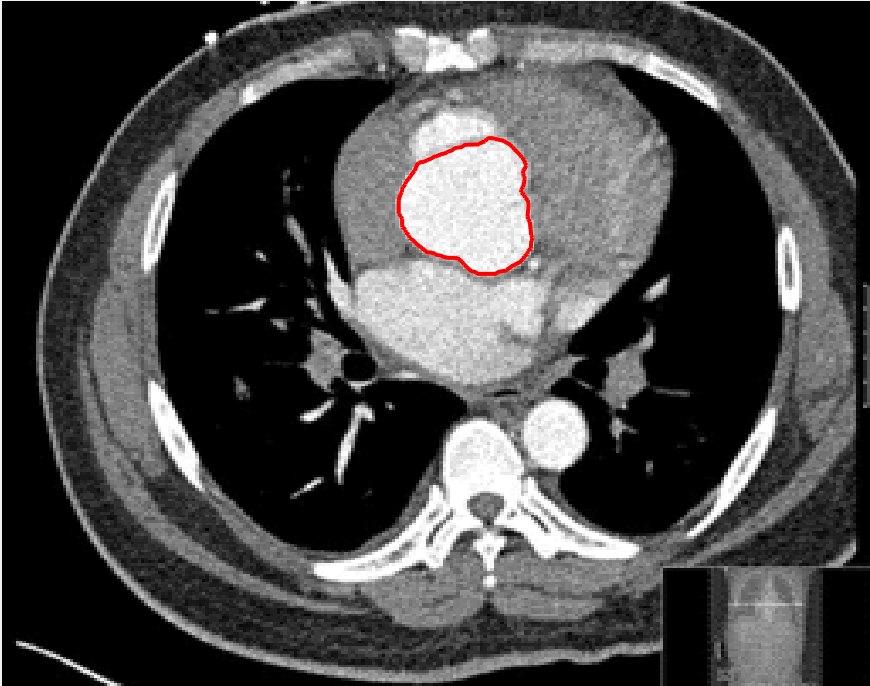}
        \caption{{\bf M7} $\lambda = 10, \theta = 1$.}
    \end{subfigure}
    }%
\caption{Visual comparison of {\bf M1 -- M7}  results for Test Image 1\label{fig:test1}. {\bf M1} segmented part of the object, {\bf M2 -- M4} failed to segment the object, {\bf M5} gave a reasonable result (though not accurate) and, {\bf M6} and {\bf M7} correctly segmented the object.  }
\end{figure}

\begin{figure}[htb!]
\centering
\makebox[\textwidth][c]{
  \begin{subfigure}[b]{1.1\textwidth}
        \includegraphics[width=0.33\textwidth,height = 4cm]{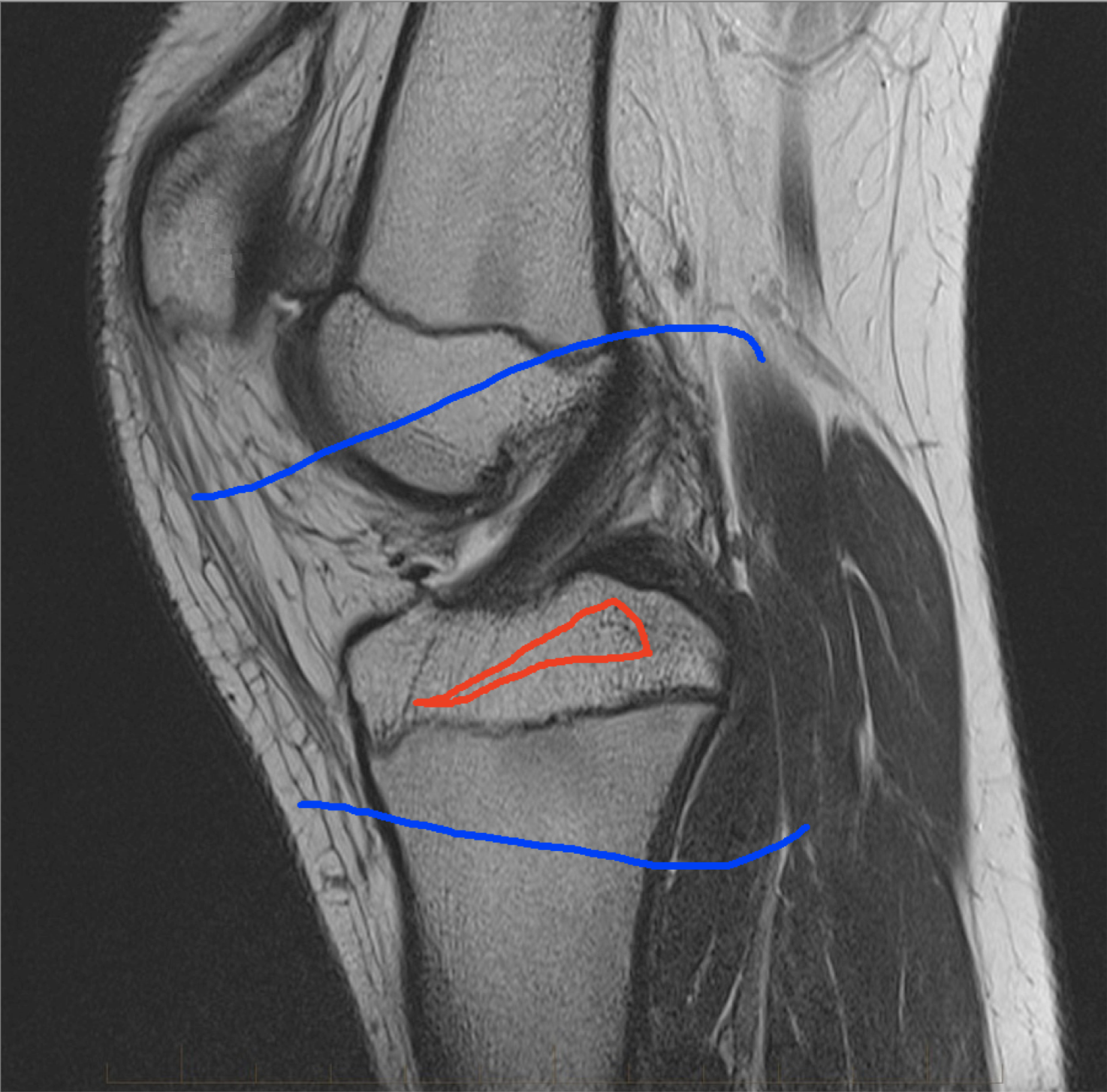}\quad
        \includegraphics[width=0.33\textwidth,height = 4cm]{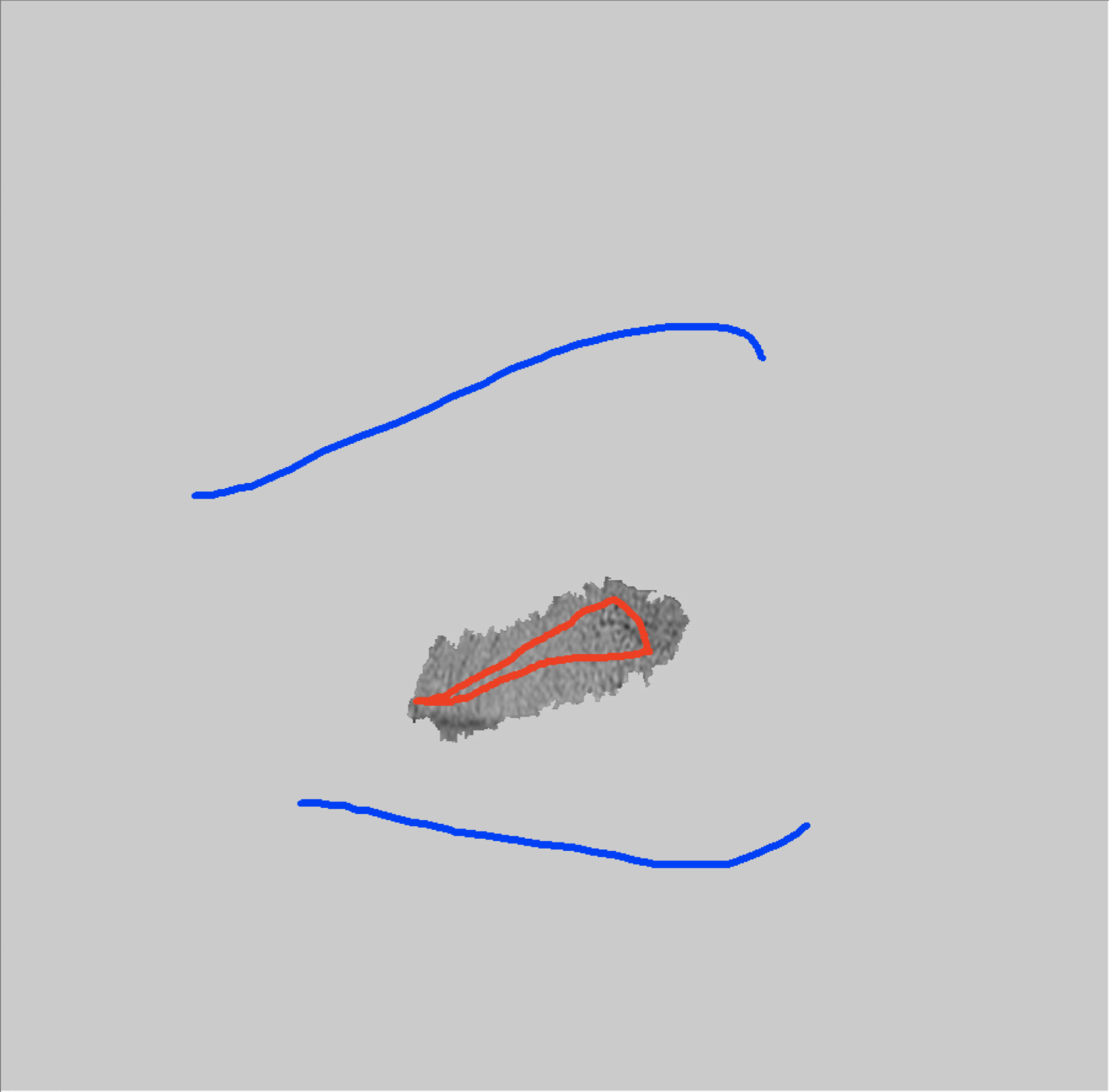}\quad
        \includegraphics[width=0.33\textwidth,height = 4cm]{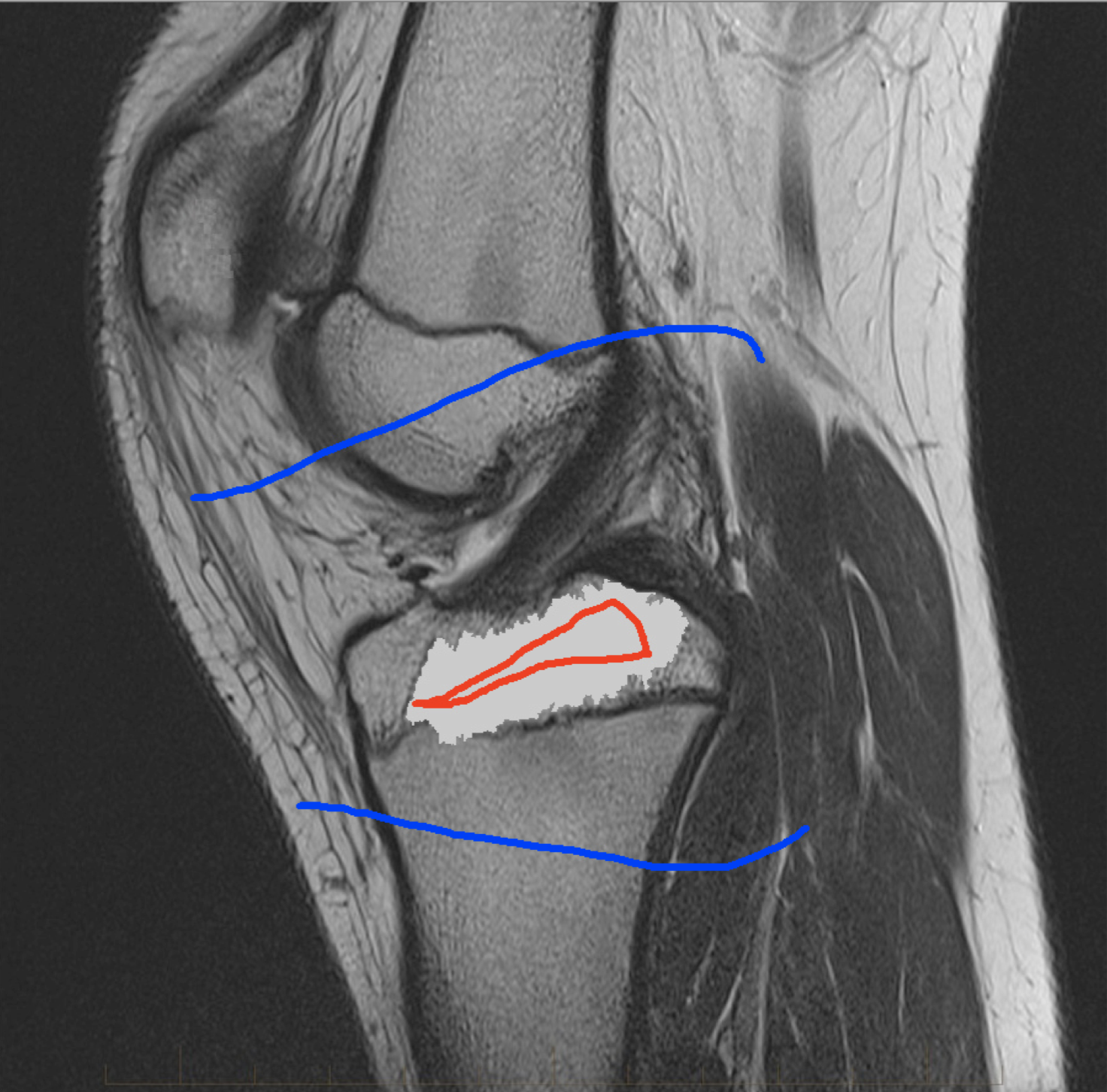}
        \caption{{\bf M1} (Left to right:) Test Image 2 with markers (red) and anti-markers (blue), foreground segmentation and background segmentation (we used published software, no parameter choice required).}
        \end{subfigure}
    }%

\vspace{0.2in}

\centering
\makebox[\textwidth][c]{
  \begin{subfigure}[b]{0.28\textwidth}
        \includegraphics[width=\textwidth]{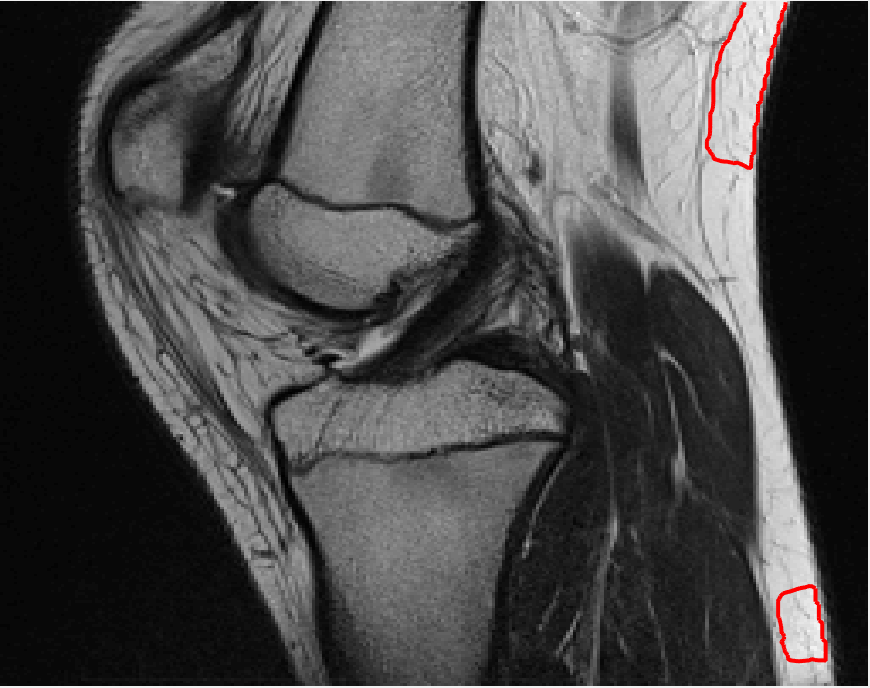}
        \caption{{\bf M2} $\lambda=1$, $\gamma=15$.}
    \end{subfigure}\quad%
      \begin{subfigure}[b]{0.28\textwidth}
        \includegraphics[width=\textwidth]{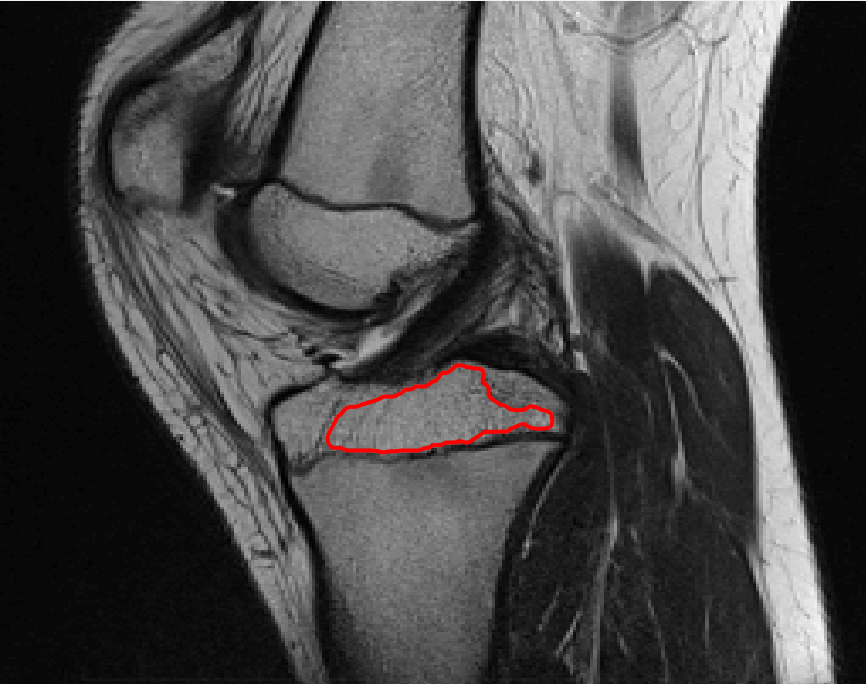}
        \caption{{\bf M3} $\lambda = 5$, $\theta = 1$.}
    \end{subfigure}\quad%
      \begin{subfigure}[b]{0.28\textwidth}
        \includegraphics[width=\textwidth]{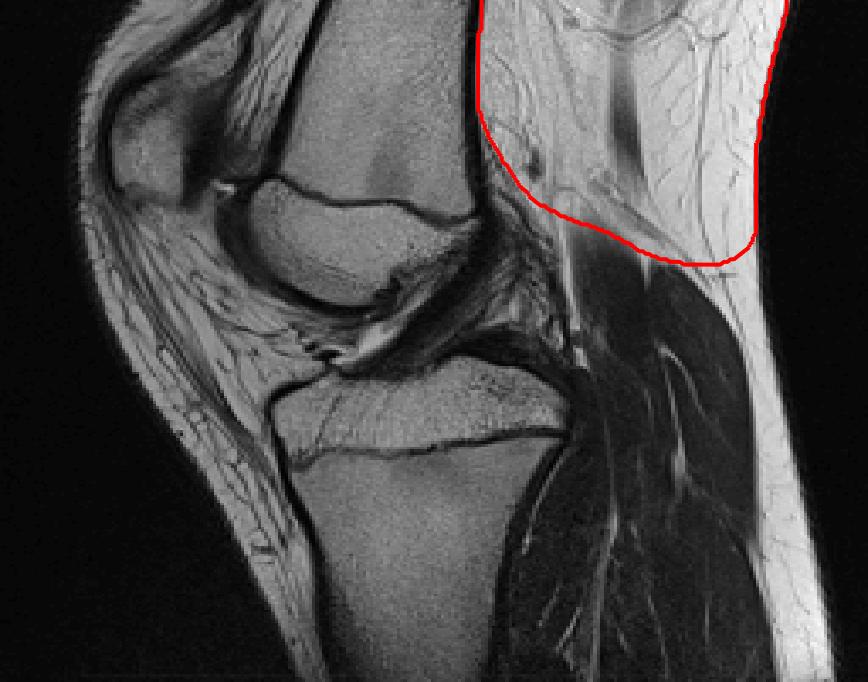}
        \caption{{\bf M4} $\lambda = 1/8$.}
    \end{subfigure}
    }\\
    \makebox[\textwidth][c]{
      \begin{subfigure}[b]{0.28\textwidth}
        \includegraphics[width=\textwidth]{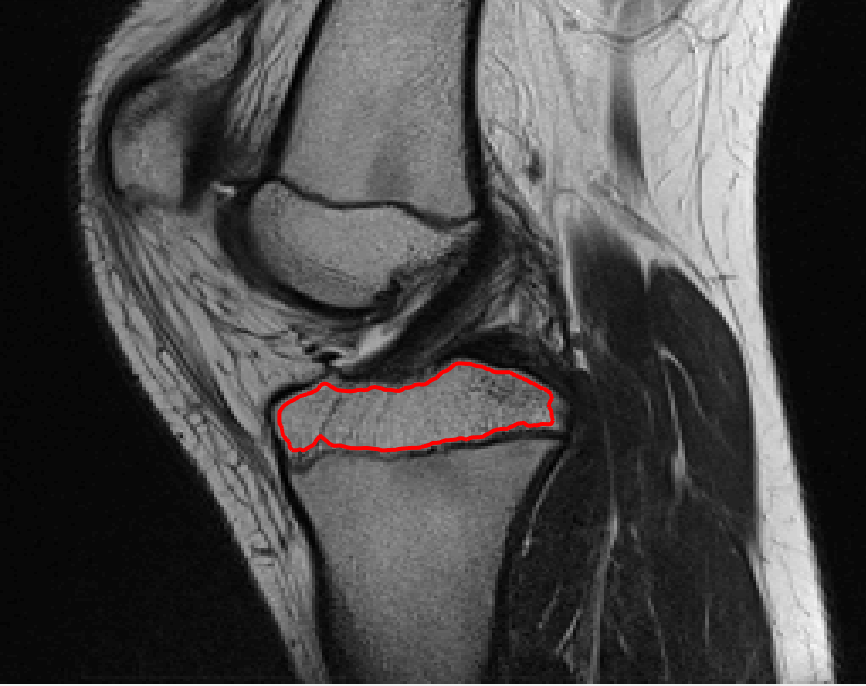}
        \caption{{\bf M5} $\lambda = 1, \gamma = 15, \theta = \frac{1}{10}$.}
    \end{subfigure}\quad%
      \begin{subfigure}[b]{0.28\textwidth}
        \includegraphics[width=\textwidth]{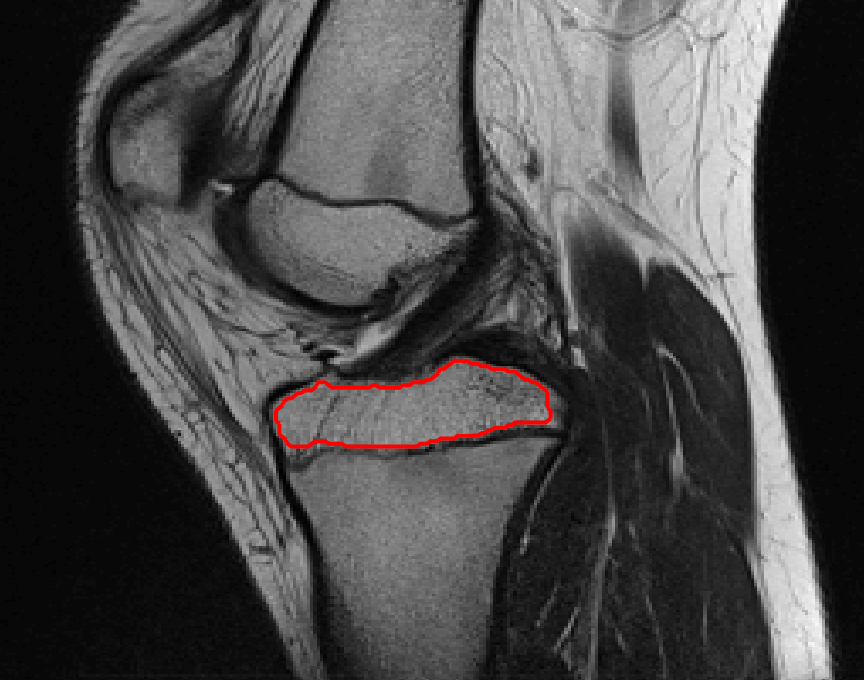}
        \caption{{\bf M6} $\lambda = 15, \theta = 1$.}
    \end{subfigure}\quad%
          \begin{subfigure}[b]{0.28\textwidth}
        \includegraphics[width=\textwidth]{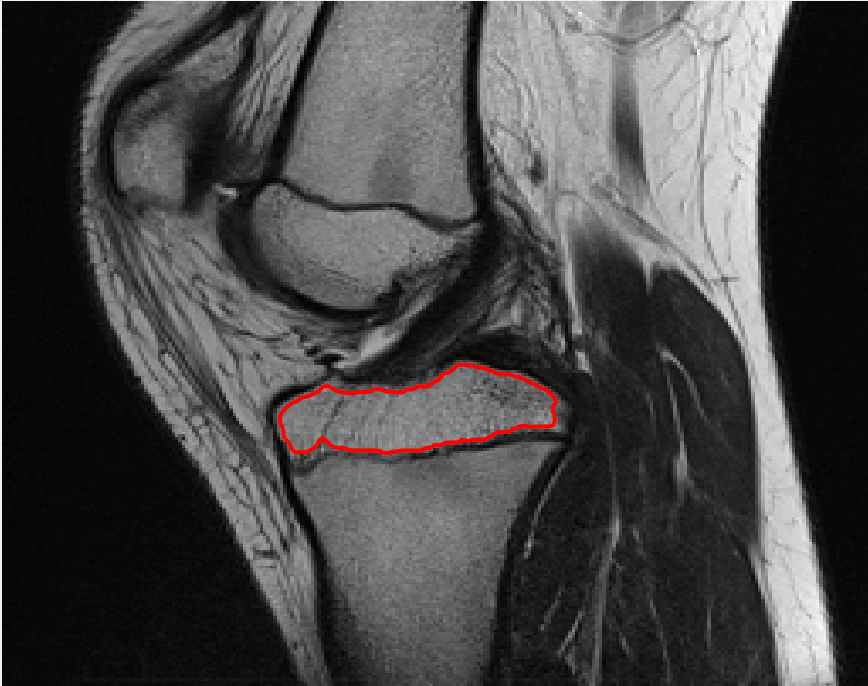}
        \caption{{\bf M7} $\lambda = 10, \theta = 1$.}
    \end{subfigure}
    }\\
\caption{Visual comparison of {\bf M1 -- M7}  results for Test Image 2.  {\bf M1} segmented part of the object, {\bf M2 -- M4} failed to segment the object, {\bf M5}, {\bf M6} and {\bf M7} correctly segmented the object.}\label{fig:test2}
        \end{figure}

\begin{figure}[htb!]
\centering
\makebox[\textwidth][c]{
  \begin{subfigure}[b]{0.35\textwidth}
        \includegraphics[width=\textwidth]{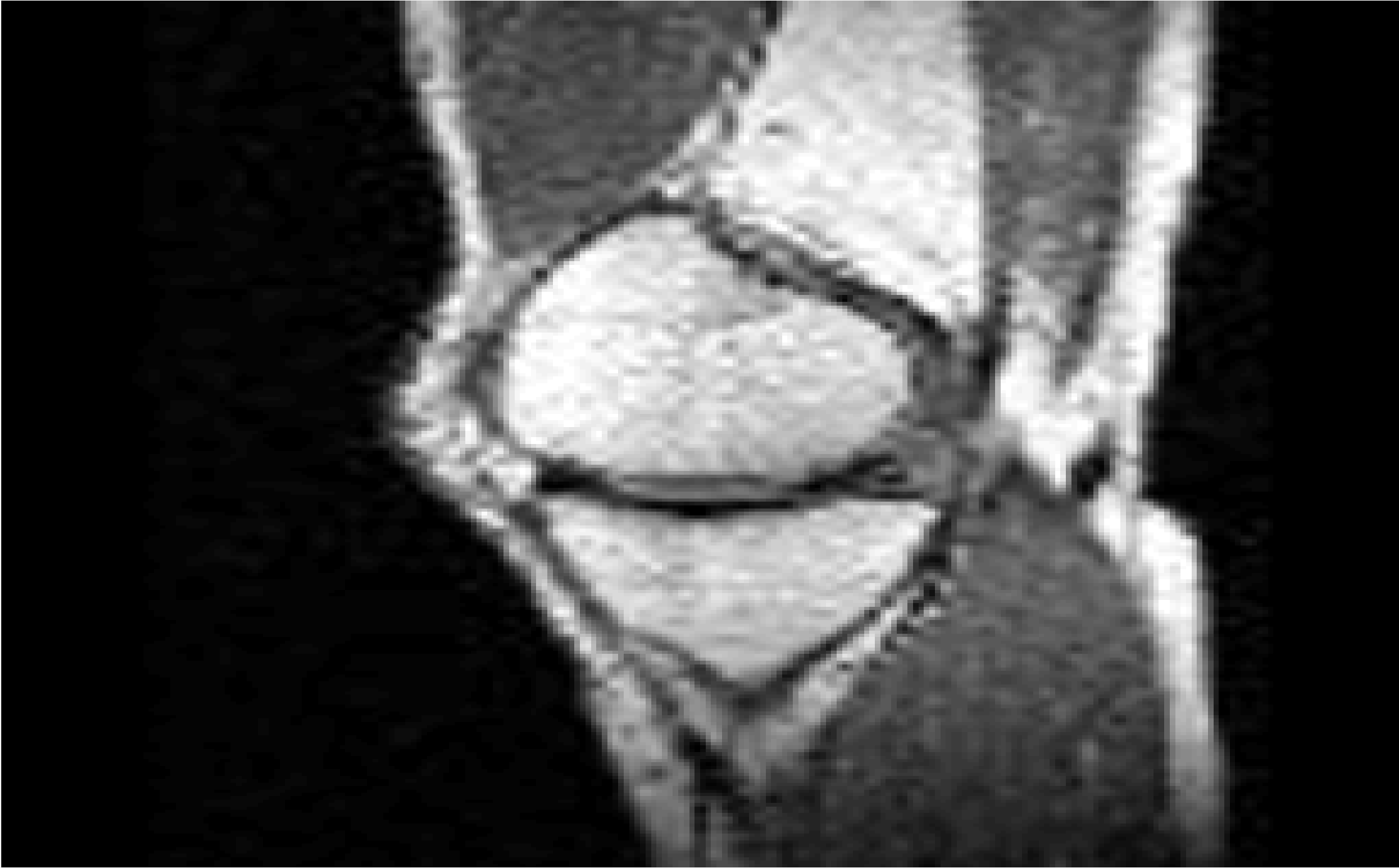}
        \caption{Original Image.}
    \end{subfigure}\quad%
      \begin{subfigure}[b]{0.35\textwidth}
        \includegraphics[width=\textwidth]{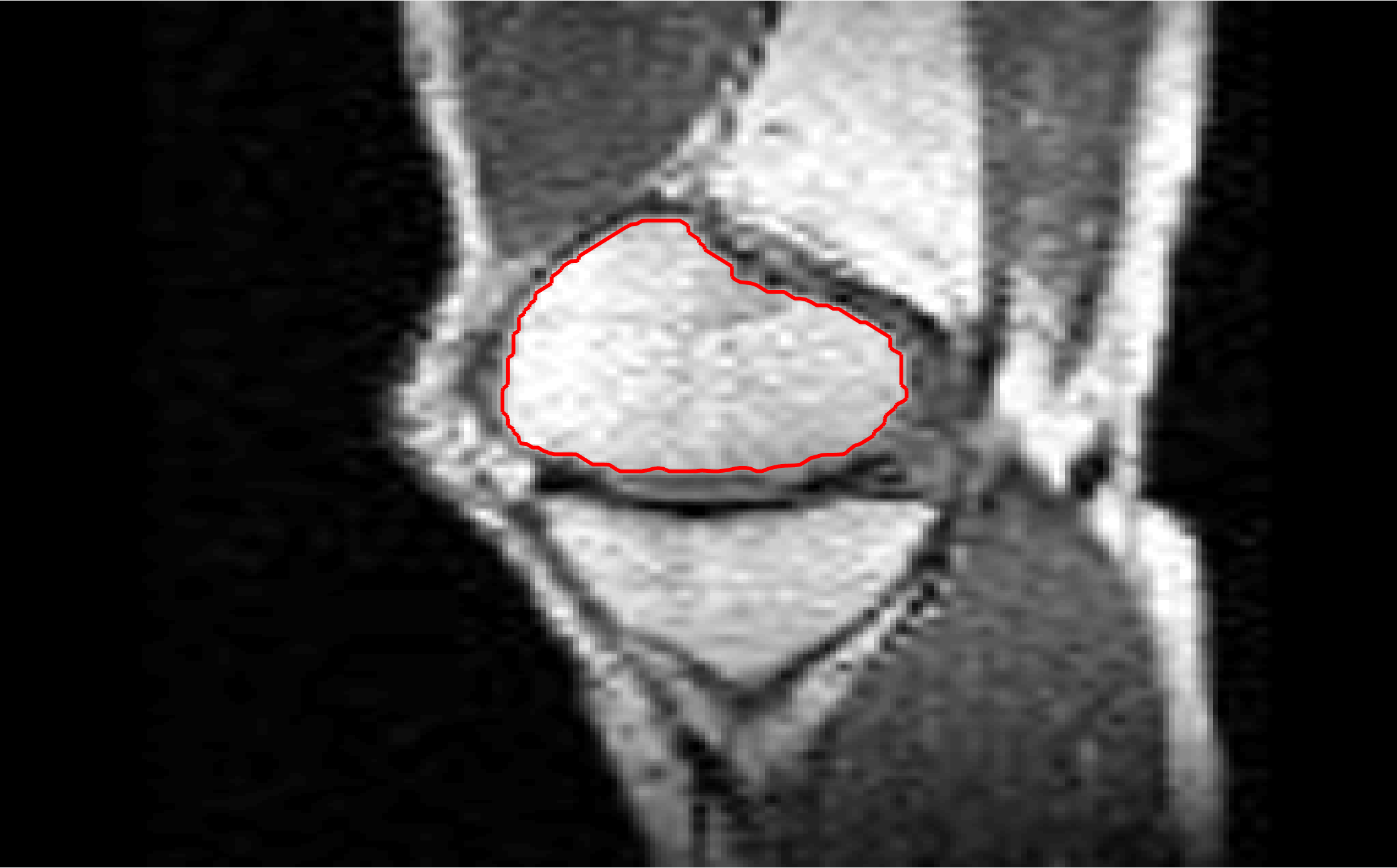}
        \caption{Ground Truth Segmentation.}
    \end{subfigure}
   }\\
   \makebox[\textwidth][c]{
         \begin{subfigure}[b]{0.39\textwidth}
        \includegraphics[width=\textwidth]{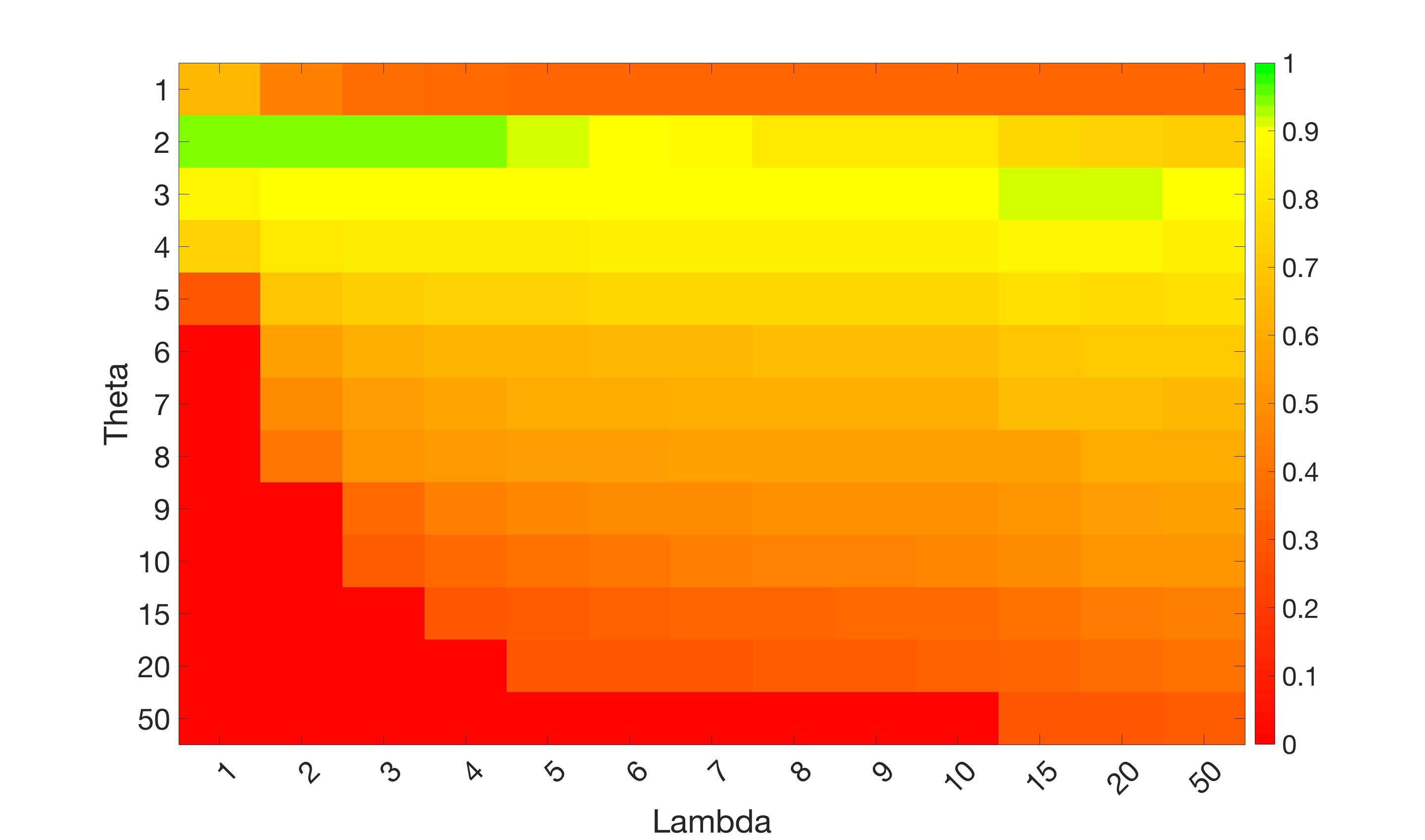}
        \caption{{\bf M3} TC values for various $\lambda$ and $\theta$ values.}
          \end{subfigure}\quad%
      \begin{subfigure}[b]{0.39\textwidth}
        \includegraphics[width=\textwidth]{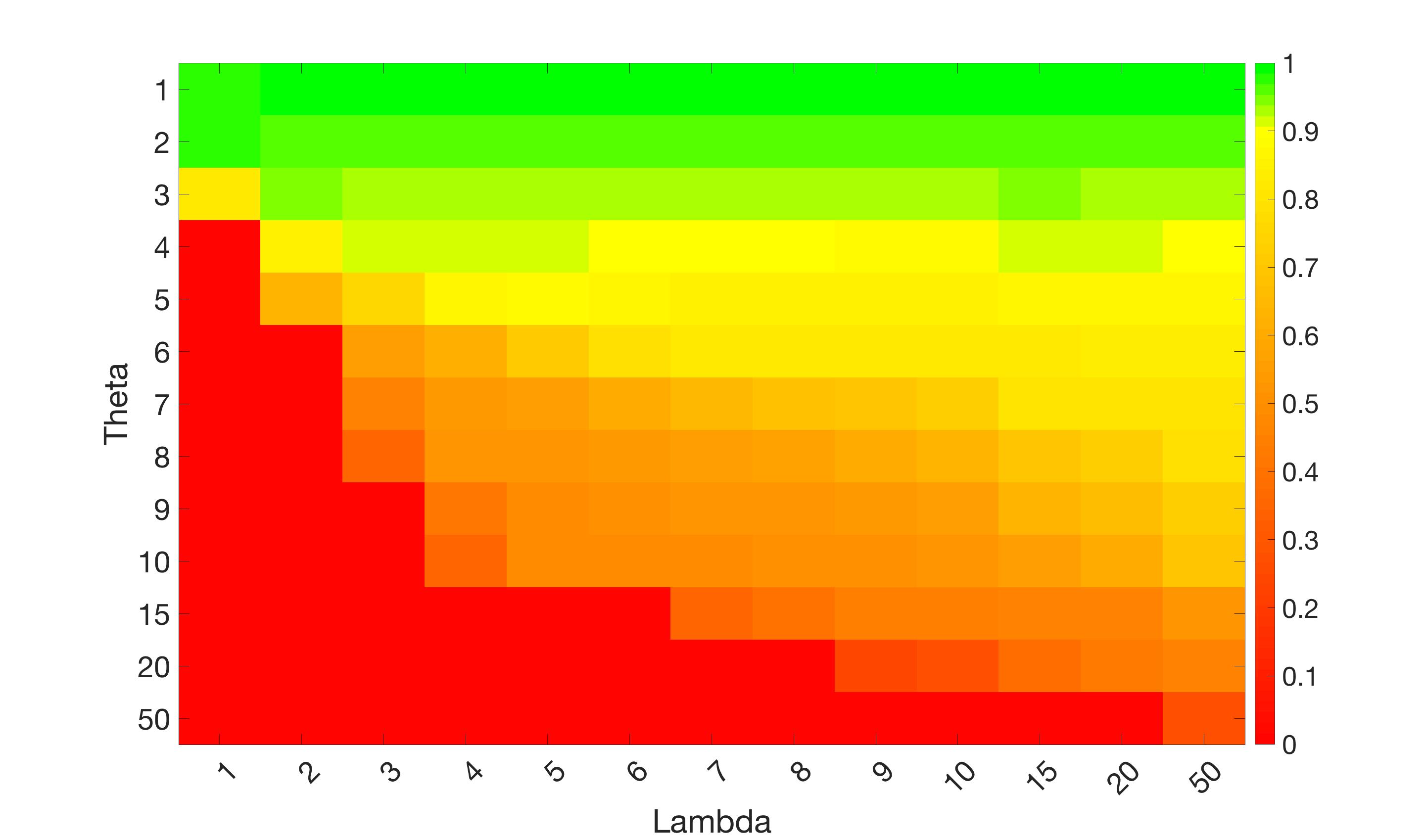}
        \caption{{\bf M6} TC values for various $\lambda$ and $\theta$ values.}
          \end{subfigure}\quad%
                \begin{subfigure}[b]{0.39\textwidth}
        \includegraphics[width=\textwidth]{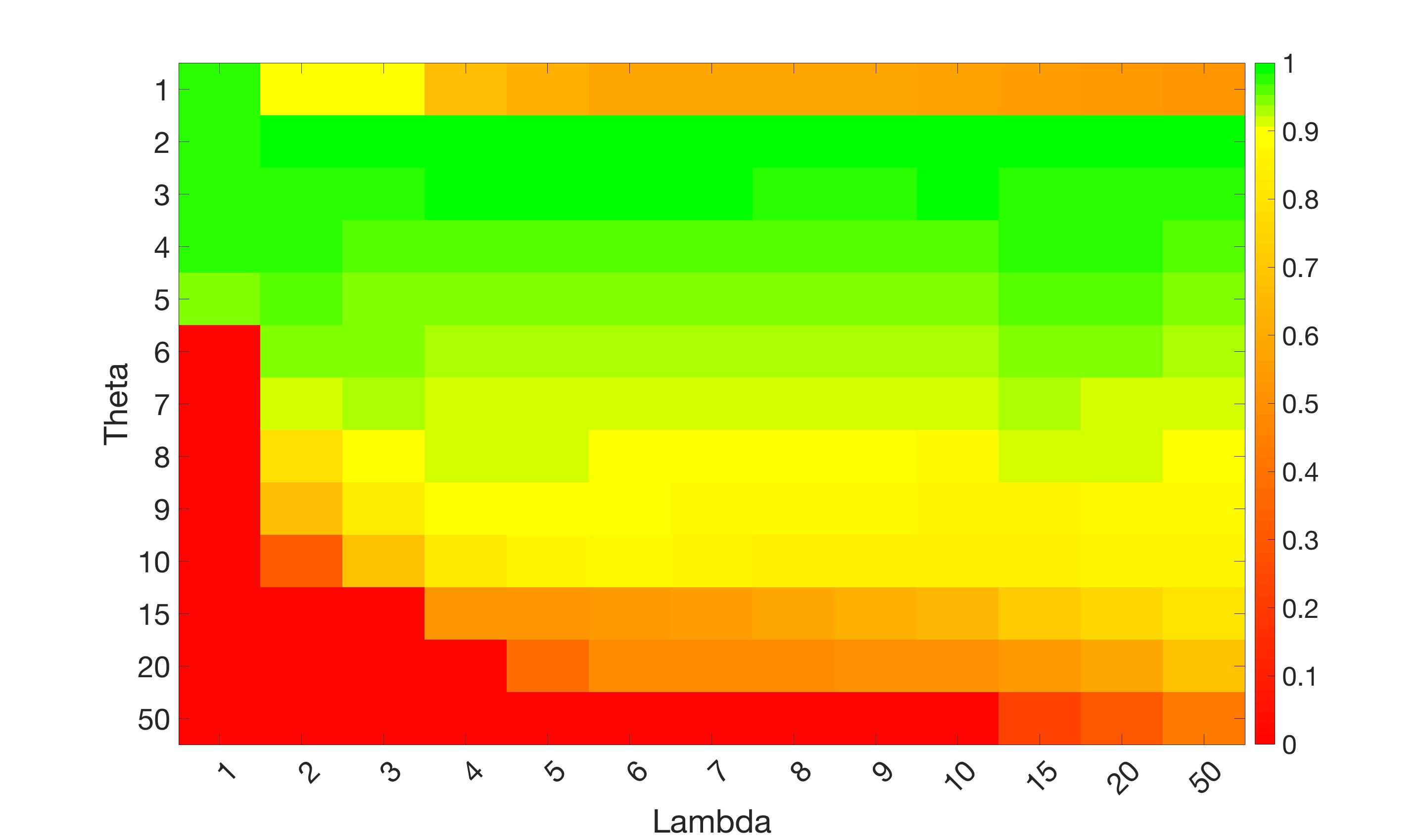}
        \caption{{\bf M7} TC values for various $\lambda$ and $\theta$ values.}
          \end{subfigure}
    }%
\caption{Parameter maps for {\bf M3}, {\bf M6} and {\bf M7}  \label{fig:kneeparams}}
\end{figure}

\begin{figure}[htb!]
\centering
\makebox[\textwidth][c]{
  \begin{subfigure}[b]{0.35\textwidth}
        \includegraphics[width=\textwidth,height=4cm]{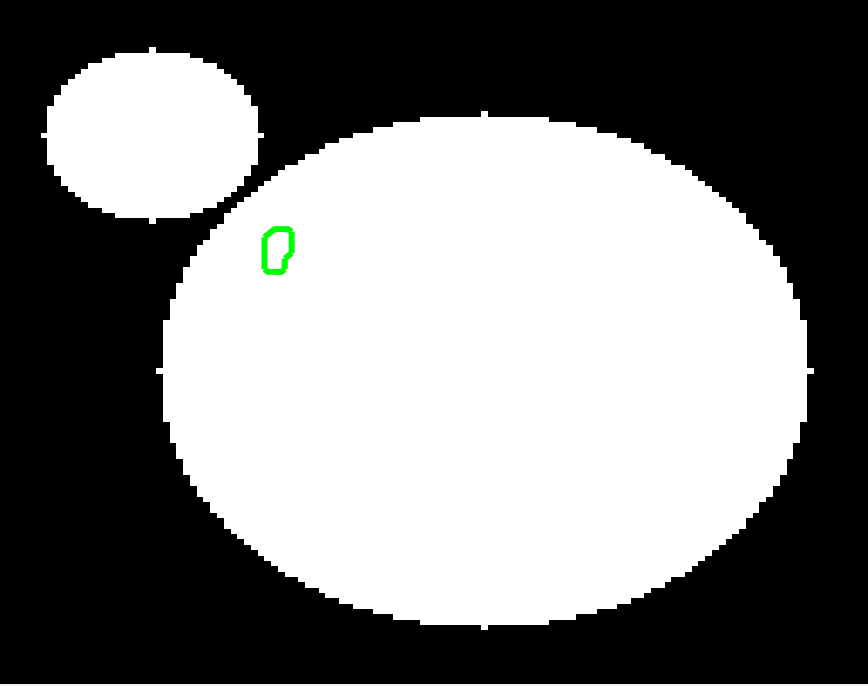}
        \caption{Original image with marker set.}
    \end{subfigure}\quad\quad\quad%
      \begin{subfigure}[b]{0.35\textwidth}
        \includegraphics[width=\textwidth,height = 4cm]{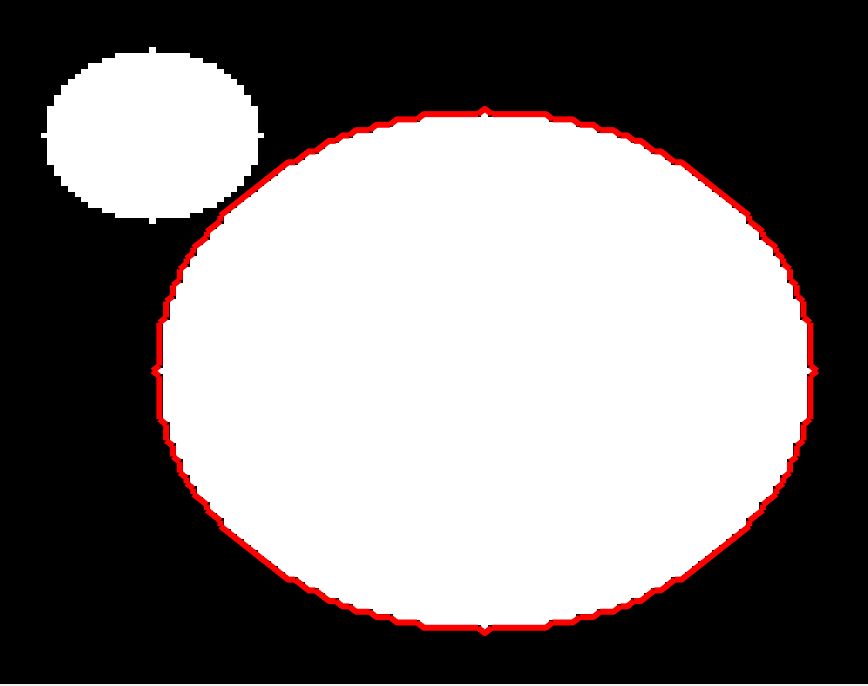}
        \caption{Ground truth segmentation.}
    \end{subfigure}    }%

\vspace{0.2in}

\centering
\makebox[\textwidth][c]{
      \begin{subfigure}[b]{0.39\textwidth}
        \includegraphics[width=\textwidth]{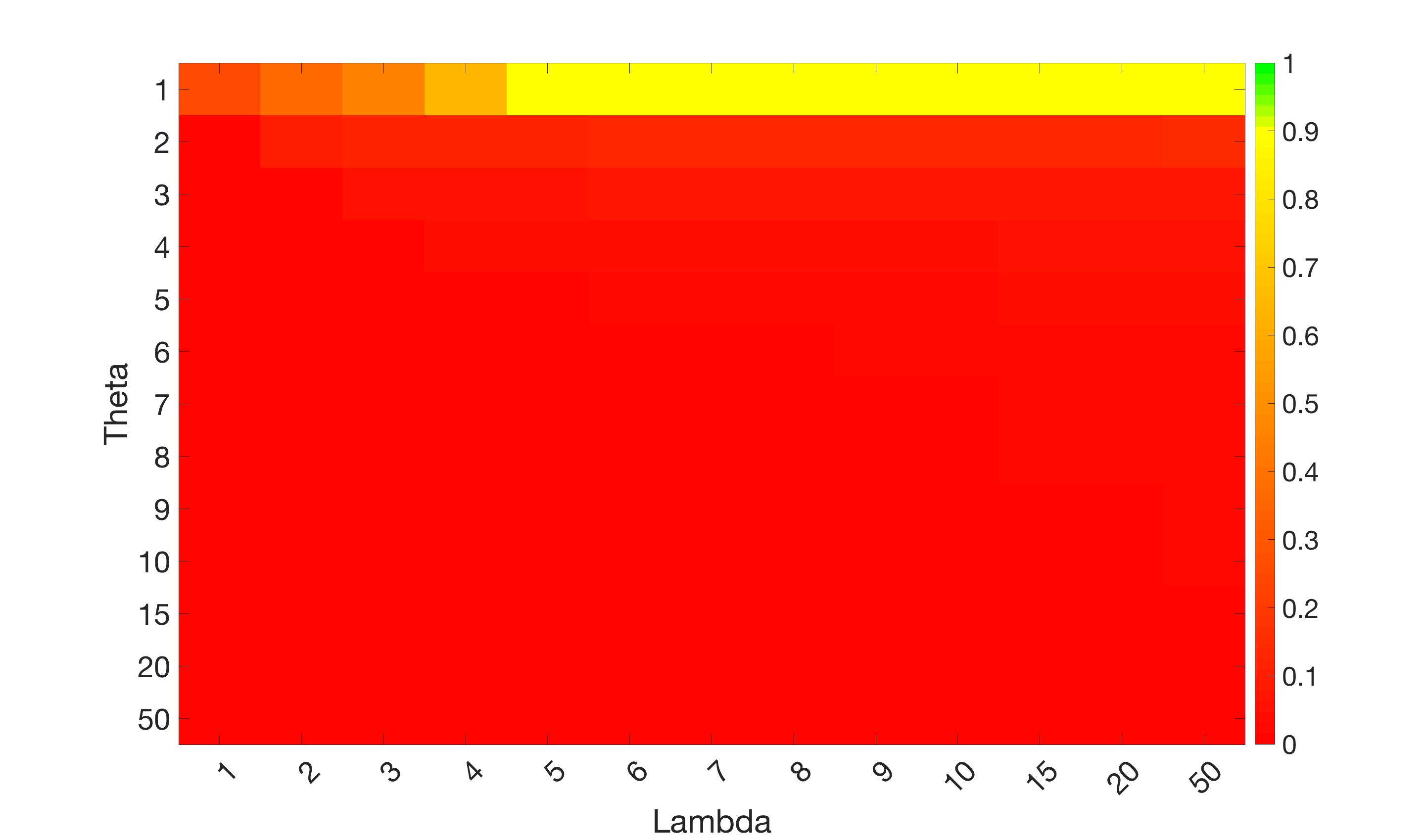}
        \caption{{\bf M3} TC values for various $\lambda$ and $\theta$ values.}
          \end{subfigure}
      \begin{subfigure}[b]{0.39\textwidth}
        \includegraphics[width=\textwidth]{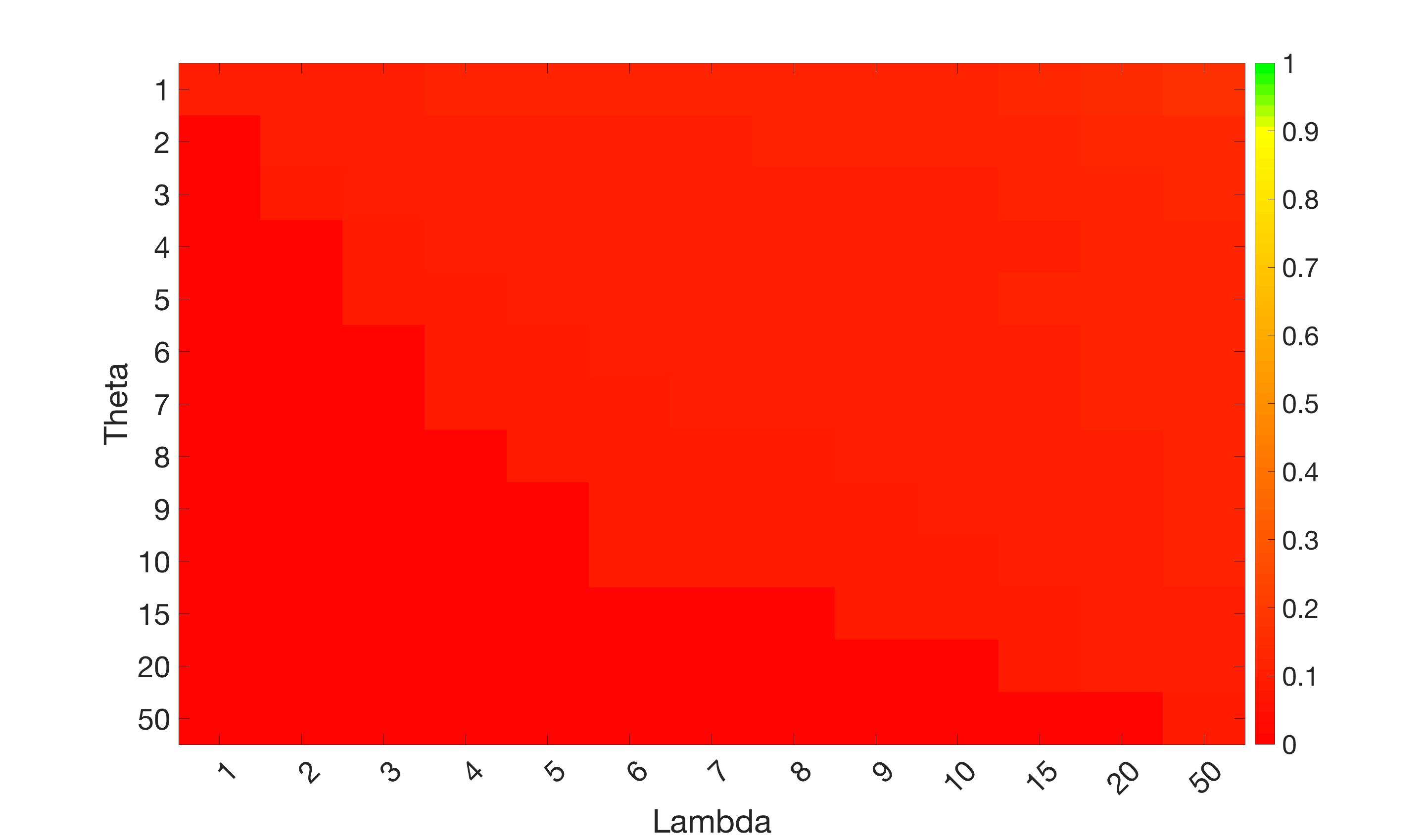}
        \caption{{\bf M6} TC values for various $\lambda$ and $\theta$ values.}
          \end{subfigure}
      \begin{subfigure}[b]{0.39\textwidth}
        \includegraphics[width=\textwidth]{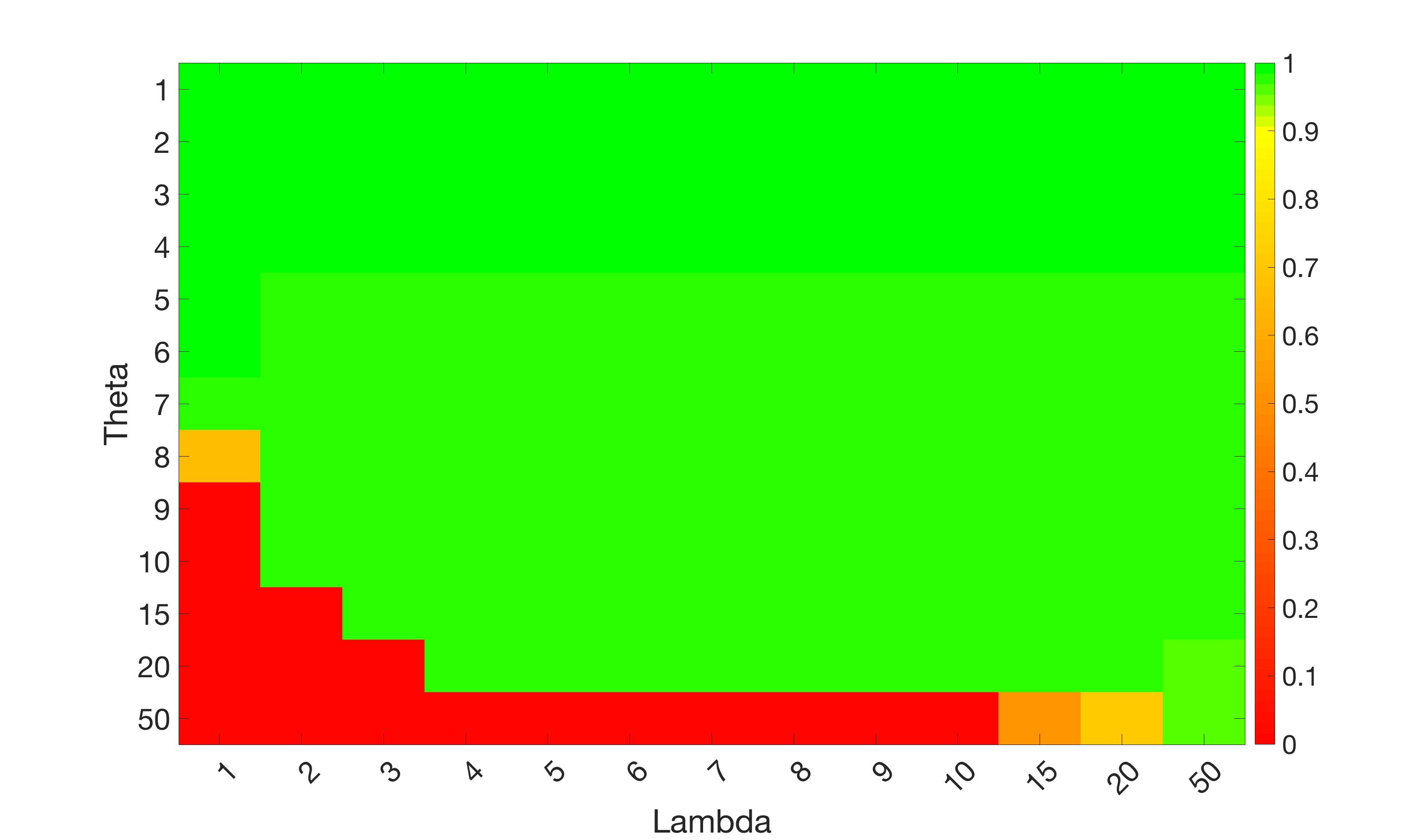}
        \caption{{\bf M7} TC values for various $\lambda$ and $\theta$ values.}
          \end{subfigure}
}%
%
\caption{Parameter maps for {\bf M3}, {\bf M6} and {\bf M7}  \label{fig:circparams}}

\end{figure}

In this test we give the segmentation results for models {\bf M1 -- M7} for the two challenging test images shown in Figure~\ref{fig:testimages}. The marker and anti-marker sets used in the experiments are also shown in this figure. After extensive parameter tuning, the best final segmentation results for each of the models are shown in Figures~\ref{fig:test1} and \ref{fig:test2}. For {\bf M1 -- M4} we obtain incorrect segmentations in both cases.
In particular, the results of {\bf M2} and {\bf M4} are interesting as the former gives poor results for both images, and the latter gives a reasonable result for Test Image 1 and a poor result for Test Image 2. In the case of {\bf M2}, the regularisation term includes the edge detector and the distance penalty term (see (\ref{eqn:rcfunc})). It is precisely this which permits the poor result in Figures~\ref{fig:test1}(b) and \ref{fig:test2}(b) as the edge detector is zero along the contour and the fitting terms are satisfied there (both intensity and area constraints) -- the distance term is not large enough to counteract the effect of these.
In the case of {\bf M4}, the distance term and edge detector are separated from the regulariser and are used to weight the Chan-Vese fitting terms (see (\ref{eqn:liufunc})). The poor segmentation in Figure~\ref{fig:test2}(b) is due to the Chan-Vese terms encouraging segmentation of bright objects (in this case), weighting $\omega$ enforces these terms at all edges in the image and near $\mathcal{M}$. In experiments, we find that {\bf M4} performs well when the object to segment is of approximately the highest or lowest intensity in the image, however when this is not the case, results tend to be poor.
We see that, in both cases, models {\bf M5} and {\bf M6} give much improved results to {\bf M2} and {\bf M4} (obtained by incorporating the geodesic distance penalty into each). The proposed Geodesic Model {\bf M7} gives an accurate segmentation in both cases. It remains to compare {\bf M5, M6} and {\bf M7}. We see that {\bf M5} is a non-convex model (and cannot be made convex \cite{SC}), therefore results are initialisation dependent. It also requires one more parameter than {\bf M6}  and {\bf M7}, and an accurate set $\mathcal{M}$ to give a reasonable area constraint in (\ref{eqn:rcfunc}).
These limitations lead us to conclude {\bf M6} and {\bf M7} are better choices than {\bf M5}.
In the case of {\bf M6}, it has the same number of parameters as {\bf M7} and gives good results.
{\bf M6} can be viewed as the model {\bf M7} with weighted intensity fitting terms (compare (\ref{eqn:geofunc2F}) and (\ref{eqn:convliucvfunc})). Experimentally, we find that the same quality of segmentation result can be achieved with both models generally, however {\bf M6} is more parameter sensitive than {\bf M7}. This can be seen in the parameter map in Figure~\ref{fig:kneeparams} with {\bf M7} giving an accurate result for a wider range of parameters than {\bf M6}. To show the improvement of {\bf M7} over previous models, we also give an image in Figure~\ref{fig:circparams} which can be accurately segmented with {\bf M7} but the correct result is never achieved with {\bf M6} (or {\bf M3}).
Therefore we find that {\bf M7} outperforms all other models tested {\bf M1 -- M6}.
\begin{rmk}
Models {\bf M2 -- M7} are coded in MATLAB and use exactly the same marker/anti-marker set. For model {\bf M1}, the software of Nguyen et al. requires marker and anti-marker sets to be input to an interface. These have been drawn as close as possible to match those used in the MATLAB code.
\end{rmk}


{\bf Test 2 -- Test of {\bf M7}'s sensitivity to changes in its main parameters.} In this test we demonstrate that the proposed Geodesic Model is robust to changes in the main parameters. The main parameters in (\ref{eqn:geoel2}) are $\mu, \lambda_{1},\lambda_{2},\theta$ and $\varepsilon_{2}$. In all tests we set $\mu = 1$, which is simply a rescaling of the other parameters, and we set $\lambda = \lambda_{1} = \lambda_{2}$. In the first example, in Figure~\ref{fig:kneeparams}, we compare the TC value for various $\lambda$ and $\theta$ values for segmentation of a bone in a knee scan. We see that the segmentation is very good for a larger
 range of $\theta$ and $\lambda$ values.
 For the second example, in Figure~\ref{fig:circparams}, we show an image and marker set for which the Spencer-Chen model ({\bf M3}) and modified Liu et al. model {\bf M6} cannot achieve the desired segmentation for any parameter range, but which can be attained for the Geodesic Model for a vast range of parameters. The final example, in Table~\ref{tab:epscomp}, compares the TC values for various $\varepsilon_{2}$ values with fixed parameters $\lambda = 2$ and $\theta = 2$. We use the images and ground truth as shown in Figures~\ref{fig:kneeparams} and \ref{fig:circparams}: on the synthetic circles image we obtain a perfect segmentation for all values of $\varepsilon_{2}$ tested, and in the case of the knee segmentation the results are almost identical for any $\varepsilon_{2} < 10^{-6}$, above which the quality slowly deteriorates.

\begin{table}
\begin{center}
\begin{tabular}{ | c | c | c | }
  \hline			
  $\varepsilon_{2}$ & Knee Segmentation (Figure~\ref{fig:kneeparams}) & Circle Segmentation (Figure~\ref{fig:circparams}) \\
  \hline
  $10^{-10}$ & 0.97287 & 1.00000 \\
  $10^{-8}$ & 0.97287 & 1.00000 \\
    $10^{-6}$ & 0.97235 & 1.00000 \\
      $10^{-4}$ & 0.96562 & 1.00000 \\
        $10^{-2}$ & 0.94463 & 1.00000 \\
          $10^{0}$ & 0.90660 & 1.00000 \\
            $10^{2}$ & 0.89573 & 1.00000 \\
              $10^{4}$ & 0.89159 & 1.00000 \\
  \hline
\end{tabular}
\caption{The Tanimoto Coeffcient for various $\varepsilon_{2}$ values, segmenting the images in Figures~\ref{fig:kneeparams} and \ref{fig:circparams}. \label{tab:epscomp}}
\end{center}
\end{table}

\begin{figure}[h!tbp]
\centering
\makebox[\textwidth][c]{
\includegraphics[width=0.30\textwidth,height = 4cm]{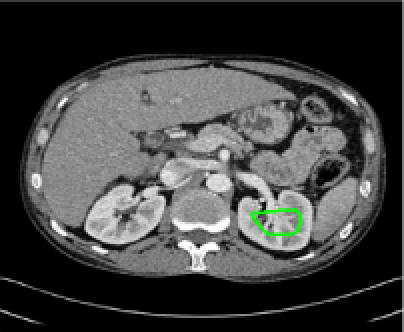}\hspace*{0.2in} 
\includegraphics[width=0.30\textwidth,height = 4cm]{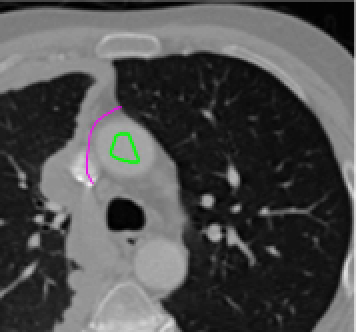}\hspace*{0.2in} 
\includegraphics[width=0.30\textwidth,height = 4cm]{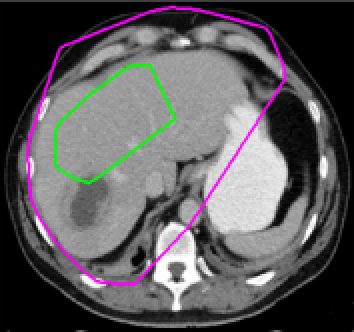}
}

\vspace*{0.2in}

\makebox[\textwidth][c]{
\includegraphics[width=0.30\textwidth,height = 4cm]{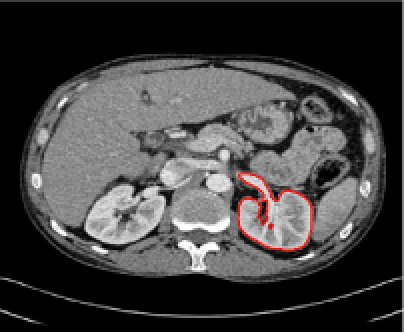}\hspace*{0.2in} 
  \includegraphics[width=0.30\textwidth,height = 4cm]{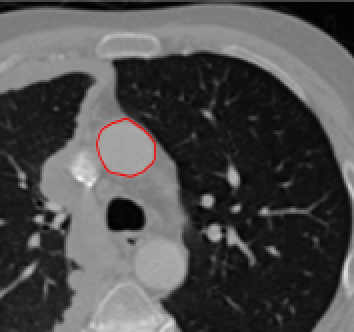}\hspace*{0.2in} 
\includegraphics[width=0.30\textwidth,height = 4cm]{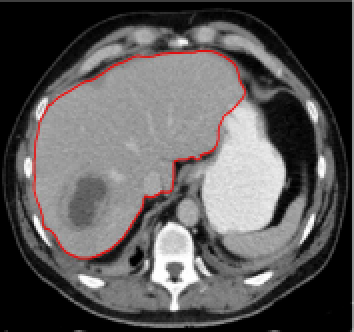}
}

\vspace*{0.2in}

\makebox[\textwidth][c]{
\includegraphics[width=0.30\textwidth,height = 4cm]{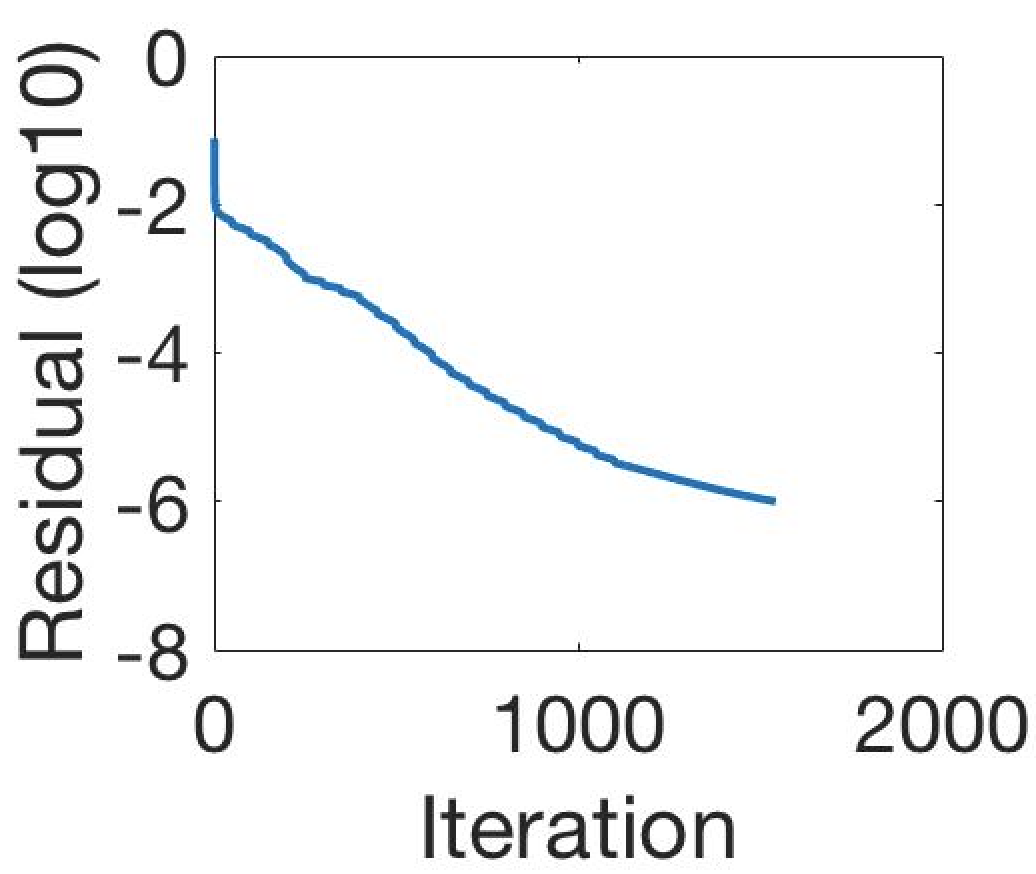}\hspace*{0.2in} 
  \includegraphics[width=0.30\textwidth,height = 4cm]{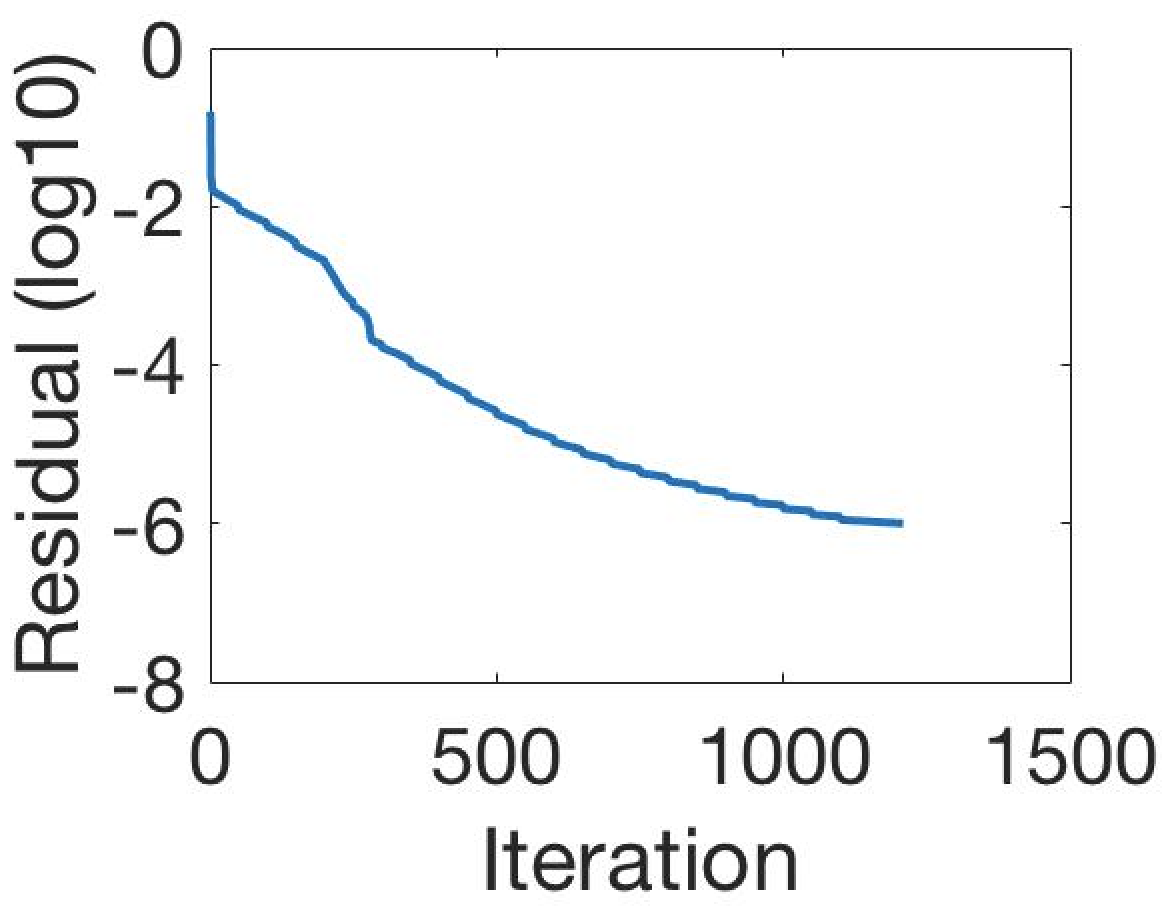}\hspace*{0.2in} 
\includegraphics[width=0.30\textwidth,height = 4cm]{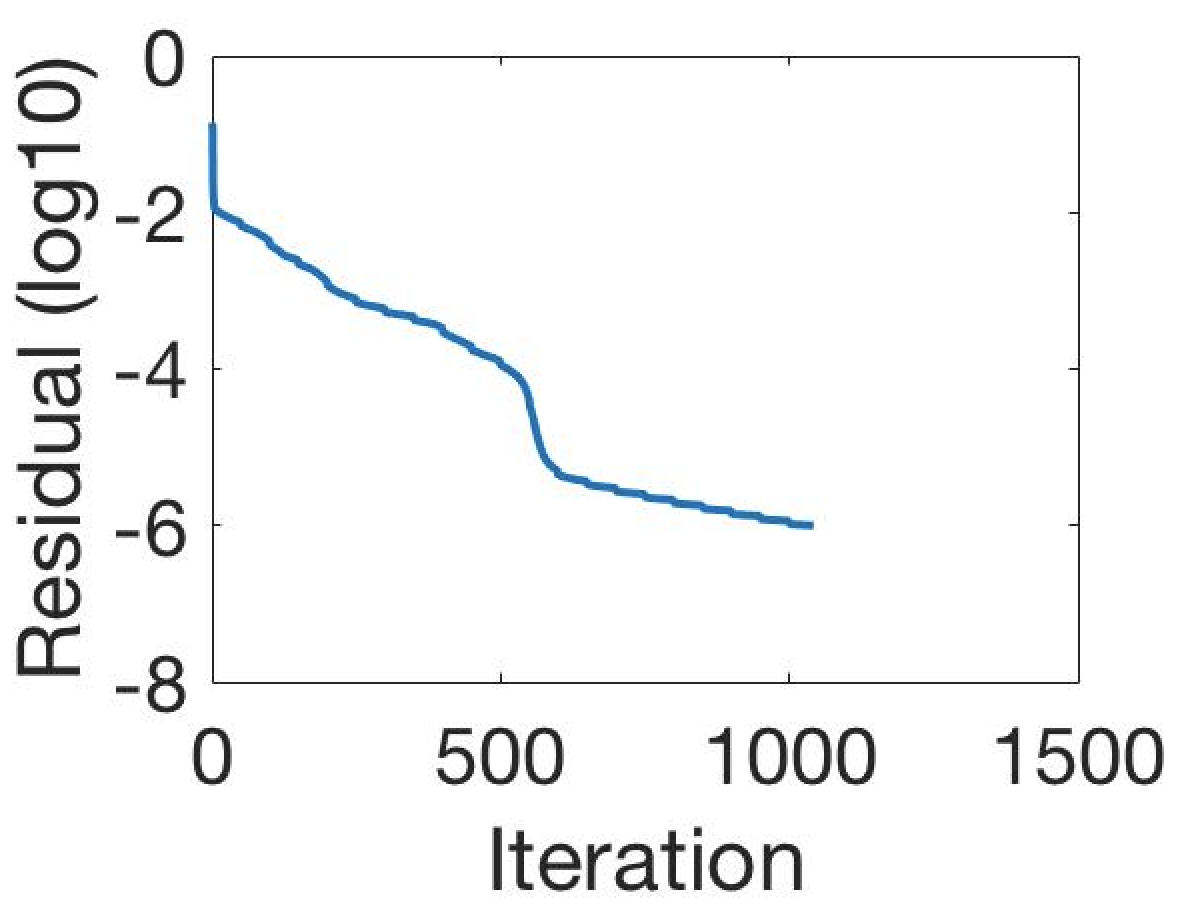}
}




\caption{Three further test results obtained using our Geodesic Model {\bf M7}, all with parameters $\theta =5$, $\lambda = 5$. The  first row  shows the original image with the marker set (plus anti-marker set),
the second row the final segmentation result and the final row shows the residual history.
 \label{fig:georesults}}

\end{figure}

{\bf Test 3 -- Further Results from the Geodesic Model {\bf M7}.} In this test
we give some medical segmentation results obtained using the Geodesic Model {\bf M7}. The results are shown in Figure~\ref{fig:georesults}.
In the final two columns
 we use anti-markers to demonstrate how to overcome blurred edges and low contrast edges in an image. These are challenging and it is pleasing to see the correctly segmented results.

\begin{figure}[htb!]
\centering
\makebox[\textwidth][c]{
\includegraphics[width=0.24\textwidth,height = 4cm]{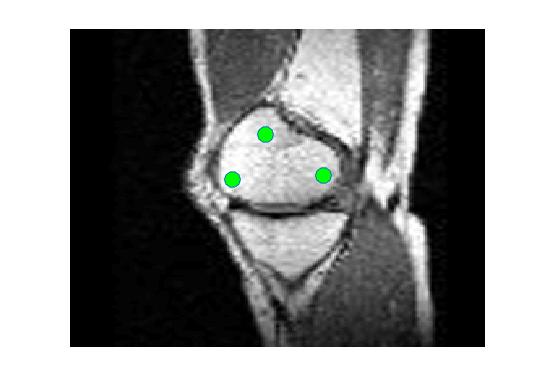}\hspace*{-0.2in}
\includegraphics[width=0.24\textwidth,height = 4cm]{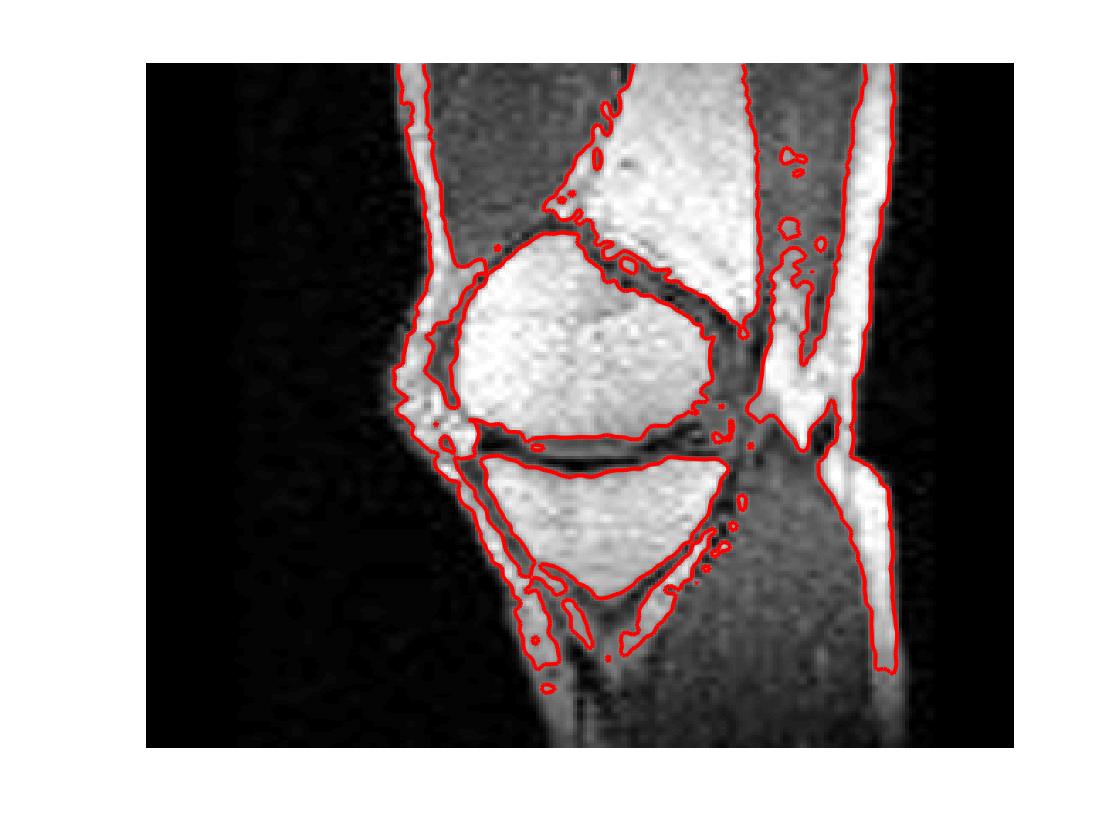}\hspace*{-0.2in}
\includegraphics[width=0.24\textwidth,height = 4cm]{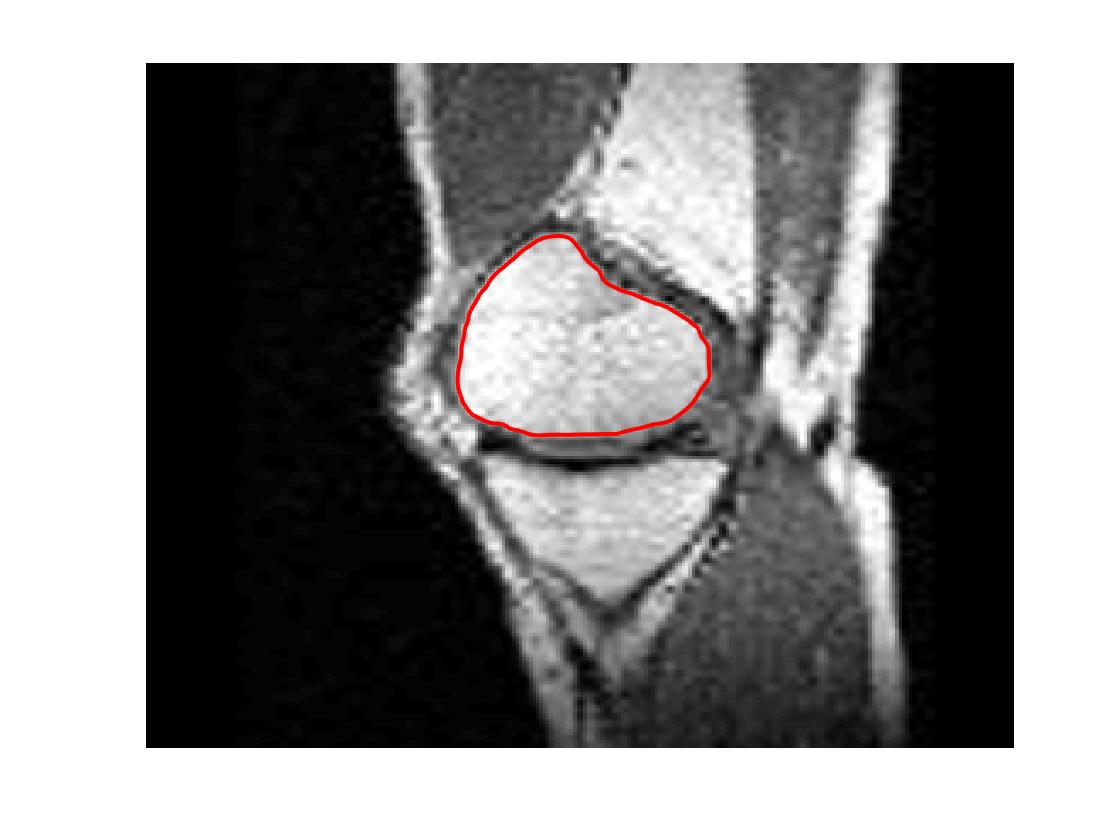}\hspace*{-0.2in}
\includegraphics[width=0.24\textwidth,height = 4cm]{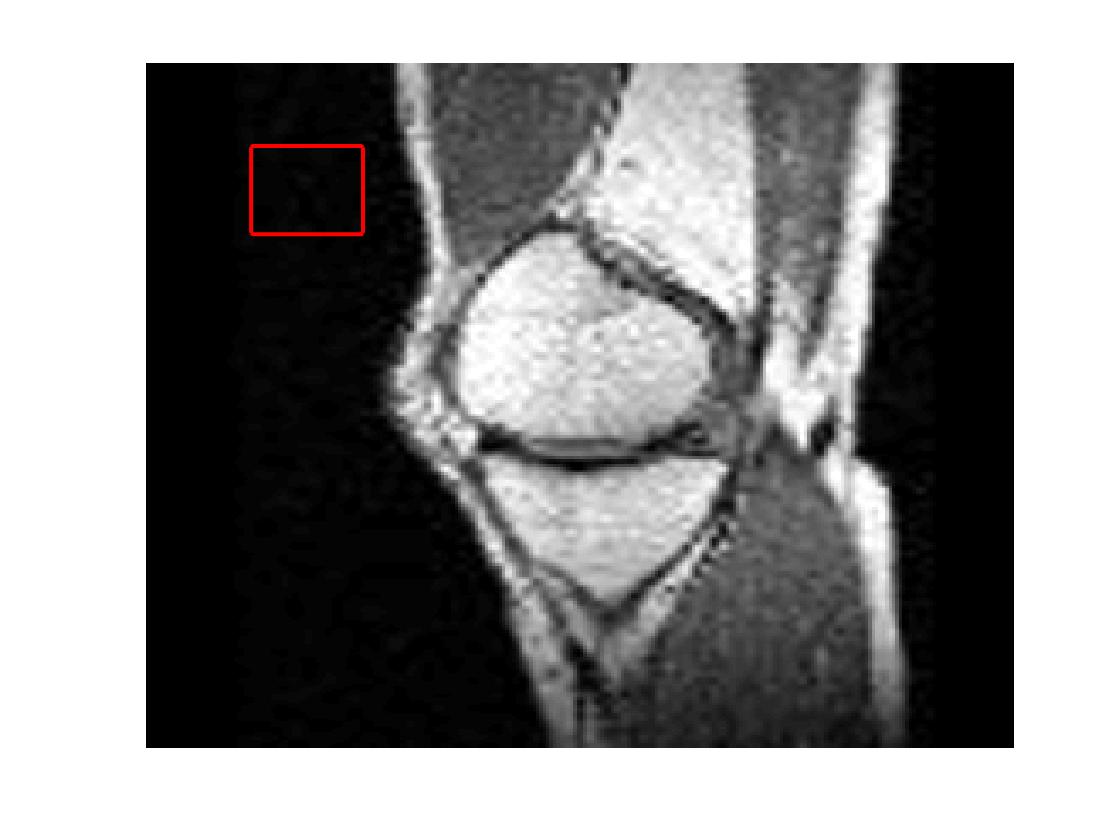}\hspace*{-0.2in}
\includegraphics[width=0.24\textwidth,height = 4cm]{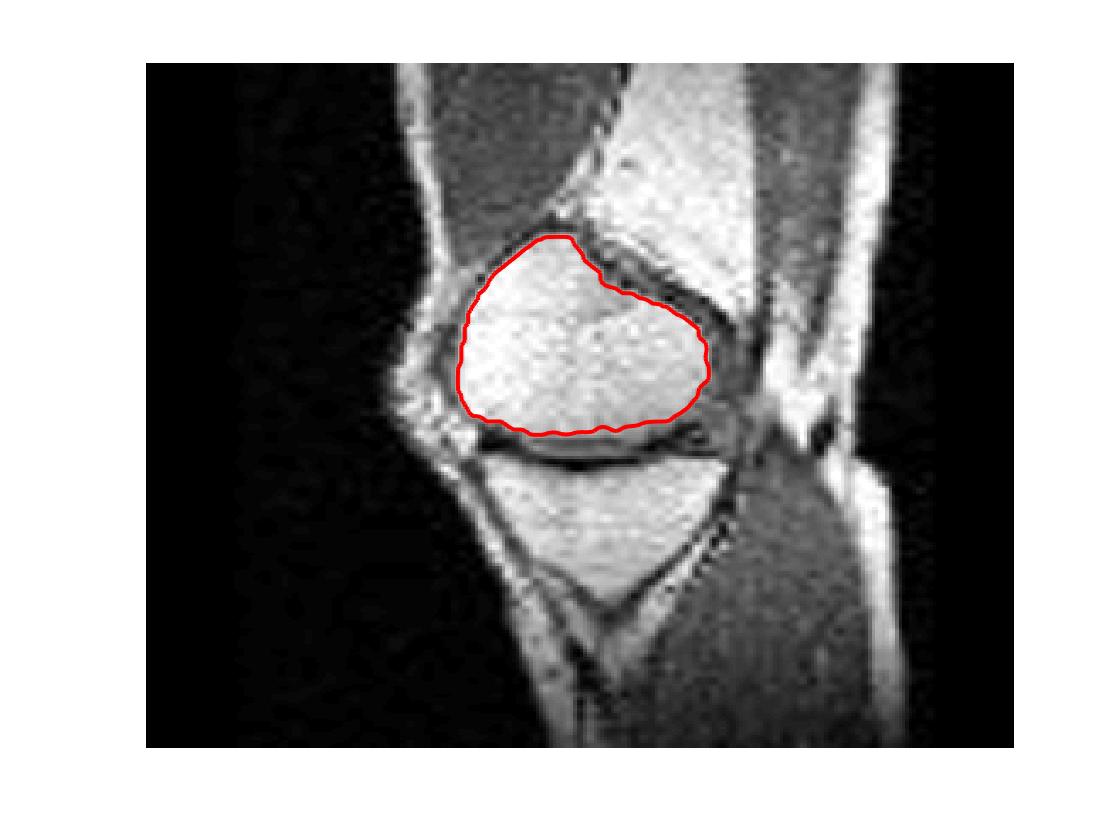}
}
\centerline{
(i)\hspace*{0.175\textwidth}
(ii)\hspace*{0.175\textwidth}
(iii)\hspace*{0.175\textwidth}
(iv)\hspace*{0.175\textwidth}(v)}
\caption{Test $4$ on {\bf M7}'s initialisations ($\theta = 5, \lambda = 5$). \ (i)  The original image with marker set indicated; (ii) Initialisation $1$ using the image itself; (iii) Segmentation result from Initialisation $1$;
(iv) Initialisation $2$ away from the object to be segmented; (v)  Segmentation $2$ from initialisation $2$.
Clearly {\bf M7} gives the same result.} \label{fig:differentinit}

\centering
\makebox[\textwidth][c]{
\includegraphics[width=0.35\textwidth,height = 3cm]{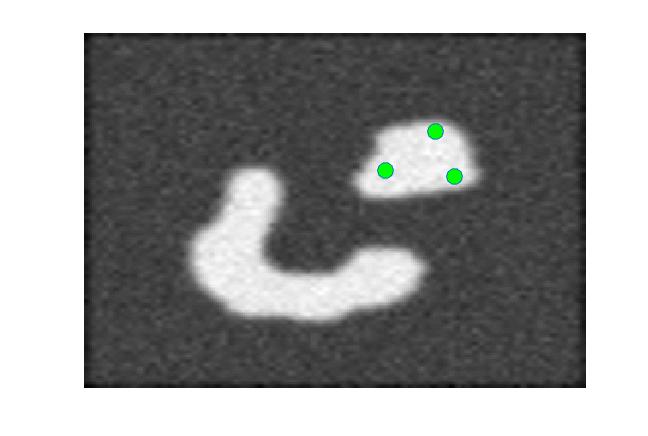}\hspace*{-0.2in}
\includegraphics[width=0.35\textwidth,height = 3cm]{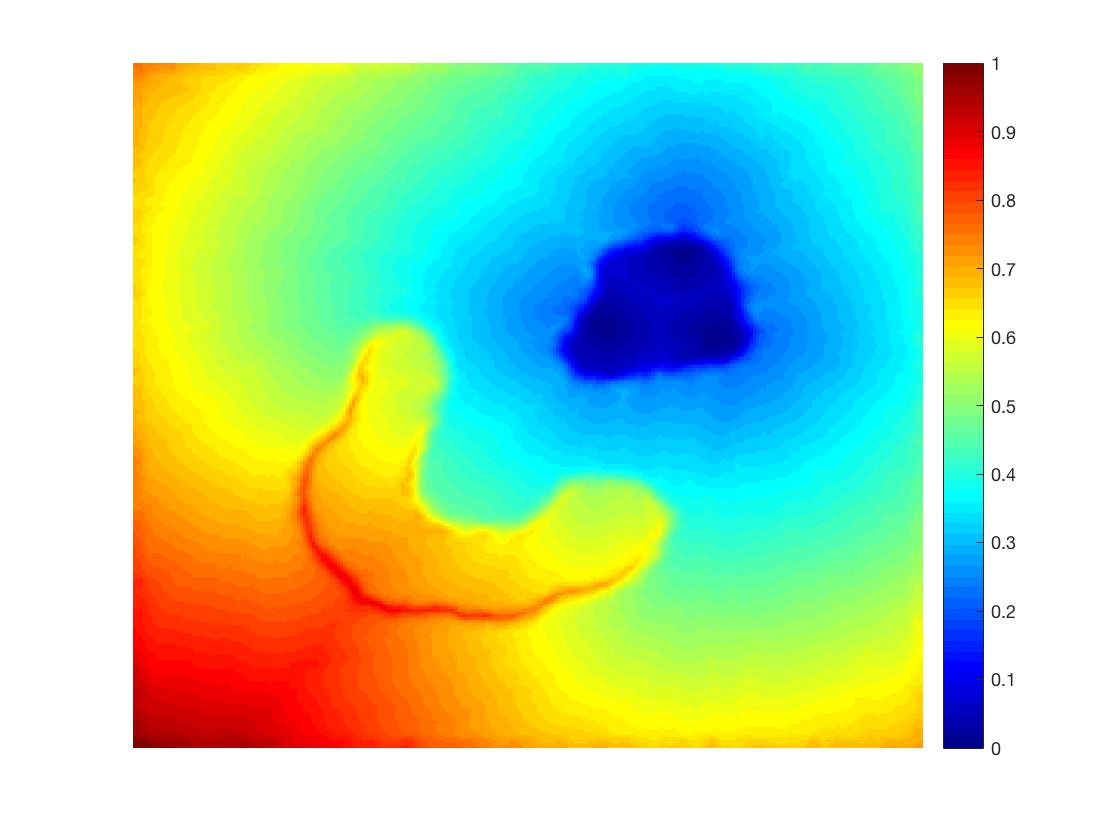}\hspace*{-0.2in}
\includegraphics[width=0.35\textwidth,height = 3cm]{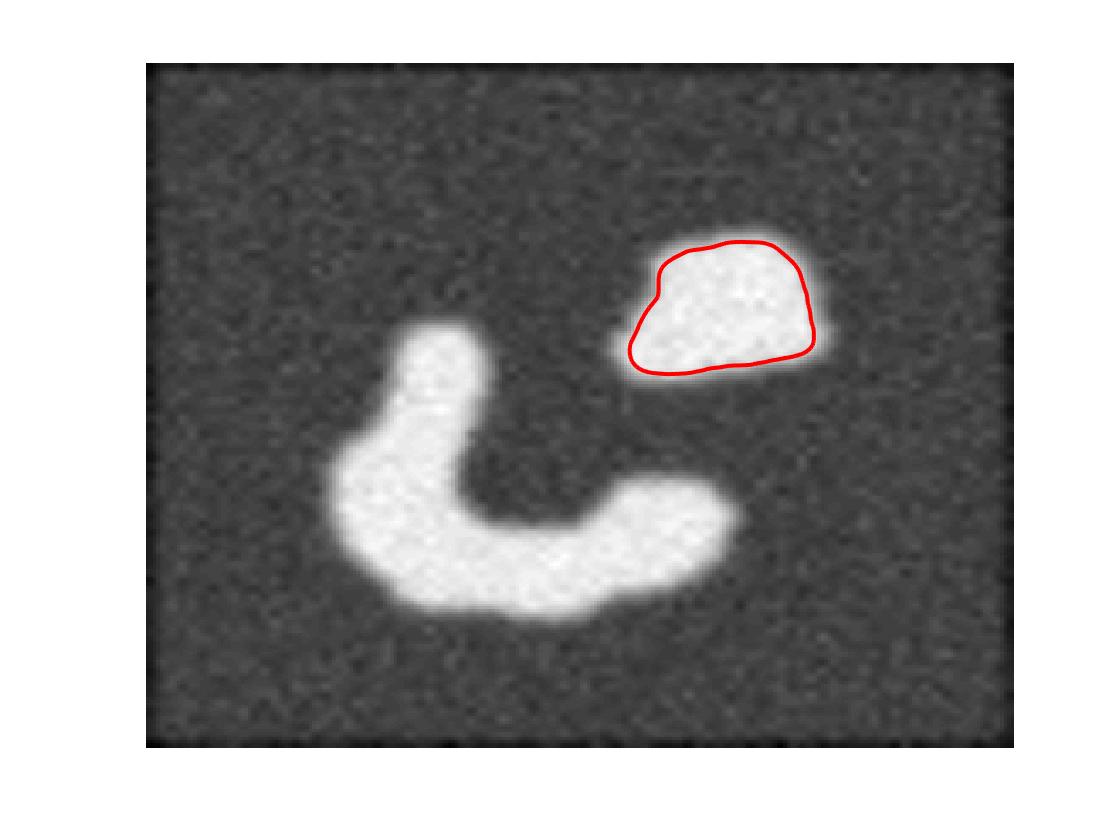}
}
\centering
\makebox[\textwidth][c]{
\includegraphics[width=0.35\textwidth,height = 3cm]{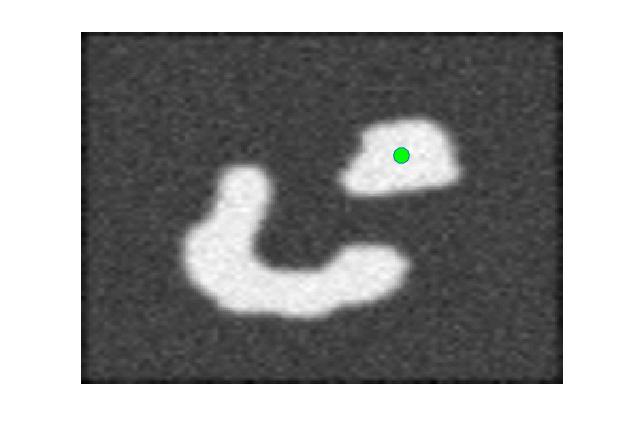}\hspace*{-0.2in}
\includegraphics[width=0.35\textwidth,height = 3cm]{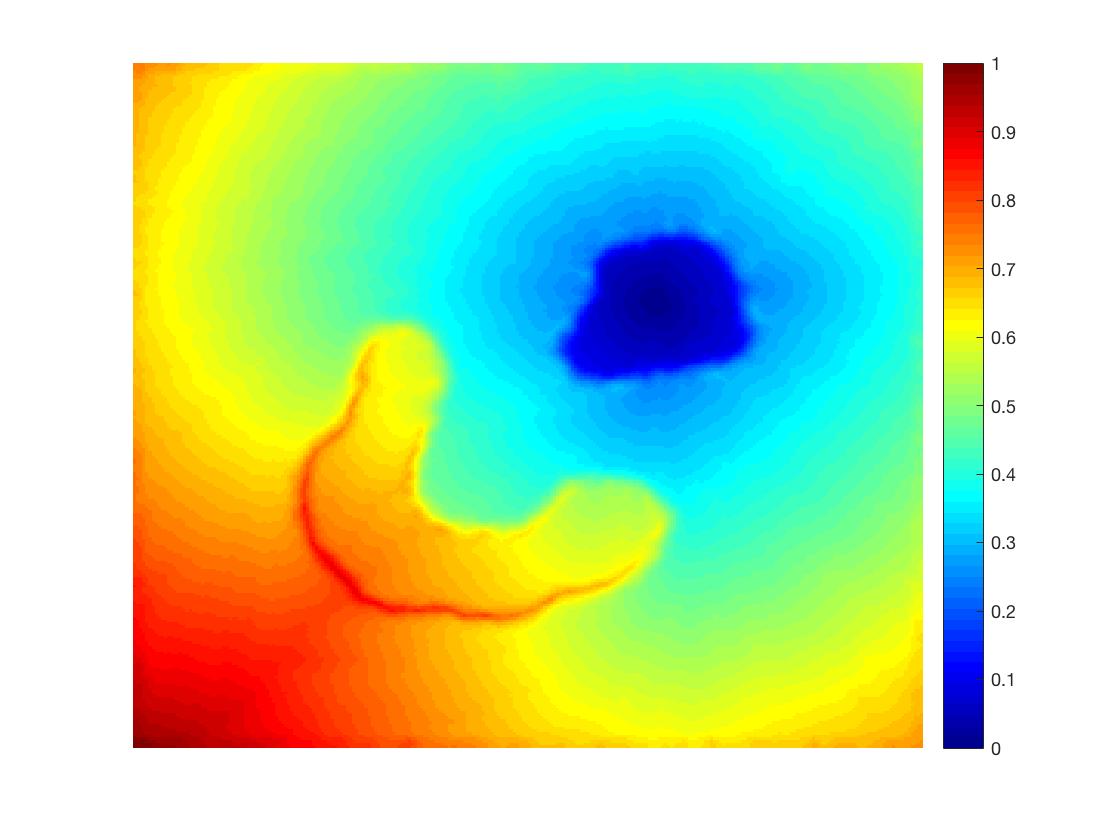}\hspace*{-0.2in}
\includegraphics[width=0.35\textwidth,height = 3cm]{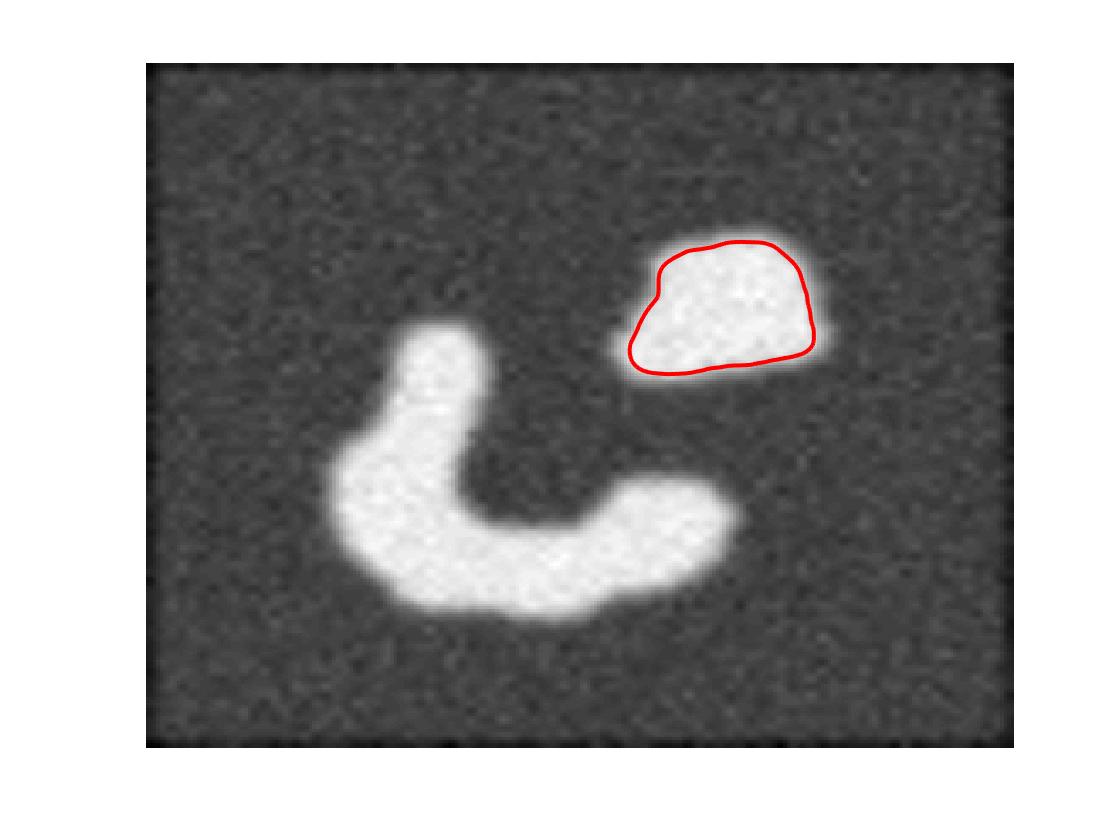}
}
\caption{Test $4$ on {\bf M7}'s marker set ($\theta = 5, \lambda = 3$).
Row $1$  shows the original image with $3$ marker points, the normalised geodesic distance map and the final segmentation result.
Row $2$  shows the original image with $1$ marker point, the normalised geodesic distance map and the final segmentation result. Clearly the second and third columns are the same
for different marker points. Thus {\bf M7} is robust. }
 \label{fig:differentmarkers}
\end{figure}

{\bf Test 4 -- Initialisation and Marker Set Independence.} In the first example, in Figure~\ref{fig:differentinit}, we see how the convex Geodesic Model {\bf M7} gives the same segmentation result regardless of initialisation, as expected of a convex model. Hence the model is flexible in implementation. From many experiments it is found that using the polygon formed by the marker points as the initialisation converges to the final solution faster than using an arbitrary initialisation. In the second example, in Figure~\ref{fig:differentmarkers}, we show intuitively how   Model {\bf M7} is robust to the number of markers and the location of the markers within the object to be segmented.
The Euclidean distance term, used in the Spencer-Chen model {\bf M3}, is sensitive to the position and number of marker points, however, regardless of where the markers are chosen, and how many are chosen, the geodesic distance map will be almost identical.


\section{Conclusions}

In this paper a new convex selective segmentation model has been proposed, using geodesic distance as a penalty term. This model gives results that are unachievable by alternative selective segmentation models and is also more robust to the parameter choices. Adaptations to the penalty term have been discussed which make it robust to noisy images and blurry edges whilst also penalising objects far from the marker set (in a Euclidean distance sense). A proof for the existence and uniqueness of the viscosity solution to the PDE given by the Euler-Lagrange equation for the model has been given (which applies to an entire class of image segmentation PDEs). Finally we have confirmed the advantages of using the geodesic distance with some experimental results.
Future works will look for further extension of selective segmentation to other frameworks such as using high order regularizers \cite{tai13,duan16} where only incomplete theories exist.

\section*{Acknowledgements}
The first author wishes to thank the UK EPSRC, the Smith Institute for Industrial Mathematics, and the Liverpool Heart and Chest Hospital for
supporting the work through an Industrial CASE award. Thanks also must go to Professor Joachim Weickert (Saarland, Germany)  for fruitful discussions at the early stages of this work.
The second author is grateful to the UK
EPSRC for the grant EP/N014499/1.



\section*{{\bf Appendix A} --- Proof that Condition (I7) Holds
                   in Theorem \ref{thm:no2}}

Using the assumption in (\ref{eqn:assumpI3}), we write
\begin{equation*}
\begin{gathered}
\begin{aligned}
(Xr,r) + (Ys,s) = r^{T}Xr + s^{T}Ys  &\le \mu_{1} \begin{bmatrix}
r^{T} &  s^{T}
\end{bmatrix}
\begin{bmatrix}
I & -I \\
-I & I \\
\end{bmatrix} \begin{bmatrix}
r \\
s\\
\end{bmatrix}+\mu_{2} \begin{bmatrix}
r^{T} &  s^{T}
\end{bmatrix}
\begin{bmatrix}
I &0\\
0& I \\
\end{bmatrix}\begin{bmatrix}
r \\
s\\
\end{bmatrix}\\
& = \mu_{1}|r-s|^{2} + \mu_{2}\left(|r|^{2}+|s|^{2}\right).\\
\end{aligned}
\end{gathered}
\end{equation*}
Note that matrix
$A$ from (\ref{eqn:defA}) is a real symmetric matrix and decomposes as $A = QDQ^{T} = QD^{1/2}D^{1/2}Q^{T} = BB^{T}$ with $Q$ orthonormal and $B = QD^{1/2}$. Successively define $r = B(p)e_{i}$ and $s = B(p)e_{i}$ for all $(e_{i})$, an orthonormal basis, and obtain
\begin{equation*}
\begin{gathered}
\begin{aligned}
(Xr,r) = r^{T}Xr &= \sum_{i}(Be_{i})^{T}X(Be_{i}) = \sum_{i}e_{i}^{T}B^{T}XBe_{i} = \text{trace}(B^{T}XB) = \text{trace}(A(x,p)X).
\end{aligned}
\end{gathered}
\end{equation*}
Therefore, we can write
\begin{equation*}
\begin{gathered}
\begin{aligned}
\text{trace}(A(x,p)X) + \text{trace}(A(y,q)Y)  &= (XB(p)e_{i},B(p)e_{i}) + (YB(q)e_{i},B(q)e_{i})\\
&\le \mu_{1} | B(p)e_{i} - B(q)e_{i} |^{2} + \mu_{2}\left( |B(p)e_{i}|^{2} +|B(q)e_{i}|^{2}   \right)\\
& = \mu_{1}\text{trace}\left( \left( B(p)- B(q)\right)^{T} \left( B(p) - B(q)\right)  \right) + \mu_{2}\left( G(x) + G(y) \right).\\
\end{aligned}
\end{gathered}
\end{equation*}
We now focus on reformulating the first term, we start by decomposing $A(x,p)$ as follows
\[
A(x,p) =
\begin{bmatrix}
\frac{p_{1}}{|p|} & -\frac{p_{2}}{|p|} \\
\frac{p_{2}}{|p|} & \frac{p_{1}}{|p|} \\
\end{bmatrix}
\begin{bmatrix}
0 & 0\\
0& G(x)\\
\end{bmatrix}
\begin{bmatrix}
\frac{p_{1}}{|p|} & \frac{p_{2}}{|p|} \\
-\frac{p_{2}}{|p|} & \frac{p_{1}}{|p|} \\
\end{bmatrix}
=
\begin{bmatrix}
\frac{p_{1}}{|p|} & -\frac{p_{2}}{|p|} \\
\frac{p_{2}}{|p|} & \frac{p_{1}}{|p|} \\
\end{bmatrix}
\begin{bmatrix}
0 & 0\\
0& \sqrt{G(x)}\\
\end{bmatrix}
\begin{bmatrix}
0 & 0\\
0& \sqrt{G(x)}\\
\end{bmatrix}
\begin{bmatrix}
\frac{p_{1}}{|p|} & \frac{p_{2}}{|p|} \\
-\frac{p_{2}}{|p|} & \frac{p_{1}}{|p|} \\
\end{bmatrix}
\]
so we have $A=BB^{T}$ where
\[
B(p) =
\begin{bmatrix}
0 & -\frac{p_{2}}{|p|}\sqrt{G(x)} \\
0 & \frac{p_{1}}{|p|}\sqrt{G(x)}  \\
\end{bmatrix}.
\]
Using this we compute
\begin{equation*}
\begin{gathered}
\begin{aligned}
\text{trace}\left( \left( B(p)- B(q)\right)^{T} \left( B(p) - B(q)\right)  \right) = \left| \frac{p}{|p|}\sqrt{G(x)}-\frac{q}{|q|}\sqrt{G(y)}  \right|^{2}.
\end{aligned}
\end{gathered}
\end{equation*}
Substituting this in the overall trace sum we have
\begin{equation*}
\begin{gathered}
\begin{aligned}
\text{trace}(A(x,p)X) + \text{trace}(A(y,q)Y)  \le  \mu_{1}\left| \frac{p}{|p|}\sqrt{G(x)}-\frac{q}{|q|}\sqrt{G(y)}  \right|^{2} + 2\mu_{2}\theta. \\
\end{aligned}
\end{gathered}
\end{equation*}
as $G(x)<\theta$ ($G$ is bounded) for all $x\in\Omega$. Focussing on the first term in this expression we compute
\begin{equation*}
\begin{gathered}
\begin{aligned}
\left| \frac{p}{|p|}\sqrt{G(x)}-\frac{q}{|q|}\sqrt{G(y)}  \right|^{2}  &= \left| \frac{p}{|p|}\sqrt{G(x)}-\frac{p}{|p|}\sqrt{G(y)}+ \frac{p}{|p|}\sqrt{G(y)} -\frac{q}{|q|}\sqrt{G(y)}  \right|^{2} \\
& \le 2\left(\sqrt{G(x)}-\sqrt{G(y)}\right)^{2} + 2G(y)\left|\frac{p}{|p|}-\frac{q}{|q|}\right|^{2} \\
& \le 2\left(\sqrt{G(x)}-\sqrt{G(y)}\right)^{2} + 8\theta\rho(p,q)^{2} \\
\end{aligned}
\end{gathered}
\end{equation*}
where $\rho = \min\left(\frac{|p-q|}{\min(|p|,|q|)},1\right)$. This uses inequality $\left|\frac{p}{|p|}-\frac{q}{|q|}\right|^{2} \le 2\rho(p,q)$ (see \cite{RC15,gout2005segmentation,RC16,guillot2007existence,le2004imagerie,lavdie2013variational}).
We now note that $g(s) = \frac{1}{1+s^{2}}$ is Lipschitz continuous with Lipschitz constant $\frac{3\sqrt{3}}{8}$. 

{\bf Note.} In the Geodesic Model we fix $G(x) = g(|\nabla z|)$. Therefore, assuming $G(x)$ and $\sqrt{G(x)}$ as Lipschitz requires us to assume that the underlying $z$ is a smooth function \cite{gout2005segmentation}. Thankfully, $z$ is typically provided as a smoothed image after some filtering (e.g. Gaussian smoothing) and we can assume regularity of $z$.

\begin{rmk}
It is less clear that $\sqrt{G(x)}$ is Lipschitz, we now prove it explicitly.
Firstly, it is relatively easy to prove that
\[
\sqrt{G(x)}-\sqrt{G(y)} \le \frac{2}{3\sqrt{3}}\bigg|\left|\nabla z(x)\right|-\left|\nabla z(y)\right|\bigg|
\]
by letting $K(s) = \sqrt{g(s)}$ and we find $\sup\limits_{s}|K'(s)| = \frac{2}{3\sqrt{3}}$. We now need to prove that $|\nabla z(x)|$ is Lipschitz also. Take $h(x) = |\nabla z(x)|$, then by a remark in \cite{gout2005segmentation}, we can conclude $\exists\,\zeta <\infty$ such that
\[
\bigg|\left|\nabla z(x)\right|-\left|\nabla z(y)\right|\bigg| \le \zeta | x-y|
\]
and so $\sqrt{G(x)}$ is Lipschitz with constant $\frac{2}{3\sqrt{3}}\zeta$.
\end{rmk}

After some computations we obtain
\[
\left| \frac{p}{|p|}\sqrt{G(x)}-\frac{q}{|q|}\sqrt{G(y)}  \right|^{2}  \le  2\left(\frac{2}{3\sqrt{3}}\zeta\right)^{2}\left|x-y\right|^{2} + 8\theta\rho(p,q)^{2} = \frac{8}{27}\zeta^{2}\left|x-y\right|^{2} + 8\theta\rho(p,q)^{2}.
\]
Following the results in \cite{RC15,gout2005segmentation,RC16,guillot2007existence,le2004imagerie,lavdie2013variational} we have
\[
\left| \nabla G(x)-\nabla G(y)\right|\left| p \right| < \kappa |p||x-y| \le \kappa\max(|p|,|q|)|x-y|.
\]
so overall
\[
\langle \nabla G(x), p \rangle - \langle \nabla G(y), q \rangle  \le  \kappa\max(|p|,|q|)|x-y| + \eta |p-q|\\
\]
where $|\nabla G(y)| < \eta <\infty$. Finally, we note that $-\left(|p|-|q| \right) = |q|-|p| \le \Big| |q|-|p| \Big| \le |p-q|$. If we now write
\begin{equation*}
\begin{gathered}
\begin{aligned}
- \left( F(t,x,u,p,X)-F(t,y,u,q,-Y) \right) =&  \mu\left(\text{trace}(A(x,p)X) + \text{trace}(A(y,q)Y) \right)\\
& + \mu\left( \langle \nabla G(x), p \rangle - \langle \nabla G(y), q \rangle \right)\\
&- \left(|p|-|q|\right)k(u) - |p|f(x) + |q|f(y)\\
\le & ~\mu\mu_{1}\left(\frac{8}{27}\zeta^{2}|x-y|^{2} + 8\theta\rho(p,q)^{2} \right)+ 2\mu\mu_{2}\theta\\
& + \mu\kappa\max(|p|,|q|) |x-y|+ \mu\eta |p-q|\\
& - \left(|p|-|q|\right)\left(k(u) + 2\max\limits_{x\in\Omega}f(x) \right)\\
\le & ~\mu\mu_{1}\left(\frac{8}{27}\zeta^{2}|x-y|^{2} + 8\theta\rho(p,q)^{2} \right)+ 2\mu\mu_{2}\theta\\
& + \mu\kappa\left(\max\left(|p|,|q|\right)+1\right) |x-y|+ \mu\eta |p-q|+ C_{1}|p-q|.
\end{aligned}
\end{gathered}
\end{equation*}
where $C_{1} = \max\limits_{x\in\Omega}\left(k(u) + 2\max\limits_{x\in\Omega}f(x) \right)$ (we must assume $k(u),f(x)$ are bounded). Hence we have
\begin{equation*}
\begin{gathered}
\begin{aligned}
F(t,x,u,p,X)&-F(t,y,u,q,-Y) \ge\\
 &-\max\left\{\frac{8}{27}\zeta^{2}\mu,8\mu\theta,2\mu\theta,\mu\eta +C_{1} ,\mu\kappa \right\}   \left[ \mu_{1}\left(|x-y|^{2} + \rho(p,q)^{2} \right)+ \mu_{2} \right. \\
&\hspace{2.2in}+ |p-q|+ |x-y|\left(\max(|p|,|q|)+1\right) \Big]\\
\end{aligned}
\end{gathered}
\end{equation*}
and setting $\omega_{R} = \max\left\{\frac{8}{27}\zeta^{2}\mu,8\mu\theta,2\theta,\eta +C_{1} ,\mu\kappa \right\} R$, this is a non-decreasing continuous function, maps $[0,\infty)\rightarrow [0,\infty)$ and $\omega_{R}(0) = 0$ as required. We have proven that condition (I7) is satisfied.

\bibliographystyle{plain}
\bibliography{biblio_V7}

\end{document}